\documentclass[11pt]{article}

\usepackage[leqno]{amsmath}
\usepackage{amsfonts,amssymb,amsthm,amscd,amsxtra}

\usepackage{times,avant,courier}
\usepackage{eucal} 
\usepackage{mathrsfs} 
\let\goth\mathfrak
\usepackage{stmaryrd}

\usepackage{xy}
\xyoption{all}
\SelectTips{cm}{}

\usepackage{xspace}

\swapnumbers
\newtheorem{theorem}[equation]{Theorem}
\newtheorem{proposition}[equation]{Proposition}
\newtheorem{lemma}[equation]{Lemma}
\newtheorem{corollary}[equation]{Corollary}
\newtheorem{lemma-def}[equation]{Lemma-Definition}

\theoremstyle{definition}
\newtheorem{definition}[equation]{Definition}
\newtheorem{example}[equation]{Example}

\theoremstyle{remark}
\newtheorem{remark}[equation]{Remark}

\numberwithin{equation}{subsection}

\newenvironment{enumeratep}{
  \begin{enumerate}}
{\end{enumerate}}

\newcommand{\abs}[1]{\left\lvert#1\right\rvert}
\newcommand{\norm}[1]{\lVert#1\rVert}
\newcommand{\tame}[2]{\bigl(#1,#2\bigr]}

\newcommand{\tamehh}[2]{\bigl(#1,#2\bigr]_\mathit{h.h.}}
\newcommand{\tameangle}[2]{\bigl\langle#1,#2\bigr]}
\newcommand{\tate}{2\pi\smash{\sqrt{-1}}}
\newcommand{\cprod}[1]{\wedge^{#1}} 
\newcommand{\unit}{\vartimes} 
\newcommand{\onehalf}{\frac{1}{2}}
\newcommand{\onefourth}{\frac{1}{4}}

\newcommand{\eqdef}{\overset{\text{\normalfont\tiny def}}{=}}
\newcommand{\gr}{\emph{gr}} 
\renewcommand{\phi}{\varphi}
\renewcommand{\epsilon}{\varepsilon}
\newcommand{\bei}{Be\u\i{}linson\xspace}
\newcommand{\cech}{\v{C}ech\xspace}

\newcommand{\loccit}{loc.\ cit.\xspace}
\newcommand{\lto}{\longrightarrow}
\newcommand{\lmto}{\longmapsto}
\newcommand{\iso}{\simeq}
\newcommand{\caniso}{\cong}
\newcommand{\isoto}{\overset{\sim}{\rightarrow}}
\newcommand{\lisoto}{\overset{\sim}{\longrightarrow}}
\newcommand{\coin}{\equiv}
\newcommand{\lqi}{\overset{\simeq}{\longrightarrow}}

\newcommand{\sheaf}[1]{\mathscr{#1}}

\newcommand{\coho}[1]{\mathrm{#1}}%
\renewcommand{\H}{\coho{H}}
\newcommand{\K}{\coho{K}} 
\newcommand{\hypercoho}[1]{\mathbf{#1}}%
\newcommand{\HHH}{\hypercoho{H}}
\newcommand{\complex}[2][\bullet]{{#2}^{#1}}
\newcommand{\cover}[1]{\mathfrak{#1}}

\DeclareMathOperator{\Aut}{Aut}
\newcommand{\SheafAut}{\operatorname{\underline{\mathrm{Aut}}}}
\DeclareMathOperator{\cone}{Cone}
\DeclareMathOperator{\dd}{d}\renewcommand{\d}{\dd\mspace{-2mu}}
\newcommand{\dc}{\dd^c\mspace{-2mu}}
\newcommand{\del}{\partial}
\newcommand{\delb}{\Bar\partial}
\newcommand{\equ}{\operatorname{\mathscr{E}\mspace{-3mu}\mathit{q}}}
\DeclareMathOperator{\Hom}{Hom}

\newcommand{\SheafHom}{\operatorname{\underline{\mathrm{Hom}}}}
\DeclareMathOperator{\id}{id}
\DeclareMathOperator{\Pic}{Pic}
\DeclareMathOperator{\Pichat}{\widehat{Pic}}
\DeclareMathOperator{\Spec}{Spec}

\newcommand{\field}[1]{\ensuremath{\mathbf{#1}}}
\newcommand{\ZZ}{\field{Z}}

\newcommand{\RR}{\field{R}}
\newcommand{\CC}{\field{C}}

\newcommand{\cat}[1]{\mathsf{#1}}
\newcommand{\twocat}[1]{\goth{#1}} 
\DeclareMathOperator{\ob}{Ob}
\DeclareMathOperator{\Ob}{Ob}

\newcommand{\Top}   {\cat{Top}}
\newcommand{\sch}[1]{\cat{Sch}/{#1}}

\newcommand{\an}{^\mathit{an}}
\newcommand{\Et}{_{\mathit{\Acute{E}t}}}
\newcommand{\et}{_{\mathit{\Acute{e}t}}}

\newcommand{\zar}{_{\mathit{zar}}}

\newcommand{\stack}[1]{\mathscr{#1}}
\newcommand{\stE}{\stack{E}}
\newcommand{\stF}{\stack{F}}
\newcommand{\stP}{\stack{P}}
\newcommand{\stQ}{\stack{Q}}

\newcommand{\grstack}[1]{\mathscr{#1}}
\newcommand{\grstA}{\grstack{A}}
\newcommand{\grstB}{\grstack{B}}
\newcommand{\grstC}{\grstack{C}}
\newcommand{\grstG}{\grstack{G}}
\newcommand{\grstH}{\grstack{H}}

\newcommand{\gerbe}[1]{\mathscr{#1}}
\newcommand{\gG}{\gerbe{G}}
\newcommand{\gH}{\gerbe{H}}
\newcommand{\gQ}{\gerbe{Q}}

\newcommand{\twogerbe}[1]{\goth{#1}}
\newcommand{\tgG}{\twogerbe{G}}
\newcommand{\tgH}{\twogerbe{H}}

\newcommand{\catAut}{\operatorname{\mathscr{A}\mspace{-3.5mu}\mathit{ut}}}
\newcommand{\catHom}{\operatorname{\mathscr{H}\mspace{-4mu}\mathit{om}}}

\newcommand{\conn}{\operatorname{\mathscr{C}\mspace{-3mu}\mathit{o}}}
\newcommand{\curv}{\operatorname{\mathscr{K}}}
\newcommand{\tors}{\operatorname{\textsc{Tors}}}
\newcommand{\bitors}{\operatorname{\textsc{Bitors}}}
\newcommand{\gerbes}{\operatorname{\textsc{Gerbes}}}
\newcommand{\twotors}{\tors}

\newcommand{\sha}[2][\bullet]{\mathscr{A}_{#2}^{#1}}
\newcommand{\she}[2][\bullet]{\mathscr{E}_{#2}^{#1}}
\newcommand{\sho}[1]{\mathscr{O}_{\mspace{-2mu}#1}}
\newcommand{\shomega}[2][\bullet]{\Omega_{#2}^{#1}}

\newcommand{\deligne}[4][\bullet]{%
  \def\tempa{}%
  \def\tempb{#2}%
  \ifx\tempa\tempb
  #3(#4)^{#1}_\mathcal{D}
  \else
  #3(#4)^{#1}_{\mathcal{D},#2}
  \fi
}
\newcommand{\deltilde}[4][\bullet]{%
  \smash[t]{\widetilde{#3_{#2}(#4)}}^{#1}_{\mathcal{D}}}
\newcommand{\delH}[4][\bullet]{{\H}^{#1}_\mathcal{D}(#2, #3(#4))}
\newcommand{\delz}[2][\bullet]{\deligne[#1]{}{\ZZ}{#2}} 

\newcommand{\delalg}[3][\bullet]{\goth{D}^{#1}(#2,#3)}

\newcommand{\brydhh}[3][\bullet]{C(#3)^{#1}_{#2}}

\newcommand{\dhh}[3][\bullet]{
  \def\tempa{}%
  \def\tempb{#2}%
  \ifx\tempa\tempb
  D_\mathit{h.h.}(#3)^{#1}
  \else
  D_\mathit{h.h.}(#3)^{#1}_{#2}
  \fi
}

\newcommand{\ndhh}[3][\bullet]{
  \def\tempa{}%
  \def\tempb{#2}%
  \ifx\tempa\tempb
  \goth{D}_\mathit{h.h.}(#3)^{#1}
  \else
  \goth{D}_\mathit{h.h.}(#3)^{#1}_{#2}
  \fi
}
\newcommand{\dhhH}[3][\bullet]{{\widehat\H}^{#1}_{\mathcal{D}}(#2,#3)}

\usepackage{hyperref}

\begin{document}

\title{2-Gerbes bound by complexes of \gr-stacks, and cohomology}
\author{Ettore Aldrovandi\\
  \small Department of Mathematics, Florida State University\\
  \small Tallahassee, FL 32306-4510, USA\\
  \small \texttt{aldrovandi@math.fsu.edu}}
\maketitle%
\begin{abstract}
  We define 2-gerbes bound by complexes of braided group-like
  stacks.  We prove a classification result in terms of
  hypercohomology groups with values in abelian crossed squares
  and cones of morphisms of complexes of length 3.  We give an
  application to the geometric construction of certain elements
  in Hermitian Deligne cohomology groups.
\end{abstract}
\tableofcontents%

\section*{Introduction}
\addcontentsline{toc}{section}{Introduction}
\label{sec:introduction}

The aim of the present work is to study in some detail gerbes
and, mostly, 2-gerbes bound by complexes of groups and braided
\gr-stacks, respectively, and the cohomology groups determined by
their equivalence classes.

\subsection*{Background and motivations}
\label{sec:backgr-motiv}

The idea of a gerbe bound by a complex is of course not new: it
dates back to Debremaeker (\cite{MR0480515}) in the form of a
gerbe $\gG$ on a site $\cat{S}$ bound by a crossed module $\delta
\colon A\to B$.  Milne (\cite{math.AG/0301304}) adopts the same
idea in the special case of an \emph{abelian} crossed module.  It
is observed in \loccit that the crossed module in fact reduces to
a homomorphism of sheaves of abelian groups, and the whole
structure simplifies to that of a gerbe $\gG$ bound by the sheaf
$A$ and equipped with a functor $\gG \to \tors (B)$ which is a
$\delta$-morphism, i.e.\ compatible with the homomorphism
$\delta$ (see below for the precise definition).

Our starting point is the observation that this structure
captures the differential geometric notion of ``connective
structure'' on an abelian gerbe, introduced by Brylinski and
McLaughlin\footnote{In \cite{bry:loop} the concept is ascribed to
  Deligne.} (\cite{brymcl:deg4:I,brymcl:deg4:II,bry:quillen}, see
also\cite{bry:loop} for a version in the context of smooth
manifolds). Briefly, by suitably generalizing the familiar
concept of connection on an invertible sheaf on an analytic or
algebraic manifold $X$, they defined a connective structure on an
abelian gerbe bound by $\sho{X}^\unit$ as a functor $x
\rightsquigarrow \conn (x)$ associating to each local object $x$
over an open $U$ a $\shomega[1]{U}$-torsor, subject to a certain
list of properties reviewed in sect.~\ref{sec:examples}.  It
turns out, and we show it explicitly in sect.~\ref{sec:examples},
that this is exactly the same thing as prescribing a structure of
gerbe bound by the complex
\begin{equation*}
  \sho{X}^\unit \xrightarrow{\d\log} \shomega[1]{X}\,.
\end{equation*}

More recently, we have similarly introduced the concept of
hermitian structure on an abelian gerbe bound by $\sho{X}^\unit$,
by modeling it on the corresponding familiar notion of invertible
sheaf equipped with a fiber hermitian metric (\cite{MR2142353}).
In simplified terms, this structure is also of the type
introduced above, namely we find that in this case it can be
conveniently encoded in the structure of gerbe bound by the complex
\begin{equation*}
  \sho{X}^\unit \xrightarrow{\log\abs{\cdot}} \she[0]{X}\,,
\end{equation*}
where the latter denotes the sheaf of smooth real functions on
$X$.

It is reasonable to expect that the list can be made longer with
other interesting examples.  However, we want point out that the
real interest of these construction lies in a different direction
(or directions).  On one hand, there is the obvious interest of
being able to generalize to the case of gerbes several structures
of differential geometric interest.  On the other, there is the
fact that typically equivalence classes (suitably defined) of
these structures turn out to be classified by interesting
cohomology theories, and as a feedback we can get a geometric
characterization of the elements of these groups.  For instance,
the cohomology groups relevant in the  above examples are the
Deligne cohomology group $\delH[3]{X}{\ZZ}{2}$, and the hermitian
Deligne cohomology group $\dhhH[3]{X}{1}$.

In fact, Brylinski and McLaughlin have shown that their
constructions provide the adequate context for notable extensions
of the tame symbol map in algebraic $\K$-theory, where gerbes are
useful in order to obtain a geometric picture for some regulator
maps to Deligne-\bei cohomology (cf.\ \cite{MR1317118}).  More
importantly, they extend their framework in two directions:
(1)~they consider the case of 2-gerbes as well, and (2)~they
define appropriate notions of curvature both for gerbes and
2-gerbes bound by $\sho{X}^\unit$.  Passing from gerbes to
2-gerbes corresponds to an increase in the degree of the involved
cohomology groups, whereas introducing more levels of
differential geometric structures corresponds to cohomology
groups of higher \emph{weights.}  The geometric and the
cohomological aspects are tied together very neatly in the
following sense: the Deligne cohomology groups
$\delH[p]{X}{A}{k}$, where $A$ is a subring of $\RR$, can be
regarded as somewhat pathological in the range $p > 2k$, where
they cannot receive regulator maps from, say, absolute
cohomology.\footnote{The absolute cohomology groups in that range
  are zero.} It is reassuring that the gerbes and 2-gerbes
corresponding to the tame symbol maps and various related cup
products turn out to naturally have a connective structure (and
even curvatures), so that their classifying Deligne cohomology
groups lie in the ``safe'' range $p\leq 2k$.\footnote{There is of
course an interest in knowing that, say, $\delH[3]{X}{\ZZ}{1}$
classifies abelian gerbes bound by $\sho{X}^\unit$, however the
nice connection with regulators, etc.\ is lost.}

A similar story was developed by the author in the case of
hermitian Deligne cohomology (\cite{MR2142353}), motivated by the
existence of certain natural hermitian structures on tame
symbols.  As mentioned before, the cohomological counterpart is
given by hermitian Deligne cohomology, and there is a parallel
for 2-gerbes as well.  Namely, we have put forward a definition
of hermitian structure for 2-gerbes (to be reviewed and revised
below) bound by $\sho{X}^\unit$ and found that the corresponding
equivalence classes are in 1--1 correspondence with the elements
of the group $\dhhH[4]{X}{1}$. In particular, the gerbes and
2-gerbes corresponding to the tame symbols studied by Brylinski
and McLaughlin were found to naturally support a hermitian
structure as well.  Moreover, it was found that these structures,
namely the analytic (or algebraic) connective structure of
Brylinski and McLaughlin and the hermitian structure we
introduced are compatible in the following sense: One of the
byproducts of our work is that there is a natural notion of
connective structure canonically associated with the hermitian
structure. It was found that this new connective structure agrees
with the one of Brylinski-McLaughlin once they are mapped into an
appropriate complex of smooth forms. (Part of this theory will be
recalled and further clarified in the last part of the present
paper.)

Not quite satisfying, as the reader will have no doubt noticed,
is the fact that weights and degrees are precisely in what seems
to be the bad range.  However, a more interesting group
$\dhhH[4]{X}{2}$ does appear in the following way: in
\cite{MR2142353} we introduced a complex, denoted
$\complex{\Gamma (2)}$ (defined in
section~\ref{sec:trunc-herm-deligne}), and we (informally) argued
that the hypercohomology group $\HHH^4(X,\complex{\Gamma (2)})$
classifies 2-gerbes equipped with both a connective structure \`a
la Brylinski-McLaughlin and a hermitian structure in our sense
which are compatible as explained to above.  In \loccit we found
there is a surjection $\dhhH[4]{X}{2} \to
\HHH^4(X,\complex{\Gamma (2)})$, so classes of 2-gerbes can
indeed be lifted to a more desirable group, but a truly geometric
characterization was not provided.  Let us remark that the
interest of the group $\dhhH[4]{X}{2}$ lies in the fact that it
is the receiving target of the cup product map
\begin{equation*}
  \Pichat X \otimes \Pichat X \lto \dhhH[4]{X}{2}\,,
\end{equation*}
where $\Pichat X\iso \dhhH[2]{X}{1}$ is the group of isomorphism
classes of metrized invertible sheaves. When $X$ is a complete
curve, this map gives a cohomological interpretation of Deligne's
determinant of cohomology construction (\cite{MR89b:32038}),
which has been analyzed in various guises in
\cite{bry:quillen,math.CV/0211055}, and \cite{MR2145708} in the
singular case.

The desire to remedy the above shortcoming and enhance the
results of \cite{MR2142353}, as well as the desire to cast the
results in the form expounded at the beginning of this
introduction---suitably extended to include 2-gerbes---constitute
our motivation for the present work. The framework we have found,
that of 2-gerbes bound by a complex of braided \gr-stacks, is
quite more general than what would be minimally required for just
solving the mentioned problems, and lends itself to possible
generalizations to the non-abelian case, which we plan to address
in part in a subsequent publication. We now proceed to describe
the present results in the remaining part of this introduction.

\subsection*{Statement of the results}
\label{sec:statement-results}

For the purpose of this introduction let us informally assume
that $X$ is a smooth base scheme, or an analytic manifold, and
that $\cat{C}/X$ is an appropriate category of spaces ``over''
$X$ with a Grothendieck topology, making it into a site.

To keep track of cohomology degrees, recall that Deligne
cohomology and its variants have a built-in degree index shift.
The convention we use in this introduction and the rest of the
paper is to revert to standard cohomology degrees whenever we are
not specifically dealing with one of these specific cohomology
theories.

Our first result is a straightforward generalization of the
concept of abelian gerbe bound by a homomorphism of sheaves of
abelian groups to the case where we have a complex of abelian
groups of the form:
\begin{equation*}
  A \overset{\delta}{\lto} B \overset{\sigma}{\lto} C\,.
\end{equation*}
We find that an abelian gerbe $\gG$ bound by the above complex is
conveniently defined as an $A$-gerbe $\gG$ equipped with a
functor
\begin{equation*}
  \gG \lto \tors (B,C)\,,
\end{equation*}
where the right hand side denotes the gerbe of $B$-torsors with a
section of the associated $C$-torsor obtained by extension of the
structure group from $B$ to $C$.  We then obtain through a simple
\cech cohomology argument that equivalence classes of such gerbes
are classified by the hypercohomology group
\begin{equation*}
  \HHH^2(X,A\to B\to C)\,.
\end{equation*}
We show at the end of section~\ref{sec:complexes-length-3} that
this is the appropriate general cadre for the notion of
curvature: indeed we prove that Brylinski and McLaughlin's
original definition of a gerbe with connective structure and
``curving'' can be cast as a gerbe bound by a complex of length
3, for an appropriate choice of the groups involved.

The extension of the idea of gerbe bound by a complex to the case
of 2-gerbes is more involved, but quite interesting.

We want to consider abelian 2-gerbes, where of course the word
``abelian'' must be properly qualified. We adopt the point of
view of \loccit of calling ``abelian'' a 2-gerbe bound by a
braided \gr-stack in the following sense: It is known that the
fibered category of automorphisms of an object $x$ over $U\to X$
in a 2-gerbe is a \gr-stack.  Let $\grstA$ be a \gr-stack over
$X$.  A 2-$\grstA$-gerbe $\tgG$ is a 2-gerbe with the property
that each local automorphism \gr-stack is equivalent to (the
restriction of) $\grstA$. As we know from \cite{MR95m:18006}, if
this equivalence is natural in $x$, then $\grstA$ will be forced
to be braided, i.e.\ its group law has a non-strict commutativity
property. A special case is when $\grstA=\tors (A)$, that is, the
\gr-stack is the stack of torsors (in fact, a gerbe) over an
abelian group $A$. Then we speak of a 2-gerbe bound by $A$, or
2-$A$-gerbe.

Note that it follows from \cite{MR93k:18019,MR95m:18006} that for
an abelian 2-$\grstA$-gerbe $\tgG$ the stack of morphisms
$\catAut_U(x,y)$ of two objects over $U \to X$ has the structure
of $\grstA\vert_U$-torsor, and that $\tgG$ determines a
1-cocycle, hence a cohomology set, with values in $\tors
(\grstA)$.  Note that for any \gr-stack $\grstA$ this is a
neutral 2-gerbe, see \cite{MR92m:18019}.  By suitably decomposing
the torsors comprising this cocycle, we obtain a degree 2
cohomology set with values in $\grstA$ itself. This leads to the
familiar degree 3 cohomology group with values in $A$ in the case
$\grstA = \tors (A)$. We will find generalizations by studying
the analogous constructions for complexes of \gr-stacks, defined
below.

Thus, given an additive functor $\lambda \colon \grstA \to
\grstB$ of braided \gr-stacks we define a 2-gerbe bound by this
``complex'' as a pair $(\tgG,J)$, where $\tgG$ is a
2-$\grstA$-gerbe and $J$ is a cartesian 2-functor
\begin{equation*}
  J\colon \tgG \lto \tors (\grstB)
\end{equation*}
which is a $\lambda$-morphism, see section~\ref{sec:2-a-b} for
the precise definition.  Once the notion of morphism and then of
equivalence of such pairs are defined, we find that equivalence
classes are in 1--1 correspondence with the elements of a
cohomology set which we could provisionally write as:
\begin{equation*}
  \H^1(\tors (\grstA)\to \tors (\grstB))\,.
\end{equation*}
Once again, by suitably decomposing the torsors comprising the
1-cocycle with values in the complex
\begin{equation*}
  \lambda_\ast \colon \tors
  (\grstA)\to \tors (\grstB)
\end{equation*}
determined by $\tgG$, we obtain a
degree 2 cohomology set with values in the complex $\lambda
\colon \grstA \to \grstB$ itself.

In order to properly handle the hermitian Deligne cohomology
group we are ultimately interested in, we can further generalize
this notion to that of a 2-gerbe bound by a \emph{complex} of
\gr-stacks, that is a diagram of additive functors:
\begin{equation}
  \label{eq:18}
  \grstA \overset{\lambda}{\lto}
  \grstB \overset{\mu}{\lto} \grstC \tag{+}
\end{equation}
where the composition $\mu\circ\lambda$ is required to be
isomorphic to the null functor sending $\grstA$ to the unit
object of $\grstC$. Thus a 2-gerbe $\tgG$ is bound by the above
complex if there is a cartesian 2-functor
\begin{equation*}
  \Tilde J\colon \tgG \lto \tors (\grstB,\grstC)\,,
\end{equation*}
where the right hand side denotes the 2-gerbe of $\grstB$-torsors
which become equivalent to the trivial $\grstC$-torsor.  Then we
show that equivalence classes of such pairs $(\tgG,\Tilde J)$ are
classified by a cohomology set:
\begin{equation*}
  \H^1(\tors (\grstA)\to \tors (\grstB) \to \tors (\grstC))\,,
\end{equation*}
from which we can obtain a degree two cohomology set with
coefficients in the \gr-stack complex above. This is done in
sections~\ref{sec:2-gerbes-bound} and~\ref{sec:2-gerbes-bound-2},
where the relevant theorems are stated and proven in full.

Along the way we get interesting byproducts shedding a new light
on the notion of gerbe bound by a complex. In
section~\ref{sec:class-results} we prove that for a strictly
abelian (and not just braided) \gr-stack $\grstB$, that is, one
that arises from a homomorphism of sheaves of abelian groups, we
have the equivalence
\begin{equation*}
  \gerbes (B,H) \lisoto \tors (\grstB)
\end{equation*}
where $\grstB = \tors (B,H)$.  Then later in
section~\ref{sec:grstb-grstc-torsors}, we observe that $\tors
(\grstB, \grstC)$ introduced above is equivalent, when
$\grstC=\tors (C,K)$ with the 2-gerbe of gerbes bound by $B\to H$
which become neutral as gerbes bound by $C\to K$.

These partial results are part of a general process whereby we
make contact with ordinary hypercohomology by assuming all all
the involved \gr-stacks are strictly abelian.  Concretely, if
$\grstA=\tors (A,G)$, $\grstB = \tors (B,H)$, and $\grstC=\tors
(C,K)$ the complex of \gr-stacks we have been considering reduces
to the commutative diagram of (sheaves of) abelian groups:
\begin{equation}
  \label{eq:19}
  \vcenter{\xymatrix{%
    A \ar[d]_\delta \ar[r]^f & B \ar[d]^\sigma \ar[r]^g & C
    \ar[d]^\tau \\ 
    G \ar[r]_u & H \ar[r]_v & K}}\tag{*}
\end{equation}
The theorem we obtain in section~\ref{sec:classification-iii} is
that equivalence classes of 2-gerbes bound by the
complex~\eqref{eq:18} are classified by the standard
hypercohomology group
\begin{equation*}
  \HHH^3(X,(\text{cone of~\eqref{eq:19}})[-1])\,.
\end{equation*}
As we will see in section~\ref{sec:geom-appl}, this is exactly
the kind of cohomology group we need in order to give a geometric
construction of the elements of the hermitian Deligne cohomology
group $\dhhH[4]{X}{2}$.  In particular, in
section~\ref{sec:geom-interpr-some}, we give a reasonably
detailed construction of a 2-gerbe, denoted
$\tamehh{\sheaf{L}}{\sheaf{M}}\sphat$, whose class in
$\dhhH[4]{X}{2}$ is the cup product $[\sheaf{L},\rho]\cup
[\sheaf{M},\sigma]$ of $[\sheaf{L},\rho], [\sheaf{M},\sigma]\in
\Pichat X$.

In section~\ref{sec:2-gerbes-bound}, especially in
sections~\ref{sec:class-results} and~\ref{sec:2-gerbes-bound-1}
we prove intermediate results for the case where there is no
$\grstC$, so the diagram~\eqref{eq:19} above reduces to the left
square.

In all cases, when moving from cohomology sets with values in
complexes of gerbes of torsors to (hyper)cohomology groups with
values in cone of complexes, we compute explicit cocycles with
respect to hypercovers, rather than ordinary covers.  We find
that even in the case of groups the cocycles so obtained present
additional interesting terms.

\subsection*{Organization and contents of the paper}

Overall we have adopted a mix of bottom-up and top-down
approaches.  We have refrained from starting from the most
general statement and then working our way down.  Instead we have
adopted a sequence of successive generalizations.

Our treatment of cohomology deserves some explanations. At the
beginning, where several proofs are standard, we have adopted a
\cech point of view.  In the latter part of the paper, where we
deal with torsors over \gr-stacks, we have found worthwhile
\emph{not} to assume that decompositions with respect to \cech
covers are sufficient.  So we have actually computed cocycles
using hypercovers, adopting the same point of view and formalism
of~\cite{MR95m:18006}. Since we have dealt with hypercovers in a
rather direct way, formulas acquire a substantial decoration of
indices, which can be quite daunting.  The usual advice is to
ignore the hypercover indices on first parsing and reduce
everything to the \cech formalism and replace (hyper)cohomology
with its \cech counterpart. 

A note about sites: When dealing with categorical matters, it
comes at no additional cost to formulate everything, including
cohomology sets, for sites.  Thus usually we will assume that
gerbes and 2-gerbes are fibered over a site $\cat{S}$.  This site
will in fact be a category of objects over an object $X$, so that
we will often use the notation $\cat{C}/X$, assuming the category
$\cat{C}$ has been equipped with an appropriate Grothendieck
topology. By thinking of $X$ as the terminal object in
$\cat{C}/X$, we can conveniently denote cohomology sets as
$\H^\bullet(X,-)$ or $\HHH^\bullet(X,-)$, depending on whether we
wish to emphasize the ``hyper'' aspect.

This paper is organized as follows. In
section~\ref{sec:background-notions} we recall a few background
notions, collect some notation, and we provide a quick overview
of various Deligne-type cohomology theories needed in the rest of
the paper.

We introduce the concept of gerbe bound by a length 2 complex in
section~\ref{sec:review-gerbes-with}, where we also review the
pivotal example of connective structure in some detail. We then
proceed in section~\ref{sec:complexes-length-3} to define and
classify gerbes bound by a length 3
complex. Section~\ref{sec:2-gerbes} is dedicated to a quick
review of 2-gerbes. Unfortunately we cannot make this paper
completely self-contained without writing another book on
2-gerbes, therefore referring to the literature,
especially~\cite{MR95m:18006}, remains indispensable.
Sections~\ref{sec:2-gerbes-bound} and~\ref{sec:2-gerbes-bound-2}
then contain our main results, where we classify 2-gerbes bound
by complexes of \gr-stacks. Finally, in
section~\ref{sec:geom-appl}, we return to the realm complex
algebraic manifolds, and give some applications to hermitian
Deligne cohomology.

\subsection*{Acknowledgments}
I warmly thank Larry Breen for an email exchange clarifying a few points.

Partial support by FSU CRC under a COFRS grant no.~015290. is
gratefully acknowledged.

\section{Background notions}
\label{sec:background-notions}

\subsection{Assumptions and notations}
\label{sec:assumpt-notat}

In the following, $X$ will be a smooth scheme or a complex
analytic manifold. In the algebraic case, some results can be
stated for $X$ smooth over a base scheme $S$. Actually, in most
of the applications we will be concerned with the case when $X$
is an algebraic manifold,\footnote{By algebraic manifold we mean
  a smooth, separated scheme of finite type over \CC} hence
$S=\Spec \CC$. In this case the complex analytic manifold above
will be $X\an$, the set of complex points of $X$ with the
analytic topology, but usually we will not explicitly mark this
in the notation.

Gerbes ``over $X$'' are stacks in groupoids and, similarly,
2-gerbes are 2-cate\-gories fibered in (lax) 2-groupoids satisfying
certain conditions to be explained below, over an appropriate
site of ``spaces'' over $X$.  As explained at the end of the
introduction, whenever dealing with general categorical matters,
the specific choice of this site will be somewhat immaterial.  In
order to fix ideas, and to revert in the end to specific
cohomology theories, we will assume that we are given an
appropriate category with fiber products $\cat{C}/X$ of spaces
over $X$ equipped with a Grothendieck topology.  The main
requirement will be that the various sheaves such as $\sho{X}$,
$\shomega{X}$, etc.\ as defined with respect to $\cat{C}/X$
restrict to their usual counterparts under $U\to X$, whenever $U$
is open in the ordinary---for the Zariski or Analytic
topology---sense. More specifically, following ref.\ 
\cite{bry:loop}, if $X$ is a scheme we may as well consider the
small \'etale site $X\et$, namely $\cat{C}/X=\cat{Et}/X$, where
we denote by $\cat{Et}$ the class of \'etale maps over $X$, and
covers are jointly surjective families of \'etale maps. It is
useful to allow the same type of construction when $S=\Spec \CC$,
and we want to consider $X\an$. Namely we obtain a corresponding
``analytic'' site by mapping $U\to X$ from $X\et$ to $U\an \to
X\an$.  According to ref.\ \cite{MR559531}, this determines the
same topology as the standard analytic one.  In the latter case,
that is if $X$ is a complex manifold, $\cat{C}/X$ will be the
small $\Top$ site. Similarly, when $X$ is a scheme to be
considered with its ordinary topology, we set $\cat{C}/X=X\zar$,
the small Zariski site of $X$ whose covers are injective maps
$V\to U$ with $U$ open in $X$.  Note that in general we will not
be considering the corresponding ``big'' sites. However, the
general categorical constructions which form the main body of
this paper are going to work in that context too.~\footnote{To be
  more specific one could consider sites such as $X\Et$, the big
  \'etale site of $X$, if $X$ is a scheme, namely
  $\cat{C}/X=\sch{X}$ equipped with the \'etale topology defined
  by the class $\cat{Et}$ of \'etale maps over $X$;
  correspondingly, $\cat{C}/X=\cat{Cmplx}/X$, with the topology
  given by standard open covers, or by analytification of \'etale
  covers as described above.}

In general we will refer to the topology on $\cat{C}/X$
simply as a topology on $X$, and accordingly we will simply speak
of ``open'' sets for members $V\to U$ of a cover of $U\to X$. As
it is well known, fibered products take the place of
intersections, and we will use the standard notation of denoting
the various multiple ``intersections'' (i.e.\ fibered products)
relative to a covering $\lbrace U_i\to U\rbrace_{i\in I}$ as:
$U_{ij}=U_i\times_U U_j$, $U_{ijk}=U_i\times_U U_j\times_U U_k$,
etc.  Also in the relative case of $X$ over a base $S$,
$\cat{C}/X$ will be obtained by restriction from
$\cat{C}/S$. However, our notation will not always
explicitly reflect this.

\subsubsection{Often used  notations.}
\label{sec:gener-assumpt-notat}

For a subring $A$ of $\RR$ and an integer $p$, $A(p) =
(\tate)^p\,A$ is the Tate twist of $A$. We identify $\CC/\ZZ(p)
\iso \CC^\unit$ via the exponential map $z \mapsto \exp
(z/(\tate)^{p-1})$, and $\CC\iso \RR(p) \oplus \RR (p-1)$, so
$\CC / \RR(p) \iso \RR (p-1)$.  The projection $\pi_p\colon \CC
\to \RR(p)$ is given by $\pi\sb{p} (z) = \onehalf ( z + (-1)\sp p
\Bar z)$, for $z\in \CC$, and similarly for any other complex
quantity.

If $E$ is a set (or group, ring, module...), then $E_X$ denotes
the corresponding constant sheaf of sets (or groups, rings,
modules...).

If $X$ is a scheme or complex manifold, $\shomega{X}$ denotes the
corresponding (algebraic or analytic) de~Rham complex.  We set
$\sho{X} \coin \shomega[0]{X}$ as usual.  $\she{X}$ denotes the
de~Rham complex of sheaves of $\RR$-valued smooth forms on the
underlying smooth manifold. Furthermore, $\sha{X} =
\she{X}\otimes_\RR \CC$, and is $\she{X} (p)$ the twist $\she{X}
\otimes_\RR \RR(p)$.  Also, $\sha[p,q]{X}$ will denote the sheaf
of smooth $(p,q)$-forms, and $\sha[n]{X} = \bigoplus_{p+q=n}
\sha[p,q]{X}$, where the differential decomposes in the standard
fashion, $\d=\del + \delb$, according to types. We also introduce
the imaginary operator $\dc = \del -\delb$~\footnote{We omit the
  customary factor $1/(4\pi \sqrt{-1})$} and we have the rules
\begin{displaymath}
  \d\pi_p(\omega) = \pi_p(\d\omega)\,,\quad
  \dc\pi_p(\omega) = \pi_{p+1}(\dc\omega)
\end{displaymath}
for any complex form $\omega$.  Note that we have $2\del\delb =
\dc\d$.

The standard Hodge filtrations on $\shomega{X}$ and $\sha{X}$ are
as follows: $F^p\shomega{X}\coin \sigma^p \shomega{X}$ is the
sharp truncation in degree $p$:
\begin{equation*}
  0\lto \dotsm \lto 0 \lto
  \shomega[p]{X}\lto \dotsm \lto \shomega[\dim X]{X}\,,
\end{equation*}
whereas $F^p\!\sha{X}$ is the total complex of:
\begin{math}
  \bigoplus_{r\geq p} \sha[r,\bullet -r]{X}\,.  
\end{math}

\subsection{Various Deligne complexes and cohomologies}
\label{sec:vari-deligne-compl}

Standard references on Deligne cohomology are:
\cite{MR86h:11103,esn-vie:del}.

For a subring $A\subset \RR$ and an integer $p$, the Deligne
cohomology groups of weight $p$ of $X$ with values in $A$ are the
hypercohomology groups:
\begin{equation}
  \label{eq:1}
  \delH{X}{A}{p} \eqdef \HHH^{\,\bullet} (X,\deligne{X}{A}{p})\,,
\end{equation}
where $\deligne{X}{A}{p}$ is the complex
\begin{align}
  \label{eq:2}
  \deligne{X}{A}{p} & =
  \cone \bigl( A(p)_X\oplus F^p\shomega{X} \lto \shomega{X}
  \bigr)[-1] \\
  \label{eq:3}
  & \lqi \bigl( A(j)_X
  \overset{\imath}{\lto} \sho{X}
  \overset{\d}{\lto} \shomega[1]{X}
  \overset{\d}{\lto} \dotsm
  \overset{\d}{\lto} \shomega[{p-1}]{X}\bigr) \,.
\end{align}
where the map in the cone is the difference of the two inclusions
and $\lqi$ denotes a quasi-isomorphism.  The complex
in~\eqref{eq:3} is the one we will normally use in what follows.

When $A=\RR$, Deligne cohomology groups can be computed using
other complexes quasi-isomorphic to~\eqref{eq:2}
or~\eqref{eq:3}, in particular:
\begin{equation}
  \label{eq:4}
  \deltilde{}{\RR}{p} =
  \cone \big(F^p\!\sha{X} \rightarrow \she{X}(p-1)\big)[-1]\,.
\end{equation}
(See the references quoted above for a proof.)

The \emph{Hermitian} variant of Deligne cohomology is obtained by
considering the hypercohomology groups
\begin{equation}
  \label{eq:5}
  \dhhH{X}{p} \eqdef \HHH^{\,\bullet} (X,\brydhh{X}{p})
\end{equation}
of the complex
\begin{equation}
  \label{eq:6}
  \brydhh{X}{p} = \cone \bigl(
  \ZZ(p)_X \bigoplus (F^p\!\sha{X}\cap \sigma^{2p}\she{X}(p))
  \lto \she{X}(p)
  \bigr)[-1]\,,
\end{equation}
introduced by Brylinski in
\cite{bry:quillen}. We proved in~\cite{math.CV/0211055} that it
is quasi-isomorphic to the complex:
\begin{equation}
  \label{eq:7}
  \dhh{X}{p} = 
  \cone \bigl(
  \delz{p}\oplus (F^p\!\sha{X}\cap \sigma^{2p}\she{X}(p))
  \lto \deltilde{}{\RR}{p}
  \bigr)[-1]\,.
\end{equation}
The interest of~\eqref{eq:7} lies in the fact that the second
hypercohomology group of $\dhh{X}{1}$ provides a characterization
of the \emph{canonical connection} associated to a hermitian line
bundle (\cite{math.CV/0211055,MR2142353}). We will also need a
leaner version of the complex~\eqref{eq:7} introduced
in~\cite{MR2145708}, namely:
\begin{equation}
  \label{eq:8}
  \ndhh{X}{p} =
  \cone \bigl(
  \delz{p} \overset{\rho_p}{\lto} \sigma^{<2p}\delalg{\sha{X}}{p}
  \bigr)[-1] \,.
\end{equation}
Here $\delalg{\sha{X}}{p}$ is the \emph{Deligne Algebra} over the
complex $\sha{X}$, discussed in full in
\cite{MR99d:14015,math.AG/0404122,gonch:JAMS}, and $\sigma^{<2p}$
denotes its sharp truncation in degrees above $2p$, so that:
\begin{equation}
  \label{eq:9}
  \sigma^{<2p}\delalg[n]{\sha{X}}{p} =
  \begin{cases}
    0 & n=0\,,\\
    \she[n-1]{X}(p-1)\bigcap
    \bigoplus_{\substack{p'+q'=n-1\\ p'< p, q'< p}}
    \sha[p',q']{X} & n\leq 2p-1\,.
  \end{cases}
\end{equation}
The differential is $-\pi\circ\d$, where $\pi$ is the projection
that simply chops off the degrees falling outside the scope
of~\eqref{eq:9}. Using~\eqref{eq:3}, the map $\rho_p$ is:
\begin{equation*}
    \rho_p^n = 
  \begin{cases}
    0 & n=0\,,\\
    (-1)^n\pi_{p-1} & 1\leq n\leq p\,.
  \end{cases}
\end{equation*}
\begin{example}
  In the following we will be concerned almost exclusively with
  the complexes of weight $p=1$ and $p=2$. Explicitly, we have:
  \begin{equation}
    \label{eq:10}
    \ndhh{X}{1} = \bigl(
    \ZZ(1)_X \overset{\imath}{\lto} \sho{X}
    \xrightarrow{\pi_0} \she[0]{X}\bigr)  \,,
  \end{equation}
  whereas the complex $\ndhh{X}{2}$ is the cone (shifted by 1) of
  the map:
  \begin{equation}
    \label{eq:11}
    \vcenter{\xymatrix@=2.3pc{%
        \ZZ(2)_X \ar[r]^\imath & \sho{X} \ar[r]^\d
        \ar[d]^{-\pi_1} & \shomega[1]{X} \ar[d]^{\pi_1} \\
        & \she[0]{X}(1) \ar[r]^{-\d} & \she[1]{X}(1)
        \ar[r]^(.4){-\pi\circ \d}  &
        \she[2]{X}(1)\cap \sha[1,1]{X}}}
  \end{equation}
\end{example}
\begin{remark}
  \label{rem:1}
  Using the complex~\eqref{eq:10}, one shows that
  \begin{equation*}
    \dhhH[2]{X}{1} \iso \Pichat X\,,
  \end{equation*}
  the group of isomorphism classes of line bundles with hermitian
  metric. This follows from an easy \cech argument, as in
  \cite{esn:char}. Thus the same type of argument, using the
  complex $\dhh{X}{1}$, implies the uniqueness of the canonical
  connection, see \cite{MR2142353}.
\end{remark}
We conclude this review section by observing that all complexes
introduced so far possess a product structure (or several
mutually homotopic such structures), additive with respect to the
weights, so that we have graded commutative cup products
\begin{equation*}
  \delH[k]{X}{A}{p}\otimes \delH[l]{X}{A}{q}
  \overset{\cup}{\lto} 
  \delH[k+l]{X}{A}{p+q}
\end{equation*}
and
\begin{equation*}
    \dhhH[k]{X}{p}\otimes \dhhH[l]{X}{q}
  \overset{\cup}{\lto} 
  \dhhH[k+l]{X}{p+q}\,.
\end{equation*}
The reader should refer to the literature cited in this section
for more details and explicit formulas for the products.

\section{Gerbes with abelian band}
\label{sec:review-gerbes-with}

In the following we recall a few definitions about gerbes. The
canonical reference is \cite{MR49:8992}, whereas a detailed
exposition adapted to spaces is \cite{MR95m:18006}. We will need
the abelian part of the whole theory, for which a readable
account is to be found in \cite{bry:loop}.

Let $\cat{C}$ be a category with finite fibered products,
equipped with a Grothendieck topology.  A \emph{gerbe} $\gG$ over
$\cat{C}$ is a stack in groupoids $p\colon \gG\to \cat{C}$ such
that:
\begin{enumerate}
\item $\gG$ is \emph{locally non-empty,} namely there exists a
  cover $U\to X$ such that $\ob (\gG_U)$ is non-empty;
\item $\gG$ is \emph{locally connected,} that is, for each pair
  of objects of $\gG$, there is a cover $\phi\colon V\to U$
  such that their inverse images are isomorphic. In other words
  if $x,y\in \ob\gG_U$, then $\Hom_U(\phi^\ast x,\phi^\ast
  y)$ is non-empty.
\end{enumerate}
For an object $x\in \ob \gG_U$, the sheaf $\SheafAut(x)$ is a
sheaf of groups on $\cat{C}/U$. (Recall that over $\phi\colon
V\to U$, we have $\SheafAut(x)(V) = \Aut_V(\phi^\ast x)$.) Let
now $A$ be a sheaf of groups on $\cat{C}$: We say that $\gG$ is
an \emph{$A$-gerbe} if for each object $x$ with $\phi (x)=U$ as
above there is a \emph{natural} isomorphism
\begin{equation*}
  a_x\colon \SheafAut (x) \lisoto A\vert_U\,.
\end{equation*}
The naturality in $x$ will force the group $A$ to be abelian, and
in the following we will restrict our attention to this case. The
sheaf $A$ will be referred to as the \emph{band} of the gerbe
$\gerbe{G}$. We also say that $\gG$ is \emph{bound} by $A$. (In
the general---non-abelian---case, the band $L(A)$ will have a
more complicated definition, as the various sheaves $A\vert_{U}$
are glued along $U\times_XU$ only up to inner automorphisms. In
the abelian case this is immaterial and we can abuse the language
and call $A$ the band of $\gerbe{G}$.)

A morphism $\lambda : \gG \to \gH$ is a cartesian functor between
the underlying fibered categories, and it is an equivalence if it
is an equivalence of categories. Moreover, if $\gG$ is an
$A$-gerbe, and $\gH$ is a $B$-gerbe, with a group homomorphism
$f\colon A\to B$, then the morphism $\lambda$ will have to
satisfy the obvious commutative diagrams. Such a morphism is
called an $f$-morphism.

An $f$-morphism for which $f$ is an isomorphism is automatically an
equivalence. So is, in particular, a morphism between two $A$-gerbes
$\gG$ and $\gG'$. So if $A$ is abelian, it follows from
\cite{MR49:8992} that $A$ classes of equivalences of $A$-gerbes are
classified by $\H^2(X,A)$, the standard second cohomology group of $X$
in the derived functor sense. See also, e.g. \cite{bry:loop}, for a
proof in the \cech setting.

\subsection{Gerbes bound by a complex}
\label{sec:gerbes-bound-complex}

We are going to use the notion of gerbe bound by a length two
complex $A \to B$ of sheaves of abelian groups over
$\cat{C}/X$, as in \cite{math.AG/0301304}.  Let us review
the formal definition:
\begin{definition}\label{def:1}
  Let $A$ and $B$ be two sheaves of abelian groups on
  $\cat{C}/X$, and $\delta\in \Hom(A,B)$, so that
  $A\overset{\delta}{\lto} B$ is a complex of length two. A gerbe
  $\gG$ bound by $A\to B$ is an $A$-gerbe over $\cat{C}/X$
  equipped with a $\delta$-morphism of gerbes
  \begin{equation*}
    \mu\colon \gG\to \tors(B)\,.
  \end{equation*}
\end{definition}
\noindent
(Notice that $\tors(B)$ is a $B$-gerbe, so the notion of
$\delta$-morphism makes sense.)

More generally, one would have
the notion of a gerbe $\gG$ bound by a sheaf of \emph{crossed
  modules,} as per Debremaeker's original definition in ref.\ 
\cite{MR0480515}. If $(A,B,\delta)$ be a crossed module, where
$\delta\colon A\to B$ is a group homomorphism, compatible with
the action of $B$ over $A$, a gerbe bound by it is an $A$-gerbe
$\gG$ with a $\delta$-morphism $\lambda$ above, together with
other data relative to the stacks of automorphisms of local
objects, see ref.\ \cite{MR0480515}.  When both $A$ and $B$ are
abelian, the crossed module becomes simply a complex, and
everything reduces to the data in the previous definition.

As usual, a morphism of complexes $(f,g)\colon (A,B,\delta) \to
(A',B',\delta')$ is a commutative diagram of group homomorphisms:
\begin{equation*}
    \xymatrix{%
      A \ar[r]^\delta \ar[d]_f & B\ar[d]^g \\
      A' \ar[r]^{\delta'} & B'
    }
\end{equation*}
If $\gG$ and $\gG'$ are bound by $(A,B)$ and $(A',B')$,
respectively, then we have a corresponding notion of
$(f,g)$-morphism as follows:
\begin{definition}\label{def:2}
  An $(f,g)$-morphism from $\gG$ to $\gG'$ consists of:
  \begin{enumerate}
  \item an $f$-morphism $\lambda\colon \gG\to \gG'$;
  \item a natural isomorphism of functors
    \begin{equation*}
      \alpha \colon g_\ast\circ \mu
      \Longrightarrow \mu'\circ \lambda
    \end{equation*}
    from $\gG$ to $\tors(B')$.
  \end{enumerate}
\end{definition}
In the definition $g_\ast$ is the $g$-morphism $\tors(A)\to
\tors(B)$ induced by $g$ in the obvious way.

For completeness, let us also mention that we also have the
notion of morphism of morphisms,
see~\cite{math.AG/0301304}. Namely, let $(\lambda_1,\alpha_1)$
and $(\lambda_2,\alpha_2)$ be two morphisms $(\gG,\mu) \to
(\gG',\mu')$. A morphism $m\colon (\lambda_1,\alpha_1) \to
(\lambda_2,\alpha_2)$ is a natural transformation $m\colon
\lambda_1 \Rightarrow \lambda_2$ such that the following is
verified:
\begin{equation*}
  (\mu \ast m)\circ \alpha_1=\alpha_2\,.
\end{equation*}
With these notions the gerbes bound by a complex of length 2 form
a 2-category.  In particular, when $A'=A$ and $B=B'$ we denote
this 2-category by $\gerbes (A,B)$.

\subsubsection{Classification of $(A,B)$-gerbes.}
\label{sec:classification-a-b}

Once again, consider the special case $A'=A$ and $B'=B$, with $f$
and $g$ being the respective identity maps. Then we speak of an
$(A,B)$-morphism, and in particular of a
$(A,B)$-\emph{equivalence} if the underlying functor
$\lambda\colon \gG\to \gG'$ is an equivalence in the usual sense.
$(A,B)$-\emph{equivalence} is an equivalence relation, and the
set of equivalence classes is
\begin{math}
  \HHH^2(X,A\to B)\,.
\end{math}
While this can be defined in general (see ref.\ \cite{MR0480515})
in the abelian case it turns out to coincide with the second
hypercohomology group with values in the complex $A\to B$ in the
standard sense (cf.\ \cite{math.AG/0301304}).

\subsubsection{The canonical $(f,g)$-morphism.}
\label{sec:canonical-f-g}

Given a commutative diagram of group homomorphisms as above,
there is a canonical $(f,g)$-morphism
\begin{equation*}
  (f,g)_\ast \colon \gerbes (A,B) \lto \gerbes (A',B')\,,
\end{equation*}
given by extension of the band.  Namely, if $\gG$ is an
$A$-gerbe, there is a well-defined procedure giving an $A'$-gerbe
which we may call $f_\ast (\gG)$. Since locally $\gG_U\iso \tors
(A\vert_U)$, then $f_\ast (\gG)_U$ is simply given by standard
extension of the structure group.  Now, if $(\gG,\mu)$ is an
$(A,B)$-gerbe, then $(\gG,g_\ast\circ\mu)$ is an $(A,B')$-gerbe
and locally the functor $g_\ast\circ\mu$ will be isomorphic to
$g_\ast\circ\delta_\ast$ (see in particular the proof of
Thm.~\ref{thm:6} below for more details\footnote{This
  construction will not be used until
  sect.~\ref{sec:grstb-grstc-torsors} and it is only dependent on
  the arguments of sect.~\ref{sec:class-results}, in particular
  the proof of Thm..~\ref{thm:6}.}).  The latter will be
replaced, by commutativity induced from the commutative square of
group homomorphisms, by $\delta'_\ast \circ f_\ast$, which glues
back to a functor $\mu' \colon f_\ast (\gG)\to \tors (B')$.

This construction is universal in the sense that an
$(f,g)$-morphism can be written by the composition of
$(f,g)_\ast$ followed by a unique (up to equivalence)
$(A',B')$-morphism.

An alternative characterization of $(A,B)$-gerbes will appear in
sect.~\ref{sec:class-results}, when we discuss 2-gerbes bound by
complexes.

\subsection{Examples}
\label{sec:examples}

The following are few examples of Gerbes bound by complexes of
length 2 which are relevant from the point of view of extending
differential geometric structures to gerbes.

We will first review the definition of connection---or
\emph{connective structure}---on a $\sho{X}^\unit$-gerbe
according to Brylinski and McLaughlin (see, e.g.\ 
\cite{brymcl:deg4:I,brymcl:deg4:II}, or \cite{bry:loop} for the
smooth case).
\begin{definition}\label{def:3}
  Let $\gerbe{G}$ be a $\sho{X}^\unit$-gerbe.
  A \emph{connective structure} $\conn{}$ on $\gerbe{G}$ is the
  datum of a $\shomega[1]{U}$-torsor $\conn(x)$ for any object
  $x\in \gerbe{G}_U$, where $U\subset X$, subject to the
  following conditions.
  \begin{enumerate}
  \item For every isomorphism $f : x \to y$ in $\gerbe{G}_U$
    there is an isomorphism
    \begin{equation*}
      f_* \equiv \conn(f) : \conn(x) \lto \conn(y)
    \end{equation*}
    of $\shomega[1]{U}$-torsors. In particular, if $f\in
    \Aut(x)\iso \sho{X}^\unit\vert_U$, we require:
    \begin{equation}
      \label{eq:12}
      \begin{aligned}
        f_* : \conn(x) & \lto \conn(x)\\
        \nabla &\longmapsto \nabla +\d \log f
      \end{aligned}
    \end{equation}
    where $\nabla$ is a section of $\conn(x)$.
  \item If $g: y \to z$ is another morphism in $\gerbe{G}_U$,
    then $(gf)_* \iso g_*f_*$.
  \item The correspondence must be compatible with the
    restriction functors and natural transformations. Namely, if
    $\imath^*: \gerbe{G}_U \to \gerbe{G}_V$ is the restriction
    functor corresponding to the morphism $\imath : V \to U$ in
    $\cat{C}/X$, then there is a natural isomorphism
    $\alpha_\imath: \imath^*\circ \conn \Rightarrow \conn{}\circ
    \imath^*$ such that the diagram:
    \begin{equation*}
      \xymatrix@+1pc{%
        \imath^*\conn(x) \ar[d]_{\imath^*(f_*)}
        \ar[r]^{\alpha_\imath(x)} &
        \conn(\imath^*x) \ar[d]^{(\imath^*f)_*}\\
        \imath^*\conn(y) \ar[r]^{\alpha_\imath(y)} &
        \conn(\imath^*y)
      }
    \end{equation*}
    commutes. Moreover given $\jmath: W\to V$ and the
    corresponding $\alpha_\jmath$, there must be the obvious
    pentagonal compatibility diagram with the natural
    transformations $\phi_{\imath,\jmath}: \jmath^*\imath^*\to
    (\imath\jmath)^*$ arising from the structure of fibered
    category over $X$. That is, given the object $x$, we have the
    commutative diagram:
    \begin{equation*}
      \xymatrix{
        \conn_W (\jmath^*\imath^*x)
        \ar[d]_{\phi_{\imath,\jmath}(x)_*}\ar[r]^{\alpha_\jmath} &
        \jmath^* \conn_V(\imath^*x)
        \ar[r]^{\jmath^*\alpha_\imath}&
        \jmath^*\imath^* \conn_U(x)
        \ar[d]^{\phi_{\imath,\jmath}(\conn_U(x))}\\
        \conn_W((\imath\jmath)^*x)\ar[rr]^{\alpha_{\imath\jmath}} & &
        (\imath\jmath)^*\conn_U(x)
      }
    \end{equation*}
    mapping to a corresponding one with $y$.
  \end{enumerate}
\end{definition}
The following is a reformulation of the conditions in
Definition~\ref{def:3}:
\begin{proposition}\label{prop:1}
  A connective structure on the $\sho{X}^\unit$-gerbe $\gerbe{G}$
  amounts to the datum of a structure of gerbe bound by the
  complex
  \begin{equation*}
    \Gamma :
    \sho{X}^\unit \xrightarrow{\d\log}  \shomega[1]{X}\,.
  \end{equation*}
\end{proposition}
\begin{proof}
  That the various conditions in Definition~\ref{def:3} define a
  cartesian functor
  \begin{equation*}
    \conn : \gerbe{G} \lto \tors (\shomega[1]{X})
  \end{equation*}
  is just a matter of unraveling the definition of cartesian
  functor. Moreover, eq.~\eqref{eq:12} implies that $\conn$ is in
  fact a $\d\log$-morphism.
\end{proof}
According to the general results, $\sho{X}^\unit$-gerbes with
connective structure are classified by the hypercohomology group
\begin{equation*}
  \HHH^2(X,\sho{X}^\unit \xrightarrow{\d\log}  \shomega[1]{X})\,.
\end{equation*}
Via the quasi-isomorphisms:
\begin{equation*}
  \bigl(\sho{X}^\unit \xrightarrow{\d\log}  \shomega[1]{X}\bigr)[-1]
  \lqi
  \bigl(\ZZ(2) \lto \sho{X} \overset{\d}{\lto} \shomega[1]{X}\bigr)
  \lqi \delz{2}\,,
\end{equation*}
where
\begin{math}
  \delz{k} 
\end{math}
is the weight $k$ Deligne complex, we have that the classifying
group is isomorphic to the \emph{Deligne cohomology group}
\begin{equation*}
  \delH[3]{X}{\ZZ}{2}\,.
\end{equation*}

\subsection{Further examples}
\label{sec:further-examples}

Several variations on the theme established in
Definition~\ref{def:3} and Proposition \ref{prop:1} have been
considered, typically by providing the necessary modifications in
Definition~\ref{def:3}. Following the idea embodied in
Proposition~\ref{prop:1} they can be restated in terms of gerbes
bound by a complex.

In ref.\ \cite{MR2142353} we have introduced a notion of
hermitian structure and a variant of connective structure valued
in the Hodge filtration.  We consider these examples next.

\subsubsection{Hermitian Structures.}
\label{sec:hermitian-structures}
Consider the complex:
\begin{equation*}
  \sho{X}^\unit \xrightarrow{\abs{\cdot}^2} \she[+]{X}
\end{equation*}
where $\she[+]{X}$ is the sheaf of \emph{smooth} functions valued
in $\RR_{>0}$, the connected component of $1$ in $\RR^\unit$.  A
$\sho{X}^\unit$-gerbe $\gG$ is said to have a \emph{hermitian
  structure} (cf.\ ref.\ \cite[Definition 5.2.1]{MR2142353}) if
it has the structure of a gerbe bound by $(\sho{X}^\unit,
\she[+]{X})$.

Classes of equivalences of $\sho{X}^\unit$-gerbes equipped with
hermitian structures are therefore classified by the group
\begin{equation*}
  \HHH^2(X,\sho{X}^\unit \xrightarrow{\abs{\cdot}^2} \she[+]{X})
  \iso 
  \dhhH[3]{X}{1}\,.
\end{equation*}
Recall that the latter is the third \emph{Hermitian Deligne
  cohomology group} of weight 1, and the isomorphism is induced
by the quasi-isomorphism
\begin{equation*}
  \bigl(\ZZ(1) \to \sho{X} \xrightarrow{\Re}  \she[0]{X}\bigr)
  \lqi
  \bigl(\sho{X}^\unit \xrightarrow{\abs{\cdot}^2}
  \she[+]{X}\bigr)[-1]\,,
\end{equation*}
where the first is the corresponding Hermitian Deligne complex.

\subsubsection{$F^1$-connections.}
\label{sec:f1-connections}
A slight modification of the notion of connective structure
recalled in sect.~\ref{sec:examples} is to consider the length 2
complex (\cite{MR2142353}):
\begin{equation*}
  \sho{X}^\unit \xrightarrow{\del\log} F^1\!\sha[1]{X}\,.
\end{equation*}
Note that $F^1\!\sha[1]{X}=\sha[1,0]{X}$, so this is called a
``type $(1,0)$ connective structure'' in \cite{MR2142353}.

\subsubsection{Compatibility.}
\label{sec:compatibility}
We have the obvious map $\del\log \colon \she[+]{X} \to
F^1\!\sha[1]{X}$ and the morphism of complexes
\begin{equation*}
    \xymatrix{
      \sho{X}^\unit \ar[r]^{\abs{\cdot}^2} \ar@{=}[d] &
      \she[+]{X} \ar[d]^{\del\log} \\
      \sho{X}^\unit \ar[r]^{\del\log} & F^1\!\sha[1]{X}}
\end{equation*}
The notion of compatibility between a hermitian and a type
$(1,0)$ connective structures on $\gG$ amounts to an
$(\id,\del\log)$-morphism.  In fact, it is the canonical one in
the sense of sect.~\ref{sec:canonical-f-g}.  The equivalence with
\cite[5.3.2]{MR2142353}, is merely a question of unraveling
Definition~\ref{def:2} for the case at hand. The classifying
group was identified in \cite{MR2142353} with $\dhhH[3]{X}{1}$,
computed using the complex $\dhh{X}{1}$.
\begin{remark}
  \label{rem:4}
  It was found that the notion of connection compatible with a
  given hermitian structure as defined in \loccit not the same as
  the one used by Brylinski and others (see, e.g.\ 
  \cite[Proposition 6.9 (1)]{bry:quillen}).  Here we can further
  elucidate the remarks at the end of \cite{MR2142353} by
  pinpointing the geometric difference: the notion of
  compatibility used by Brylinski involves solely the structure
  of $(\she[+]{X},\she[1]{X}(1))$-gerbe, whereas the definition
  we put forward uses the notion of \emph{morphism} of gerbes
  bound by a complex.  The latter remembers, so to speak, the
  structure of $\sho{X}^\unit$-gerbe.
\end{remark}

\section{Gerbes bound by complexes of length 3}
\label{sec:complexes-length-3}

\subsection{$(B,C)$-torsors}
\label{sec:b-c-torsors}

First, recall that for a given complex $B\overset{\sigma}{\lto}C$
of non-necessarily abelian groups, an $(B,C)$-\emph{torsor} (see
\cite{MR546620,MR92m:18019}) is a pair $(P,s)$ where $P$ is an
$B$-torsor and $s$ a section of $\sigma_*(P)\eqdef P\cprod{B}C$.
A morphism between two pairs $(P,s)$ and $(P',s')$ is a morphism
$f\colon P \to P'$ of $B$-torsors such that $\sigma_*(f)(s)=s'$.
With these definitions the $(B,C)$-torsors form a category, in
fact a gerbe, $\tors (B,C)$, and we denote by $\HHH^1(X,B\to C)$
the set of isomorphism classes. There is an obvious forgetful
functor $\tors (B,C) \lto \tors (B)$, and a corresponding map of
cohomology sets $\HHH^1(X,B\to C)\lto \H^1(X,B)$.

When $B$ and $C$ are abelian, which is the case of interest here,
the cohomology set classifying isomorphism classes of
$(B,C)$-torsors is isomorphic to the standard hypercohomology
group.

Suppose we are given a map of complexes
\begin{equation*}
    \xymatrix{%
      B\ar[r]^\sigma \ar[d]^g & C \ar[d]^h \\
      B' \ar[r]^{\sigma'} & C'
    }
\end{equation*}
then we obtain a functor
\begin{equation*}
  (g,h)_* \colon \tors(B,C) \lto \tors(B',C')\,,
\end{equation*}
which is defined as follows. To an object $(P,s)$ of $\tors(B,C)$
we associate the pair $(g_*P,h_*(s))$, where $g_*P =
P\cprod{B}B'$. This is well defined, since 
\begin{math}
  \sigma'_*g_*P \caniso h_*\sigma_* P\,.
\end{math}
Then it is immediate to verify that morphisms
\[
(P,s)\to (P',s')
\]
in $\tors(B,C)$ are brought to morphisms in $\tors(B',C')$.

The following alternative characterization will be useful in the
following. Using \cite[III.1.6.1]{MR49:8992}, it is easily seen
that the structure of $(B,C)$-torsor on $P$ corresponds to the
datum of a $C$-equivariant map:
\begin{equation}
  \label{eq:13}
\begin{aligned}
  \sigma_\ast (P) & \lto \SheafHom_B (P,C) \\
  t &\longmapsto [s \mapsto t^{-1}\sigma_\ast s]
\end{aligned}
\end{equation}
where $\SheafHom_B$ denotes (right) $B$-equivariant maps, and $C$
is considered as a right $B$-space via $\sigma$.

\subsection{$(A,B,C)$-gerbes}
\label{sec:a-b-c}

Let $A\overset{\delta}{\to} B\overset{\sigma}{\to} C$ be a
complex of abelian groups on $\cat{C}/X$, and let $p\colon
\gG\to \cat{C}/X$ be a gerbe with band $A$.
\begin{definition}
  \label{def:4}
  We say that $\gG$ is bound by the complex $A\to B\to C$, or
  that is an $(A,B,C)$-gerbe, if there is morphism
  \begin{equation*}
    \Tilde\mu \colon \gG \lto \tors(B,C)
  \end{equation*}
  such that $\gG$ is an $(A,B)$-gerbe for the $\delta$-morphism
  defined by the composition of $\Tilde\mu$ with the forgetful
  functor $\tors(B,C)\to \tors(B)$.
\end{definition}
In other words, the structure of $(A,B,C)$-gerbe on $\gG$ is a
factorization of the morphism $\mu$ defining the structure of
$(A,B)$-gerbe through $\tors(B,C)$.  For an object $x\in \ob
\gG_U$, denote
\begin{equation*}
  \Tilde\mu (x) = (\mu (x),\nu (x))\,,
\end{equation*}
where $\mu = \mathit{forget}\circ \Tilde\mu$, and $\nu (x)$ is a
section of $\sigma_*(\mu(x))$.

Next, we can consider the notion of morphism of two such gerbes
along the same lines as for $(A,B)$-gerbes. Thus, let us be given
a morphism of complexes of abelian sheaves over
$\cat{C}/X$:
\begin{equation*}
    \xymatrix{%
      A \ar[r]^\delta \ar[d]^f & B \ar[r]^\sigma \ar[d]^g &
      C \ar[d]^h \\
      A'\ar[r]^{\delta'} & B' \ar[r]^{\sigma'} & C'
    }
\end{equation*}
Let $\gG$ and $\gG'$ be two gerbes bound by $(A,B,C)$ and
$(A',B',C')$, respectively.
\begin{definition}
  \label{def:5}
  An $(f,g,h)$-morphism from $\gG$ to $\gG'$ consists of:
  \begin{enumerate}
  \item an $f$-morphism $\lambda\colon \gG \lto \gG'$;
  \item a natural isomorphism of functors
    \begin{equation*}
      \Tilde\alpha \colon (g,h)_\ast\circ \Tilde\mu
      \Longrightarrow \Tilde\mu'\circ \lambda
    \end{equation*}
    from $\gG$ to $\tors(B',C')$ such that the composition
    (=pasting) $F'*\Tilde\alpha$ with the forgetful functor
    $F'\colon \tors(B',C') \lto \tors(B')$ is the natural
    isomorphism associated to an $(f,g)$-morphism as in
    Definition~\ref{def:2}.
  \end{enumerate}
\end{definition}
The second condition in the definition can be explained as
follows. Consider the diagram
\begin{equation*}
    \xymatrix@+1em{%
      \gG \ar[r]^-{\Tilde\mu}_-{}="tmu" \ar[d]_\lambda^{}="lambda" &
      \tors(B,C) \ar[r]^F \ar[d]^{(g,h)_\ast}
      \ar@{=>}@/^/ ^{\Tilde\alpha} "tmu";"lambda" &
      \tors(B) \ar[d]^{g_\ast} \\
      \gG' \ar[r]^-{\Tilde\mu'} &
      \tors(B',C') \ar[r]^{F'} &
      \tors(B')
    }
\end{equation*}
Pasting with $F'$ gives
\begin{align*}
  F'*\Tilde\alpha \colon F'\circ(g,h)_\ast\circ \Tilde\mu
  &\Longrightarrow F'\circ\Tilde\mu'\circ \lambda\\
  \intertext{that is,}
  F'*\Tilde\alpha \colon
  g_\ast\circ F\circ \Tilde\mu & \Longrightarrow \mu'\circ\lambda
\end{align*}
We require this to coincide with the isomorphism $\alpha$ in
Definition~\ref{def:2}.

Again, we call this morphism an \emph{equivalence}, or more
precisely, an $(f,g,h)$-equivalence, if so is the underlying
functor $\lambda\colon \gG\lto \gG'$. In particular, this is the
case when $A'=A$, $B'=B$, $C'=C$ and $f$, $g$, and $h$ are the
identity map, which we refer to as an
$(A,B,C)$-equivalence. Being equivalent in this sense is an
equivalence relation, and we have:
\begin{proposition}
  \label{prop:2}
  Classes of equivalences of $(A,B,C)$-gerbes are classified by
  the hypercohomology group
  \begin{equation*}
    \HHH^2(X,A\to B\to C)\,.
  \end{equation*}
\end{proposition}
\begin{proof}
  We will just sketch how to obtain the class corresponding to a
  gerbe $\gG$ on $\cat{C}/X$ bound by the complex $A\to
  B\to C$.  Let us proceed under the assumption that working with
  \cech\ cohomology is sufficient. Thus, let $(U_i\to X)_{i\in
    I}$ be a cover for $X$ and assume that $\gG$ is decomposed
  \cite{MR95m:18006} by the choice of objects $x_i\in \Ob
  \gG_{U_i}$ and morphisms $\phi_{ij} \colon x_j\vert_{U_{ij}}
  \to x_i\vert_{U_{ij}}$.
  
  For each object $x_i$ the functor $\Tilde\mu \colon \gG \lto
  \tors (B,C)$ gives us a pair $\Tilde\mu (x_i) = (\mu (x_i), \nu
  (x_i))$, where $\nu (x_i) \in\Gamma (\sigma_\ast (x_i))$. Then,
  from the morphism $\phi_{ij}$ we obtain the morphism of torsors
  \begin{equation*}
    (\phi_{ij})_\ast \equiv \mu (\phi_{ij})
    \colon \mu (x_j) \lto \mu (x_i)
  \end{equation*}
  so that
  \begin{equation}\label{eq:14}
    \nu (x_i) = \sigma_\ast ((\phi_{ij})_\ast) (\nu (x_j))\,.
  \end{equation}
  The decomposition $(x_i,\phi_{ij})$ of $\gG$ gives a cocycle
  $(a_{ijk})\in \mathrm{Z}^2((U_i\to X),A)$ in the usual way,
  \cite{MR95m:18006}, \cite[IV.3.5.1]{MR49:8992}. Furthermore,
  let $(s_i)_{i\in I}$ be a collection where $s_i$ is a section
  of the $B\vert_{U_i}$-torsor $\mu (x_i)$. It follows that a
  cochain $(b_{ij})$ with values in $B$ is defined by
  \begin{equation*}
    (\phi_{ij})_\ast (s_j) = s_i\,b_{ij}\,,
  \end{equation*}
  and the usual argument shows that
  \begin{equation}
    \label{eq:15}
    a_{ijk} = b_{ik}^{-1}b_{ij} b_{jk}\,.
  \end{equation}
  Now, since $\Tilde\mu (x_i)$ is a $(B,C)$-torsor, we have that
  \begin{equation*}
    \sigma_\ast (s_i) = \nu (x_i)\, c_i\,,
  \end{equation*}
  for an appropriate section $c_i$ of $C\vert_{U_i}$, for each
  $i\in I$. On one hand, this gives:
  \begin{equation*}
    \sigma_\ast ((\phi_{ij})_\ast (s_j)) = \nu (x_i)\, c_i\,
    \,\sigma (b_{ij})\,.
  \end{equation*}
  On the other hand, by functoriality we have
  \begin{align*}
    \sigma_\ast ((\phi_{ij})_\ast (s_j)) 
    & = \sigma_\ast (\mu(\phi_{ij})) (\sigma_\ast (s_j))\\
    & = \sigma_\ast ((\phi_{ij})_\ast) (\nu (x_j))\, c_j\,,
  \end{align*}
  and using~\eqref{eq:14} we finally obtain
  \begin{equation}
    \label{eq:16}
    c_i\, \,\sigma (b_{ij}) = c_j\,.
  \end{equation}
  Then~\eqref{eq:15}, \eqref{eq:16}, and the cocycle property for
  $(a_{ijk})$ give the desired 2-cocycle with values in the
  complex $A\to B\to C$.
\end{proof}
The alternative characterization of $(B,C)$-torsor at the end of
sect.~\ref{sec:b-c-torsors}, and the technique used in the proof
of the proposition can be put together to provide the following
alternative characterization of the notion of $(A,B,C)$-gerbe.

Let $A\overset{\delta}{\to} B\overset{\sigma}{\to} C$ be a
complex of abelian groups over $\cat{C}/X$.
\begin{lemma}
  \label{lem:1}
  The structure of $(A,B,C)$-gerbe on $\gG\lto \cat{C}/X$
  is equivalent to the following data:
  \begin{enumerate}
  \item $\mu\colon \gG \lto \tors(B)$ making $\gG$ into an
    $(A,B)$-gerbe;
  \item for each object $x\in \Ob\gG_U$ a map $\nu (x) \colon
  \mu(x)\lto C\vert_U$ such that:
  \begin{enumerate}
  \item $\nu (x) (s\,b) = \nu (x) (s)\,\sigma (b)$ for each
    section $s$ of $\mu (x)$ and $b$ of $B\vert_U$;
  \item for each morphism $f\colon x\lto y$ in $\gG_U$ a
    commutative diagram
    \begin{equation*}
        \xymatrix@-1em{%
          \mu (x) \ar[rr]^{\mu (f)} \ar @/_/[dr]_{\nu(x)}
          && \mu (y) \ar @/^/[dl]^{\nu(y)} \\
          & C\vert_U &
        }
    \end{equation*}
  \end{enumerate}
  \end{enumerate}
\end{lemma}
\begin{proof}
  The existence of the map $\nu (x)$ is simply a consequence of
  the existence of a section $\nu (x)$ of $\sigma_\ast (\mu(x))$
  in the structure of $(B,C)$-torsor of $\mu (x)$ determines a
  morphism $\mu (x) \lto C\vert_U$ according to~\eqref{eq:13}.
  
  The commutativity of the diagram follows then from the fact
  that the structure of $(B,C)$-torsor of $\mu(x)$ implies that
  $\nu (y) = \sigma_\ast \mu (f) (\nu (x))$.
\end{proof}
A different characterization of $(A,B,C)$-gerbes in terms of
torsors over a morphism of \gr-stacks will appear in
sect.~\ref{sec:grstb-grstc-torsors}, when we will be discussing
2-gerbes bound by complexes (of \gr-stacks).

\subsection{Examples: Curvings}
\label{sec:examples:-curvings}

The main example we want to consider, is that of a curving on a
$\sho{X}^\unit$-gerbe $\gG$ equipped with a connective structure.
The concept, introduced by Brylinski (\cite{bry:loop}), but
attributed to Deligne, is the analogous of the curvature of a
connection on a line bundle.

$\gG$ possesses a connective structure if it is a gerbe bound by
$\sho{X}^\unit \xrightarrow{\d\log} \shomega[1]{X}$. We can move
one step forward and consider instead the longer complex:
\begin{equation}
  \label{eq:17}
  \sho{X}^\unit \xrightarrow{\d\log}
\shomega[1]{X} \overset{\d}{\lto} \shomega[2]{X}\,.
\end{equation}
\begin{definition}
  \label{def:6}
  A \emph{curving} on $\gG$ is the structure of gerbe bound by
  the complex~\eqref{eq:17}\,.
\end{definition}
According to Lemma~\ref{lem:1}, a curving on a gerbe $\gG$ with
connective structure $\conn$ will be given by a map 
\begin{equation*}
  \curv (x) \colon \conn(x) \lto \shomega[2]{U}
\end{equation*}
for each object $x\in \ob \gG_U$, and open $U\to X$, such that
\begin{equation*}
  \curv (x) (\nabla + \alpha) = \curv(x)(\nabla) + \d\alpha\,,
\end{equation*}
where $\nabla$ is a section of $\conn (x)$ and $\alpha$ is a
section of $\shomega[1]{U}$. Moreover, if $f\colon x\lto y$ is a
morphism in $\gG_U$, then the commutative diagram in
Lemma~\ref{lem:1}  translates into
\begin{equation*}
  \curv (y)(f_\ast(\nabla)) = \curv  (x)(\nabla)\,.
\end{equation*}
By direct comparison, we can see that these are exactly the
properties of the curving listed in \cite{bry:loop}, hence our
definition agrees with the one in \loccit

It follows from the classification result above that we have a
gerbe $\gG$ equipped with connective structure and curving
defines a class in the hypercohomology group:
\begin{equation*}
  \HHH^2(X, \sho{X}^\unit \xrightarrow{\d\log}
  \shomega[1]{X} \overset{\d}{\lto} \shomega[2]{X})\iso
  \delH[3]{X}{\ZZ}{3}\,.
\end{equation*}
The isomorphism with the Deligne cohomology group follows from
the quasi-isomorphisms:
\begin{equation*}
  \bigl(\sho{X}^\unit \xrightarrow{\d\log}
  \shomega[1]{X} \overset{\d}{\lto}
  \shomega[2]{X}\bigr)[-1]
  \lqi
  \bigl(\ZZ(3) \lto \sho{X} \overset{\d}{\lto} \shomega[1]{X}
  \overset{\d}{\lto} \shomega[2]{X}\bigr)\,,
\end{equation*}
the complex on the right hand side being $\delz{3}$.

\section{2-Gerbes: main definitions}
\label{sec:2-gerbes}

In this section we review some basic definitions and relevant
facts about 2-gerbes here. The standard reference is
\cite{MR95m:18006}, which should be referred to for a complete
treatment.

Recall that a 2-gerbe is a 2-stack, in particular a fibered
2-category, satisfying local non-emptiness and connectivity
requirements generalizing those of a gerbe. The general
definition of fibered 2-categories can be found in
\cite{MR0364245}. Analogously to \loccit, we will assume that
given a fibration $p \colon \tgG \lto \twocat{S}$ of
2-categories, the base 2-category is in effect a category
regarded as a \emph{discrete} 2-category---namely, one with all
2-arrows being identities. In other words,
$\twocat{S}=\operatorname{2-\twocat{Cat}}(\cat{S})$, where
$\cat{S}$ is a category. To avoid overburdening our notation, we
will simply write our fibrations as $p\colon \tgG\lto \cat{S}$,
without risk of confusion.  In the following, the category
$\cat{S}$ will in fact be the site $\cat{C}/X$, with all our
standing assumptions concerning $\cat{C}/X$ to be kept for
2-gerbes as well.

\subsection{2-Stacks}
\label{sec:2-stacks}

A 2-stack is a fibered 2-category $p\colon \tgG \lto \cat{S}$
such that:
\begin{enumerate}
\item\label{item:1} 1-arrows and 2-arrows can be glued, a fact
  that can be succinctly stated by saying that for any two
  objects $x,y\in \Ob \tgG_U$ over $U\in \Ob \cat{S}$, the
  fibered category $\catHom_U(x,y)$ is stack over
  $\cat{S}/U$;
\item\label{item:2} Objects can be glued, namely 2-descent on
  objects holds.
\end{enumerate}
(A pre-2-stack is a fibered 2-category satisfying only the first
condition above.)

Without entering into too many details, it is worthwhile making
the gluing condition on objects more explicit.  Thus, let $U$ be
an object of $\cat{S}$, and let $(U_i\to U)$ be a cover as usual.
The assignment of 2-descent data over $U$ is the assignment of a
collection of objects $x_i\in \ob \tgG_{U_i}$ such that there is
a 1-arrow:
\begin{equation*}
  \phi_{ij}\colon x_j\lto x_i
\end{equation*}
over $U_{ij}$ and a 2-arrow (in fact, a 2-isomorphism):
\begin{equation*}
    \xymatrix@R-1pc{%
      & x_j \ar@/^/[dr]^{\phi_{ij}}
      \ar@{=>}[d]^-{\alpha_{ijk}} & \\
      x_k \ar@/^/[ur]^{\phi_{jk}} 
      \ar[rr]_{\phi_{ik}} & & x_i
    }
\end{equation*}
over $U_{ijk}$ such that the following compatibility condition
holds:
\begin{equation*}
  \alpha_{ikl} \circ (\alpha_{ijk}\ast\phi_{kl})
  = \alpha_{ijl} \circ (\phi_{ij}\ast \alpha_{jkl})\,.
\end{equation*}
The assignment of the triple $(x_i, \phi_{ij}, \alpha_{ijk})$, is
called \emph{2-Descent data.}  Condition~\ref{item:2} above then
means that there exists an object $x\in \Ob\tgG_U$ with 1-arrows
\begin{equation*}
  \psi_i\colon x_i \lto x
\end{equation*}
and 2-isomorphisms
\begin{equation*}
    \xymatrix@R-1pc{%
      & x_i \ar@/^/[dr]^{\psi_i} 
       \ar@{=>}[d]^{\chi_{ij}}  & \\
      x_j \ar@/^/[ur]^{\phi_{ij}} \ar[rr]_{\psi_j} && x
    }
\end{equation*}
satisfying the now obvious compatibility conditions with the
isomorphisms $\alpha_{ijk}$. This is referred to by saying that
the 2-descent data is effective.

\subsection{2-Gerbes}
\label{sec:2-gerbes-1}

In words, a 2-gerbe $\tgG\lto \cat{S}$ is a 2-stack in
2-groupoids which is locally non-empty and connected.  A detailed
account of several variants of this definition of a 2-gerbe is
given in the text~\cite{MR95m:18006}.  Following \loccit, the
properties characterizing a 2-Gerbe are the following:
\newcounter{saveenum}
\begin{enumerate}
\item\label{item:3} $\tgG$ is \emph{locally non-empty:} assuming
  $\cat{S}=\cat{C}/X$, there exists a cover $U\to X$ such
  that $\Ob \tgG_U$ is not empty.
\item\label{item:4} $\tgG$ is \emph{locally connected:} for each
  $x,y\in \Ob \tgG_U$, for some object $U$ of $\cat{S}$, there
  exists a cover $\phi:V\to U$ such that the set of arrows from
  $x_V$ to $y_V$\footnote{Note that given $\phi : V\to U$ and an
    object $x$ above $U$ thanks to the axioms of a fibered
    2-category we can speak of ``the'' object $x_V$ above $V$
    with an arrow $x_V\to x$ above $\phi$ up to 2-equivalence.}
  is not empty.
  \setcounter{saveenum}{\value{enumi}}
\item\label{item:5}
  \emph{1-arrows are weakly invertible:} for
  any 1-arrow $f:x\to y$ in $\gG_U$, $U\in \Ob\tgG$, there is an
  inverse $g:y\to x$ up to two 2-arrows.
\item\label{item:6} 2-arrows are (strictly) invertible in $\tgG_U$.
\end{enumerate}
There are different equivalent forms of the last two axioms, as
well as local versions of all four to be obtained by considering
coverings of $U$, see~\cite{MR95m:18006} for more details. Here
we only quote the fact that condition~\ref{item:5} above is
equivalent (if condition~\ref{item:6} is also satisfied) to:
\begin{enumeratep}
  \setcounter{enumi}{\value{saveenum}}
\item\label{item:7} Given two 1-arrows $f:x\to y$ and $g: x \to
  z$ in $\tgG_U$, there exists a 1-arrow $h: y \to z$ and a
  2-arrow $\alpha : h\circ f \Rightarrow g$.
\end{enumeratep}
Finally, a note of caution: although the \emph{stack}
$\catHom_{U}(x,y)$ is locally non empty by
condition~\ref{item:4}, in general it will not be connected, so
that condition~\ref{item:5} does not quite imply that
$\catHom_{U}(x,y)$ is a gerbe. This is the case when $x=y$ for
fully \emph{abelian} 2-gerbes, to be discussed below.

\subsubsection{\emph{Gr}-stacks of automorphisms}
\label{sec:gr-stacks-autom}

To conclude these remarks of preparatory nature, let us briefly
discuss automorphisms of objects.

For any given object $x\in \Ob\tgG_U$, the stack $\catAut_U(x)$
of self-arrows of $x$ is a stack in groupoids equipped with a
strictly associative \emph{monoidal structure,} that is a functor
\begin{math}
  \catAut_U(x) \times \catAut_U(x) \lto \catAut_U(x)
\end{math}
implementing a product law on $\catAut_U(x)$. It follows from the
2-gerbe axioms that $\catAut_U(x)$ admits a choice of inverses,
compatible with descent, hence it is a \emph{group-like} stack,
or \emph{\gr-stack,} for short, cf.\
\cite{MR93k:18019,MR95m:18006,MR0338002}.

Analogously to the gerbe case, if $\grstA$ is a fixed \gr-stack on
$\cat{S}$, we define a 2-$\grstA$-gerbe to be a 2-gerbe $\tgG$
over $\cat{S}$ such that for every object $x\in \Ob\tgG_U$ there
is an equivalence
\begin{equation*}
  a_x \colon \catAut_U(x) \lisoto \grstA\vert_U\,.
\end{equation*}

\subsection{Abelian 2-gerbes}
\label{sec:abelian-2-gerbes}

A 2-$\grstA$-gerbe to be \emph{abelian} if the equivalences $a_x$
introduced above are natural in the sense specified
in~\cite[Definition 4.13]{MR95m:18006}. As shown in \loccit, this
has the consequence that $\grstA$ is \emph{braided,} that is,
there is a commutativity functor for the monoidal structure.

An \emph{additional} commutativity condition is to assume that
\begin{equation*}
  \grstA = \tors(A)\,,
\end{equation*}
for a sheaf of abelian groups $A$ over $\cat{S}$. (Since $A$ is
abelian, this is a \gr-stack under the standard contracted product
of $A$-torsors.)

As explained in \loccit, these two requirements have the
consequence that the \gr-stack $\catAut_U(x)$ is a gerbe over
$\cat{S}/U$, and in fact a neutral one, i.e.\ it is equivalent to
$\tors(A\vert_U)$, since it has the global object $\id_x$.
Automorphisms of 1-arrows are then equivalent to sections of the
sheaf of groups $A$, as in~\cite{brymcl:deg4:I}. If both
commutativity conditions hold, we commit a mild abuse of language
and say that the 2-gerbe $\tgG$ is \emph{bound} by the sheaf of
abelian groups $A$, or that it is a 2-$A$-gerbe, dropping the
typographical reference to the \gr-stack $\grstA$.

It is by now standard that the \emph{fully} abelian 2-gerbes, or
2-$A$-gerbes, are classified up to equivalence by the ordinary
cohomology group $\H^3(X,A)$.

In what follows we will limit our consideration to abelian
2-gerbes which are not, however, necessarily fully abelian.

\subsubsection{Morphisms.}
\label{sec:morphisms}

As noted, a morphism between two 2-gerbes $\tgG$ and $\tgH$ is a
cartesian 2-functor $F\colon \tgG \lto \tgH$ between the
underlying 2-stacks.

Suppose that $\tgG$ is a 2-$\grstA$-gerbe and $\tgH$ is a
2-$\grstB$-gerbe, and $\lambda \colon \grstA\lto \grstB$ is a
morphism of \gr-stacks, where we assume both $\grstA$ and
$\grstB$ at least braided. By analogy with the case of gerbes, we
will call $F$ a $\lambda$-\emph{morphism} if the obvious
commutative diagrams (up to 2-isomorphism, this time) are
satisfied. In particular, this means that $F$ must be compatible
with the morphisms $a_x$ in the sense that we have the following
diagram:
\begin{equation*}
    \xymatrix@C+1pc{
      \catAut(x) \ar[r] \ar[d]_{a_x} & \catAut(F(x))
      \ar[d]^{b_x}_(.4){}="beta"  \\
      \grstA\vert_U \ar[r]_\lambda^(.6){}="lambda"
      \ar@{=>}@/_0.4pc/ _{\nu_x}"beta";"lambda" & \grstB\vert_U
    }
\end{equation*}
for an appropriate isomorphism $\nu_x$.

In particular, we are interested in the situation where a
homomorphism $\delta \colon A \lto B$ of abelian groups is given,
and $\lambda=\delta_\ast$ is simply the induced functor:
\begin{equation*}
  \delta_\ast\colon \tors(A) \lto \tors(B)\,.
\end{equation*}
between the corresponding \gr-stacks.  In this case we will refer
to $F$ as a $\delta$-morphism, with a mild abuse of language. The
salient property of a $\delta$-morphism in this sense is that if
a section $a\in A\vert_U$ corresponds to an automorphism of a
1-arrow $f$ of $\tgG_U$, then the corresponding automorphism of
$F(f)$ in $\tgH_U$ will be $\delta (a)\in B\vert_U$.

\subsubsection{Classification.}
\label{sec:classification}

As already mentioned, a 2-$A$-gerbe is classified by an element
of the (ordinary) cohomology group $\H^3(X,A)$: Let us briefly
recall here the well-known local calculation leading to the
classification.

For simplicity, let us remain in the \cech setting, so let us
once again consider a cover $(U_i\to X)_{i\in I}$ of $X$. Now,
given a $2$-gerbe $\tgG$, let us choose a decomposition by
selecting a collection of objects $x_i$ in $\tgG_{U_i}$. There is
a $1$-arrow
\begin{displaymath}
  \phi_{ij}\colon x_j\to x_i
\end{displaymath}
between their restrictions to $\tgG_{U_{ij}}$.  Then
axiom~\ref{item:7} in sect.~\ref{sec:2-gerbes-1}, and the
abelianness assumptions imply that there exist 2-arrows such
that:
\begin{displaymath}
  \alpha_{ijk} \colon \phi_{ij}\circ \phi_{jk} \Longrightarrow
  \phi_{ik}\,.
\end{displaymath}
Over a $4$-fold intersection $U_{ijkl}$, we have two $1$-arrows
\begin{math}
  \phi_{ij}\circ \phi_{jk}\circ \phi_{kl}\colon x_l\to x_i
\end{math}
and
\begin{math}
  \phi_{il}\colon x_l\to x_i
\end{math}
and between them \emph{two} $2$-arrows, namely
\begin{math}
  \alpha_{ijl} \circ (\id_{\phi_{ij}} \ast\, \alpha_{jkl})
\end{math}
and
\begin{math}
  \alpha_{ikl} \circ (\alpha_{ijk} \ast \id_{\phi_{kl}})\,.
\end{math}
Since $2$-arrows are strictly invertible, it follows again from
the axioms that there exists a section $a_{ijkl}$ of
$\sho{X}^\unit$ over $U_{ijkl}$ such that
\begin{equation}
  \label{eq:57}
  \alpha_{ijl} \circ (\id_{\phi_{ij}} \ast \alpha_{jkl})
  = a_{ijkl}\circ
  \alpha_{ikl} \circ (\alpha_{ijk} \ast \id_{\phi_{kl}})\,.
\end{equation}
This section is a $3$-cocycle and the assignment $\tgG
\mapsto [a]$ gives the classification isomorphism.

\section{2-Gerbes bound by a complex}
\label{sec:2-gerbes-bound}

\subsection{$\grstB$-torsors}
\label{sec:grstackb-torsors}

The notion of \emph{torsor} under a \gr-stack will play a
significant role below.  The definition has been given in full
generality in~\cite[6.1]{MR92m:18019}, and~\cite{MR93k:18019}, so
here we will confine ourselves to only recall the main points.

Let $\grstB$ be a \gr-stack on $\cat{C}/X$. 
Briefly, a stack in groupoids $\stP$ will be a (right)
$\grstB$-torsor if there is a morphism of stacks
\begin{equation*}
  m \colon \stP \times \grstB \lto \stP
\end{equation*}
compatible with the group law of $\grstB$ in the sense
specified in \loccit, and such that the morphism
\begin{equation*}
  \Tilde m = (\mathrm{pr}_1,m) \colon
  \stP \times \grstB \lto \stP \times \stP
\end{equation*}
is an equivalence. As in \loccit, there will be an associativity
natural isomorphism:
\begin{equation*}
  \mu_{x,b,b'}\colon (x\cdot b)\cdot
  b'\lisoto x\cdot (b\cdot b')\,,
\end{equation*}
where $x\cdot b$ stands for $m(x,b)$. This isomorphism will have
to satisfy the standard pentagon diagram.

Having so far defined what ought to be called a
\emph{pseudo-}torsor, we need to complete the definition by
adding the condition that there exists a cover $U\to X$ such that
the fiber category $\stP_U$ is non-empty.

There are a few constructions for $\grstB$-torsors that are
generalizations of well-known ones for standard torsors which we
are going to recall now: cocycles and contracted products.

\subsubsection{Contracted product of torsors.}
\label{sec:contr-prod-tors}

The notion of contracted product for torsors over a \gr-stack is
introduced in~\cite[\S 6.7]{MR92m:18019}.

If $\stP$ (resp.\ $\stQ$) is a right (resp.\ left)
$\grstB$-torsor, the contracted product $\stP\cprod{\grstB}\!
\stQ$ is defined as follows. The objects are pairs $(x,y)\in \ob
\stP\times \stQ$.  A morphism $(x,y)\to (x',y')$, however, is an
equivalence classes of triples $(f,b,g)$, where $b\in \ob
\grstB$, and $f\colon x\cdot b\to x'$ and $g\colon y\to b\cdot
y'$ are morphisms of  $\stP$ and $\stQ$, respectively. Two
triples $(f,b,g)$ and $(f',b',g')$ are equivalent if there is a
morphism $\beta\colon b\to b'$ in $\grstB$ such that $f=f'\circ
(x\cdot \beta)$ and $g'=(\beta\cdot y')\circ g$.

Properties analogous to the familiar ones for ordinary torsors
hold. For example, one has the isomorphism
\begin{equation*}
  (x\cdot b, y) \lisoto (x,b\cdot y)\,,
\end{equation*}
given by the pair $(\id_{x\cdot b},b,\id_{b\cdot y})$.

In the following we will be considering braided (and in fact,
Picard) \gr-stacks exclusively, hence the distinction between
left and right-torsor will not matter. In principle, by analogy
with the case of standard torsors over an abelian group we could
dispense with the notation for the contracted product and denote
the product with the symbol $\stP\otimes \stQ$, instead. We will
not do so, however.

\subsubsection{Cocycles.}
\label{sec:cocycles}
A torsor $\stP$ over a (not necessarily braided) \gr-stack
$\grstB$ can also be characterized by a cocycle with respect to a
cover.

Given a cover $(U_i\to X)_{i\in I}$, the torsor $\stP$ has
non-empty fiber categories over it. Thus choose objects $x_i\in
\Ob\stP_{U_i}$. Since by definition $\stP$ is locally
(i.e.\ over the cover) equivalent to $\grstB$, it follows
that we can obtain isomorphism $x_j\lisoto x_i\cdot b_{ij}$,
where $b_{ij}$ is an object of $\grstB$ over $U_{ij}$, and
the isomorphism takes place in $\stP_{U_{ij}}$. (We are
systematically ignoring the isomorphisms resulting from the
pull-back functors.) By pulling back to $U_{ijk}$ we obtain a
1-cocycle with values in $\grstB$:
\begin{equation}
  \label{eq:23}
  \beta_{ijk}\colon b_{ij}\cdot b_{jk}\lisoto b_{ik}\,.
\end{equation}
The isomorphisms $\beta_{ijk}$ in $\grstB\vert_{U_{ijk}}$ turn
out to satisfy the obvious compatibility condition on quadruple
intersections $U_{ijkl}$, which we do not explicitly write here.
The pair $(b_{ij},\beta_{ijk})$ is the 1-cocycle with values in
the \gr-stack $\grstB$ determined by $\stP$.

\subsubsection{$\grstB$-torsors and $B$-gerbes.}
\label{sec:grstackb-torsors-b}

It arises from the general classification theory of 2-gerbes that
$\twotors (\grstB)$ is a 2-$\grstB$-gerbe.  Moreover, it follows
from the general discussion in~\cite[\S 7.2 and Proposition
7.3]{MR92m:18019} that if $\grstB=\tors (B)$, then $\tors
(\grstB)$ is equivalent to $\gerbes (B)$, the 2-gerbe of
$B$-gerbes over $X$.

It is possible to see this via the 1-cocycle pair~\eqref{eq:23}
as follows.  Recall that $\grstB=\tors (B)$ with $B$ abelian, so
we obtain a ``torsor cocycle'' in the sense
of~\cite{MR95b:18009}.  It follows that the groupoids $\tors
(B)\vert_{U_i}$ can be glued in the standard way to give a
$B$-gerbe.
\begin{remark}
  The argument just outlined is of course not specific to $B$
  being abelian. Upon replacing $\tors (B)$ with $\bitors (B)$
  everything works in general.
\end{remark}
\begin{remark}
  The 1-cocycle written above coincides with Hitchin's notion of
  ``gerbe data,''~\cite{MR1876068}. The latter lacks the
  categorical input, however.
\end{remark}

\subsection{Crossed modules of \gr-categories}
\label{sec:crossed-modules-gr}

It was observed above that the complex $\delta\colon A\lto B$ of
abelian groups ought to be considered as an abelian crossed module,
namely one where we impose strict commutativity on the associated
\gr-category. (That is, we demand it be strictly Picard.)

It turns out that a similar pattern holds in the case of a
\emph{crossed module of \gr-categories} in the sense
of~\cite[D\'efinition 2.2.8]{MR93k:18019}. It requires that there
exist additive functors
\begin{equation*}
  \lambda\colon \grstA \lto\grstB \,,\qquad
  \jmath \colon \grstB \lto \equ (\grstA)
\end{equation*}
such that the relations determined by the following diagrams
hold:
\begin{equation*}
  \xymatrix@-0.5pc{%
    & \grstB \ar@/^/[dr]^\jmath \ar@{=>}[d]^\mu  & \\
    \grstA \ar@/^/[ur]^\lambda \ar[rr]_{\imath_\grstA}&&
    \equ (\grstA)}
  \qquad
  \xymatrix{%
    \grstB \times \grstA \ar[r]^(0.6){\Hat\jmath}
    \ar[d]_{1_\grstB \times \lambda} &
    \grstA \ar[d]^\lambda_{}="lambda" \\
    \grstB \times \grstB
    \ar[r]_(0.6){\Hat\imath_\grstB}^{}="imath" & \grstB
    \ar@{=>}@/_/ _\nu "lambda";"imath"}
\end{equation*}
where $\equ (\grstA)$ denotes the \gr-stack of self-equivalences
of $\grstA$, $\imath_\grstA$ denotes the inner conjugation, and
the top and bottom horizontal arrows in the diagram to the right
are the actions of $\grstB$ on $\grstA$ and on itself induced by
$\jmath$ and the inner conjugation.

Now observe that requiring the resulting group law on
$\grstA\times \grstB$ to be commutative (up to natural
isomorphism), entails that both $\grstA$ and $\grstB$ are
braided, and that the action of $\grstB$ on $\grstA$ is trivial.
Thus, an abelian crossed module of \gr-categories will simply be
an additive functor
\begin{equation}
  \label{eq:20}
  \lambda\colon \grstA\lto \grstB\,,
\end{equation}
between braided \gr-categories. The same conclusions hold if we
replace \gr-cate\-go\-ries with \gr-stacks over $\cat{C}/X$. We will
also refer to~\eqref{eq:20} as a complex of (braided) \gr-stacks.

If both $\grstA$ and $\grstB$ have strict group laws, then they
are the \gr-categories associated to crossed modules, so we
obtain a ``crossed module of crossed modules,'' namely a
\emph{crossed square,} see~\cite{MR651845,MR93k:18019}.
Thus~\eqref{eq:20} reduces to the commutative square
\begin{equation}
  \label{eq:22}
  \vcenter{
    \xymatrix{%
      A \ar[d]_\delta \ar[r]^f & B \ar[d]^\sigma \\
      G \ar[r]_u & H}}
\end{equation}
where the vertical arrows are the crossed modules associated to
$\grstA$ and $\grstB$, respectively, and the horizontal arrows,
as well as the composite diagonal one, are also crossed
modules. There are other axioms, for which we refer the reader to
the treatment in \loccit We will not need them here, however,
because if both $\grstA$ and $\grstB$ are strictly commutative,
their associated crossed modules become complexes of abelian
groups, so that~\eqref{eq:22} becomes a commutative square of
homomorphisms of abelian groups, which is the situation we will
be interested in. Thus ``crossed square'' will be meant as a
synonym for a morphism of complexes of abelian groups.

\subsection{2-$(\grstA,\grstB)$-gerbes}
\label{sec:2-a-b}

We are now going to consider the analog of Definition~\ref{def:1}
for abelian 2-gerbes. We proceed by giving a direct
generalization of Definition~\ref{def:1}, where we replace the
complex $A\to B$ with the length 2-complex (that is a morphism)
of \gr-stacks, which we assume braided, heeding to the principle
that we climb the ladder of the higher algebraic structures by
promoting the \emph{coefficients} of cohomology from sheaves of
(abelian) groups, to \gr-stacks, etc.
\begin{definition}\label{def:7}
  A 2-gerbe bound by the complex~\eqref{eq:20} is a
  2-$\grstA$-gerbe $\tgG$ over $\cat{C}/X$, equipped
  with a $\lambda$-morphism:
  \begin{equation*}
    J \colon \tgG \lto \twotors (\grstB)\,.
  \end{equation*}
  A 2-gerbe bound by the complex~\eqref{eq:20} will be called a
  2-$(\grstA,\grstB)$-gerbe.  (Notice that $\twotors
  (\grstB)$ is a 2-$\grstB$-gerbe in an obvious way,
  hence the notion of $\lambda$-morphism makes sense.)
\end{definition}
If $\tgG$ is actually a 2-$A$-gerbe, and $\grstB=\tors (B)$,
where $B$ is a sheaf of abelian groups over $\cat{C}/X$,
with a homomorphism $\delta\colon A\lto B$, we call it a
2-$(A,B)$-gerbe, or a 2-gerbe bound by $A\to B$. (The morphism
$J$ in the definition is a $\lambda=\delta_\ast$-morphism.)

For a 2-$(A,B)$-gerbe, owing to the last remark in
sect.~\ref{sec:grstackb-torsors}, Definition~\ref{def:7} can be
recast in the form used in~\cite[Definition 5.5.1]{MR2142353} (in
a special case), which we state here as a lemma:
\begin{lemma}\label{lem:2}
  The datum of a 2-$(A,B)$-gerbe is equivalent to that of a
  Cartesian 2-functor
\begin{equation*}
  J \colon\tgG \lto \gerbes (B)
\end{equation*}
which is a $\delta$-morphism of 2-gerbes.
\end{lemma}

Morphisms of 2-gerbes bound by a complex of length 2 can be
defined by promoting Definition~\ref{def:2} to using braided
\gr-stacks and then (for those coming from abelian groups) using
Lemma~\ref{lem:2}. Specifically, analogously to what was done in
sect.~\ref{sec:gerbes-bound-complex}, consider the square of
\gr-stacks:
\begin{equation}
  \label{eq:21}
  \vcenter{\xymatrix@+1pc{
      \grstA \ar[r]^\lambda \ar[d]_\phi &
      \grstB \ar[d]^\psi_{}="psi" \\
      \grstA'\ar[r]_{\lambda'}^{}="lambdapr"
      \ar@{=>}@/_0.3pc/ _\jmath "psi";"lambdapr" &
      \grstB'}}
\end{equation}
\begin{definition}\label{def:8}
  A $(\phi,\psi)$-morphism $(F,\mu)\colon (\tgG,J)\lto
  (\tgG',J')$ consists of:
  \begin{enumerate}
  \item an $\phi$-morphism $F:\tgG\lto \tgG'$;
  \item\label{item:8} a natural transformation of 2-functors:
    \begin{equation*}
      \mu \colon \psi_{\ast}\circ J
      \Longrightarrow J'\circ F \colon \tgG \lto \tors
      (\grstB')\,, 
    \end{equation*}
    where $\psi_\ast\colon \tors (\grstB) \lto \tors
    (\grstB')$ is induced from $\psi$ in the obvious way.
  \end{enumerate}
\end{definition}
In particular, the special case where~\eqref{eq:21} is induced by
the morphism of complexes
\begin{equation*}
  (f,g)\colon (A,B,\delta) \lto
  (A',B',\delta')
\end{equation*}
of abelian groups will be referred to as an
$(f,g)$-morphism of the 2-$(A,B)$-gerbe $(\tgG,J)$ to the
2-$(A',B')$-gerbe $(\tgG',J')$. Using Lemma~\ref{lem:2},
condition~\ref{item:8} in Definition~\ref{def:7} says that we
have a natural transformation of 2-functors
\begin{equation*}
  \mu \colon g_{\ast\ast}\circ J
  \Longrightarrow J'\circ F \colon \tgG \lto \gerbes (B')\,,
\end{equation*}
where $g_{\ast\ast}\colon \gerbes (B)\lto \gerbes (B')$ is
induced from $g:B\lto B'$.

We speak of a $(\grstA,\grstB)$-morphism if
$\grstA'=\grstA$ and $\grstB'=\grstB$ and
both $\phi$ and $\psi$ are identities. We shorten this to
$(A,B)$-morphism if both \gr-stacks arise from abelian groups $A$
and $B$ in the usual way. We speak of an equivalence if the
underlying 2-functor $F$ is an equivalence of 2-stacks.

\subsection{Classification I}
\label{sec:class-results}

The classification of 2-$(A,B)$-gerbes follows the usual pattern.
The following theorem generalizes previous results on connective
and hermitian structures on 2-gerbes,
see~\cite{brymcl:deg4:II,bry:quillen} and~\cite{MR2142353}. 
\begin{theorem}
  \label{thm:2}
  Let $\delta\colon A\to B$ be a complex of abelian groups over
  $\cat{C}/X$.  Equivalence classes of 2-$(A,B)$-gerbes are
  classified by the elements of the (ordinary) hypercohomology
  group 
  \begin{equation*}
    \HHH^3(X,A\to B)\,.
  \end{equation*}
\end{theorem}
\begin{proof}
  We only need to sketch the proof, for the details can be lifted
  from the above quoted references and adapted to the present
  situation without difficulty.  Therefore let us only indicate
  how to obtain the cocycle representing the class of a given
  2-$(A,B)$-gerbe.
  
  Let us work in the \cech setting, so let $(U_i\to X)_{i\in I}$
  be a cover as usual. Let $(\tgG,J)$ be a 2-$(A,B)$-gerbe over
  $X$, and let $x_i$, $\phi_{ij}$, and $\alpha_{ijk}$ be objects,
  morphisms, and 2-morphisms providing a full decomposition of
  $\tgG$ relative to the chosen cover as in
  sect.~\ref{sec:classification}.  In addition, let us pick a
  decomposition of the \emph{gerbes} $J(x_i)$ over $U_i$ by
  choosing objects $r_i$ and arrows $\xi_{ij}\colon
  J(\phi_{ij})(r_j) \to r_i$.
  
  Over $U_{ijk}$ we obtain the following diagram in
  $J(x_i)\vert_{U_{ijk}}$:
  \begin{equation*}
      \xymatrix@C+2pc{%
        J(\phi_{ij}) \circ J(\phi_{jk}) (r_k) \ar[d]
        \ar[r]^-{J(\phi_{ij}) (\xi_{jk})}
        & J(\phi_{ij})(r_j) \ar[r]^{\xi_{ij}}
        & r_i \ar[d]^{b_{ijk}} \\
        J(\phi_{ij}\circ \phi_{jk}) (r_k)
        \ar[r]_(.6){J(\alpha_{ijk})(r_k)}
        & J(\phi_{ik}) (r_k) \ar[r]_{\xi_{ik}}
        & r_i
      }
  \end{equation*}
  which defines the section $b_{ijk}\in \SheafAut
  (r_i)\iso B\vert_{U_{ijk}}$. (The left vertical arrow comes
  from the natural transformation built in the definition of
  2-functor~\cite{MR0364245}.)
  
  Pulling back to $U_{ijkl}$ we obtain a cube determined by the
  objects $r_i,\dots, r_l$ whose faces are built from copies of
  the previous diagram. Using
  relation~\eqref{eq:57}, and the fact that $J$ is a
  $\delta$-morphism, we finally have:
  \begin{equation*}
    b^{\phantom{-1}}_{jkl}\,
    b^{-1}_{ikl}\,
    b^{\phantom{-1}}_{ijl}\,
    b^{-1}_{ijk}
    = \delta (a_{ijkl})\,,
  \end{equation*}
  which together with the cocycle relation satisfied by
  $a_{ijkl}$ (consequence of~\eqref{eq:57}), gives the desidered
  cocycle relation for $(a_{ijkl},b_{ijk})$.

  To conclude, let us hint at how the procedure is reversed.  The
  first step is to glue the local trivial 2-gerbes $\gerbes
  (A\vert_{U_i})$ via $a_{ijkl}$. This is standard,
  see~\cite{MR95m:18006,brymcl:deg4:I,brymcl:deg4:II}. Then we
  define a 2-functor $J$ by assigning to each object $x_i$ over
  $U_i$, i.e.\ an $A\vert_{U_i}$-gerbe, the trivial
  $B\vert_{U_i}$-gerbe $J(x_i)=\tors (B\vert_{U_i})$.  Over
  $U_{ijk}$, the section $b_{ijk}$ is used as an automorphism of
  an object $r_i$ of $J(x_i)$, and the cocycle condition above
  ensures compatibility.
\end{proof}

Using the results in sect.~\ref{sec:gerbes-bound-complex} about
$(A,B)$-gerbes we can informally reword the proof of the theorem
by noticing that the representative cocycle of the
2-$(A,B)$-gerbe $\tgG$ was given in terms of $(A,B)$-gerbes.  We
want to make this observation precise.

To this end, let us first observe that if $\delta\colon A\to B$
is a complex of sheaves of abelian groups, then
$\grstG=\tors(A,B)$, introduced in sect.~\ref{sec:b-c-torsors},
is a \gr-stack: the group law is given by the standard contracted
product, so for two pairs $(P,s)$ and $(Q,t)$ we have
$(P,s)\otimes (Q,t)= (P\otimes Q, st)$. In fact $\grstG$ is the
\gr-stack associated to the homomorphism $A\to B$ viewed as an
abelian crossed crossed module. Thus,
\begin{equation*}
  \grstG = \tors (A,B) \iso (A \overset{\delta}{\lto} B)\sptilde\,,
\end{equation*}
cf.\ \cite{MR92m:18019,MR95m:18006}.

The following intermediate results (in the next proposition and
theorem), are also of independent interest, as they provide an
alternative characterization of $(A,B)$-gerbes.  
\begin{proposition}
  \label{prop:4}
  Equivalence classes of $\grstG=\tors(A,B)$-torsors are
  classified by the hypercohomology group
  \begin{math}
    \HHH^2(X,A\to B)\,.
  \end{math}
\end{proposition}
\begin{proof}
  Let $\stP$ be a $\grstG$-torsor.  According to
  sect.~\ref{sec:grstackb-torsors-b} the choice of objects $x_i$
  in the fiber categories $\stP_{U_i}$ with respect to a
  cover $(U_i\to X)_{i\in I}$, determines a pair
  $(g_{ij},\gamma_{ijk})$ with values in $\grstG$ satisfying
  the cocycle identity~\eqref{eq:23}.
   
  Given the specific nature of $\grstG$, each $g_{ij}$ is an
  $(A\vert_{U_{ij}},B\vert_{U_{ij}})$-torsor, namely it
  corresponds to a pair $(P_{ij},t_{ij})$, where $P_{ij}$ is an
  $A$-torsor over $U_{ij}$, and $t_{ij}$ is a section of
  $P_{ij}\cprod{A} B$.  Moreover, $\gamma_{ijk} \colon
  P_{ij}\otimes P_{jk}\isoto P_{ik}$ (suitably restricted to
  $U_{ijk}$), and
  $\delta_\ast(\gamma_{ijk})(t_{ij}t_{jk})=t_{ik}$.
   
  It is perhaps better \emph{not} to assume at this point that
  the torsor $P_{ij}$ is trivialized, but rather consider the
  full blown hypercover $(U^\alpha_{ij},U_i)$, where
  $(U^\alpha_{ij}\to U_{ij})_{\alpha\in A_{ij}}$ is a cover, and
  assume that $s^\alpha_{ij}$ is a trivializing section of
  $P_{ij}$ over $U^\alpha_{ij}$. This choice gives rise to
  sections $a^{\alpha\beta\gamma}_{ijk}$ of
  $A\vert_{U^{\alpha\beta\gamma}_{ijk}}$ and $b^{\alpha}_{ij}$ of
  $B\vert_{U^\alpha_{ij}}$, in the usual way:
  \begin{equation*}
    \gamma_{ijk} (s^\alpha_{ij}\otimes s^\beta_{jk} )
    = s^\gamma_{ik}\, a^{\alpha\beta\gamma}_{ijk}\,,
    \quad
    t_{ij} = (s^\alpha_{ij}\cprod{} 1)\, b^\alpha_{ij}\,.
  \end{equation*}
  Then, using that $s\cdot a\cprod{} 1=s \cprod{} \delta (a) = (s
  \cprod{} 1)\cdot \delta (a)$, it is easily checked that
  $(a^{\alpha\beta\gamma}_{ijk},b^{\alpha}_{ij})$ satisfies the
  cocycle condition with values in the complex $A\to B$ with
  respect to the chosen (hyper)cover.  The rest of details (to
  check that this is well-defined on classes) are routine and
  left to the reader.

  Conversely, given a cocycle with values in $A\to B$ with
  respect to the above hypercover, we can reconstruct
  $(A\vert_{U_{ij}},B\vert_{U_{ij}})$-torsors $(P_{ij},t_{ij})$
  satisfying the cocycle condition. We can then glue the various
  $\grstG_{U_i}$ using this cocycle to obtain a
  $\grstG$-torsor on $X$. Details are again left to the reader.
\end{proof}
Now we consider the trivial 2-gerbe $\tors (\grstG)$ of torsors
over the \gr-stack $\grstG$. Also recall that $\gerbes (A,B)$
denotes the fibered 2-category of $(A,B)$-gerbes over $X$.

Proposition~\ref{prop:4}, and the fact that the same
hypercohomology group classifies $(A,B)$-gerbes as well suggest
the following theorem, which is an extension, in the abelian
context, of~\cite[Proposition~7.3]{MR92m:18019}.  To prepare
the statement, observe that there is an action
\begin{gather*}
  \tors (A) \times \tors (A,B) \lto \tors (A)\\
  \intertext{given on objects by}
  (Q,(P,t)) \lmto (P\otimes Q)\,,
\end{gather*}
where $(P,t)$ is an $(A,B)$-torsor, and $Q$ is an $A$-torsor.  Of
course, since $A$ is an abelian group, $\tors (A)$ is itself a
\gr-stack. Also, by local the triviality of torsors, an
$A$-torsor is locally isomorphic to an $(A,B)$-torsor, thereby
making $\tors (A)$ a $\tors (A,B)$-torsor.
\begin{theorem}
  \label{thm:6}
  Let $\grstG=\tors (A,B)$.  There is an equivalence (of
  2-stacks)
  \begin{equation*}
    F\colon \tors (\grstG) \lisoto \gerbes (A,B)
  \end{equation*}
  given by:
  \begin{equation*}
    F \colon \stP \lmto \tors (A) \cprod{\grstG}\!\! \stP\,.
  \end{equation*}
\end{theorem}
In fact, the equivalence in the proposition is an equivalence of
neutral (or trivial) 2-gerbes bound by $\grstG$.
\begin{proof}
  We will confine ourselves to give a description of the
  2-functor $F$, as well as its quasi-inverse, following \loccit,
  and leave the verification of the details to the reader.

  Given a cover $U\to X$, by definition we have an equivalence
  \begin{equation*}
    \stP_U \lisoto \grstG_U = \tors (A\vert_U,B\vert_U)\,.
  \end{equation*}
  Moreover, observe that for any \gr-stack $\grstG$ and for any
  stack in groupoid with $\grstG$-action $\stP$, we have an
  equivalence
  \begin{equation*}
    \stP \lisoto \stP \cprod{\grstG} \! \grstG
    \qquad
    x \lmto (x,o_\grstG)\,,
  \end{equation*}
  where $o_\grstG$ is the unit object in $\grstG$. By the same
  argument in the proof of~\cite[Proposition~7.3]{MR92m:18019},
  we have the equivalence:
  \begin{equation*}
    \tors (A\vert_U) \lisoto
    \tors (A\vert_U) \cprod{\grstG_U}\! \grstG_U \lisoto
    \tors (A\vert_U) \cprod{\grstG\vert_U}\! \stP_U\,,
  \end{equation*}
  showing that $\tors (A)\cprod{\grstG}\!\stP$ is locally
  equivalent to $\tors (A)$, hence it is an $A$-gerbe. We make it
  into an $(A,B)$-gerbe by defining  
  \begin{equation*}
    \mu \eqdef \delta_\ast \wedge 1 \colon
    \tors (A) \cprod{\grstG}\!\stP \lto \tors (B)\,. 
  \end{equation*}
  This is well-defined, since locally the definition dictates
  \begin{math}
    (Q,(P,t)) \mapsto \delta_\ast (Q)
  \end{math}
  and, using the properties of the contracted product, we have
  \begin{equation*}
    (Q,(P,t)) \lisoto (Q\cdot(P,t), (A,1)) = (P\otimes Q, (A,1))\,,
  \end{equation*}
  so that
  \begin{equation*}
    (Q,(P,t)) \lmto \delta_\ast (P\otimes Q) \iso
    \delta_\ast(P)\otimes \delta_\ast (Q) \iso \delta_\ast (Q)\,,
  \end{equation*}
  since $\delta_\ast (P)\iso B$, by definition of
  $(A,B)$-torsor. (The pair $(A,1)$ represents the unit element
  in $\grstG=\tors (A,B)$.)

  Conversely, let $(\gQ,\mu)$ be an $(A,B)$-gerbe.  Since it is
  in particular an $A$-gerbe, there is an equivalence
  \begin{equation*}
    \gQ\vert_U \iso \tors (A\vert_U)
  \end{equation*}
  with respect to a cover $U\to X$, so that locally the structure
  of $(A,B)$-gerbe becomes
  \begin{equation*}
    \tors (A\vert_U )
    \xrightarrow{\mu\vert_U}
    \tors (B\vert_U)\,.
  \end{equation*}
  In turn this is isomorphic to $\delta{}_\ast$, the ``change of
  structure group'' functor.  To see this, consider the image
  $E=\mu (A)$ of the trivial torsor.  Since $\mu$ commutes with
  the product of torsors (since $Q_1\otimes Q_2\iso Q_1\cprod{A}
  Q_2$ for $A$ abelian), it follows from $Q\iso Q\otimes A$ that
  $E\iso B$, the trivial $B$-torsor.  By local triviality over
  $U$ and the fact that $\mu$ is a $\delta$-morphisms, it follows
  that $\mu (Q)\iso\delta_\ast (Q)$.

  A calculation identical to the one carried out to show that
  $\delta_\ast\wedge 1$ is well-defined, shows that if $(P,t)$ is
  an $(A,B)$-torsor, then the morphism
  \begin{equation*}
    P\otimes - \colon \tors (A\vert_U ) \lto \tors (A\vert_U ) 
  \end{equation*}
  preserves the functor $\delta_\ast$, namely the diagram
  \begin{equation*}
    \xymatrix@C=-0.5pc{%
      \tors (A\vert_U ) \ar[rr]^{P\otimes -}
      \ar@/_1pc/[dr]_(.3){\delta_\ast} & &
      \tors (A\vert_U ) \ar@/^1pc/^(.3){\delta_\ast}[dl]\\
      & \tors (B\vert_U )}
  \end{equation*}
  commutes.  In other words, tensoring with an $(A,B)$-torsor is
  locally a morphism of $(A,B)$-gerbes.  Moreover, since
  \emph{any} equivalence
  \begin{math}
    \nu \colon \tors (A\vert_U ) \to \tors (A\vert_U ) 
  \end{math}
  can be realized as $Q \mapsto P_{\nu}\otimes Q$ for an
  appropriate torsor $P_{\nu}$, compatibility with the previous
  diagram forces $P_{\nu}$ to be an $(A,B)$-torsor.  Denoting by
  $\equ$ the stack of equivalences, the foregoing proves that the
  correspondence
  \begin{equation*}
    \gQ \lmto \equ (\tors (A), \gQ)
  \end{equation*}
  gives the required quasi-inverse equivalence to $F$.
\end{proof}
\begin{remark}
  \label{rem:3}
  The theorem gives another perspective on the canonical morphism
  introduced in sect.~\ref{sec:canonical-f-g}.  Namely, if we
  have a morphism \eqref{eq:20} of Picard \gr-stacks coming from
  the crossed square~\eqref{eq:22}, from the theorem we obtain a
  morphism
  \begin{equation*}
    \gerbes (A,G)  \lto \gerbes (B,H)
  \end{equation*}
  as the conjugate $F_{\grstB}\circ \lambda_\ast \circ
  F^\ast_{\grstA}$ of the induced morphism
  \begin{equation*}
    \lambda_\ast \colon \tors (\grstA) \lto \tors (\grstB)\,,
  \end{equation*}
  where $F_\bullet$ is the appropriate equivalence from
  Theorem~\ref{thm:6} and $F^\ast_\bullet$ its quasi-inverse.  It
  is immediately seen that this morphism corresponds to the
  canonical morphism $(f,u)_\ast$.
\end{remark}
We return to 2-gerbes. The following proposition generalizes
sect.~\ref{sec:classification} and Theorem~\ref{thm:2}, and it
can be considered as the analog of Proposition~\ref{prop:4} to
the case of 2-gerbes.
\begin{proposition}
  \label{prop:5}
  Let $\grstG=\tors (A,B)$. Equivalence classes of
  2-$\grstG$-gerbes are classified by the hypercohomology
  group
  \begin{math}
    \HHH^3(X,A\to B)\,.
  \end{math}
\end{proposition}
\begin{proof}
  Most of the ingredients of the proof can be extracted from the
  cocycle analysis in~\cite{MR95m:18006}, c.f.\ in particular
  \S 4.7.

  Let $\tgG$ be a 2-$\grstG$-gerbe.  Given a cover $(U_i\to
  X)_{i\in I}$, the choice of objects $x_i\in \tgG_{U_i}$
  determines, by analogy with sect.~\ref{sec:grstackb-torsors-b},
  $\grstG$-torsors $\stE_{ij}=\catHom
  (x_j\vert_{U_{ij}},x_i\vert_{U_{ij}})$ over $U_{ij}$.  Note
  that $\stE_{ij}$ is a $\grstG$-torsor, rather than a
  \emph{bi}torsor, thanks to the fact that $\grstG$ is
  braided. The torsors $\stE_{ij}$ satisfy the following cocycle
  condition: we have equivalences
  \begin{subequations}
    \label{eq:24}
    \begin{equation}
      \label{eq:25}
      g_{ijk}\colon
      \stE_{ij} \cprod{\grstG}\!\!
      \stE_{jk} \lisoto
      \stE_{ik}
    \end{equation}
    and natural transformations (isomorphisms):
    \begin{equation}
      \label{eq:26}
      \nu_{ijkl}\colon  g_{ikl} \circ (g_{ijk}\wedge 1)
      \Longrightarrow g_{ijl}\circ (1\wedge g_{jkl})
    \end{equation}
  \end{subequations}
  arising from the pentagonal 2-cell determined by starting at
  \begin{equation*}
    \bigl( \stE_{ij} \cprod{\grstG}\!\!
    \stE_{jk}\bigr) \cprod{\grstG}\!\!
    \stE_{kl}\,,
  \end{equation*}
  and associating with the help of~\eqref{eq:25}. Moreover, the
  morphisms $\nu_{ijkl}$ satisfy the appropriate coherence
  condition extracted from~\eqref{eq:26} over $U_{ijklm}$.

  Notice that a section of, say, $\stE_{ij}$ over
  $U^\alpha_{ij}\to U_{ij}$ is a 1-arrow 
  \begin{math}
    f^\alpha_{ij} \colon
    x_j\vert_{U^\alpha_{ij}} \to
    x_i\vert_{U^\alpha_{ij}}\,,
  \end{math}
  and similarly for the other indices. Therefore the restriction
  $g^{\alpha\beta\gamma}_{ijk}$ of $g_{ijk}$ to
  $U^{\alpha\beta\gamma}_{ijk}$ can be identified with an object
  of $\grstG\vert_{U^{\alpha\beta\gamma}_{ijk}}$. The same
  reasoning leads to the identification of the restriction of
  $\nu_{ijkl}$, with the appropriate decoration of upper indices,
  with an arrow of (a corresponding restriction of) $\grstG$.
  Finally, we note that the equivalence in eqn.~\eqref{eq:25} is
  given by the composition of 1-arrows and 2-arrows in $\tgG$.
  Thus eqns.~\eqref{eq:24} can be interpreted as giving a cocycle
  condition for $(g_{ijk},\nu_{ijkl})$ with values in $\grstG$.
  
  Now, since $\grstG=\tors(A,B)$, is the stack associated to
  the abelian crossed module (i.e.\ complex of abelian groups)
  $\delta\colon A\to B$, the corresponding sheaf of groupoids
  will be
  \begin{equation*}
    \xymatrix{
      A\times B \ar@<0.55ex>[r]^(.6)s \ar@<-0.55ex>[r]_(.6)t & B
    }
  \end{equation*}
  with source and target maps given by $s(a,b)=b$ and
  $t(a,b)=\delta (a)b$, so that (neglecting the upper indices)
  the object $g_{ijk}$ can be identified with a section $b_{ijk}$
  of $B$, and the morphism $\nu_{ijkl}$ with a section $a_{ijkl}$
  of $A$. Now~\eqref{eq:26} reads:
  \begin{equation*}
    \delta (a_{ijkl})\, b_{ijk}\, b_{ikl} = b_{jkl}\, b_{ijl}
  \end{equation*}
  which is the desired relation. Putting it together with the
  cocycle condition for $a_{ijkl}$ determined by the coherence
  condition on the $\nu_{ijkl}$ alluded to above, provides the
  required 3-cocycle with values in the complex $A\to B$.
\end{proof}
Methods similar to the approach of the proof of
Theorem~\ref{thm:6} give the following theorem. We omit the
proof.
\begin{theorem}
  \label{thm:1}
  Let again $\grstG=\tors (A,B)$.  Then a
  2-$\grstG$-gerbe is equivalent to a 2-$(A,B)$-gerbe, where
  the equivalence takes place in the appropriate 3-category.
\end{theorem}
The upshot of the foregoing unfortunately rather lengthy
discussion can be summarized as follows. Given a complex of
abelian groups $\delta : A\to B$, the following two structures on
a 2-$A$-gerbe $\tgG$ are equivalent:
\begin{enumerate}
\item 2-gerbe bound by $\delta : A\to B$, and:\label{item:9}
\item 2-gerbe bound by $\grstG=\tors(A,B)$.\label{item:10}
\end{enumerate}
They correspond to the following crossed squares of the
type~\eqref{eq:22}:
\begin{equation*}
  \text{item \ref{item:9}:}\;
  \vcenter{\xymatrix{%
    A \ar[r]^\delta \ar[d] & B \ar[d] \\
    1 \ar[r] & 1}}
  \qquad
  \text{item \ref{item:10}:}\;
  \vcenter{\xymatrix{%
    A \ar[r] \ar[d]^\delta & 1 \ar[d] \\
    B \ar[r] & 1}}
\end{equation*}
where for case~\ref{item:9} we consider $A$ and $B$ as crossed
modules $A\to 1$ and $B\to 1$, whereas case~\ref{item:10}
corresponds to the crossed module $\lambda \colon \grstG\to
\grstH$ where $\grstH$ is associated to $1\to 1$. The
equivalence can be traced to the symmetry of the crossed square.

Next, we are going to explore the case when the crossed
square~\eqref{eq:22} is non-trivial.

\subsection{Classification II}
\label{sec:2-gerbes-bound-1}

Our first step is to address the case of a 2-gerbe bound by a
crossed module of braided \gr-stacks~\eqref{eq:20} in greater
generality than in the preceding sections. Note that there is an
obvious induced map:
\begin{equation}
  \label{eq:27}
  \lambda_\ast \colon \tors(\grstA) \lto
  \tors (\grstB)\,,
\end{equation}
given by $\stP\to \stP\cprod{\grstA}\!\!
\grstB$. It is convenient to have the following definition
at hand:
\begin{lemma-def}
  \label{lem-def:1}
  Given a cover $\cover{U}_X=(U_i\to X)_{i\in I}$, a
  \emph{1-cocycle with values in~\eqref{eq:27}} is the datum
  of $\grstA$-torsors $\stE_{ij}$ over $U_{ij}$ and
  $\grstB$-torsors $\stF_i$ over $U_i$, such that the
  cocycle condition~\eqref{eq:24} holds for the
  $\stE_{ij}$'s, and moreover there are equivalences of
  $\grstB$-torsors
  \begin{subequations}
    \label{eq:28}
    \begin{equation}
      \label{eq:29}
      \xi_{ij} \colon
      \lambda_\ast (\stE_{ij}) \cprod{\grstB}\!\!
      \stF_j \lisoto \stF_i
    \end{equation}
    and natural transformations (isomorphisms):
    \begin{equation}
      \label{eq:30}
      m_{ijk}\colon \xi_{ij}\circ (1\wedge \xi_{jk})
      \Longrightarrow
      \xi_{ik}\circ (\lambda_\ast (g_{ijk})\wedge 1)\,.
    \end{equation}
  \end{subequations}
  The natural transformations $m_{ijk}$ are subject to the
  following coherence condition:
  \begin{multline}
    \label{eq:31}
    \xi_{il}\ast \lambda_\ast(\nu_{ijkl})\circ
    m_{ijl}\ast (1\wedge \lambda_\ast (g_{jkl})\wedge 1) \circ
    \xi_{ij}\ast m_{jkl}\\
    = m_{ikl}\ast (\lambda_\ast(g_{ijk})\wedge 1\wedge 1)\circ
    m_{ijk}\ast (1\wedge 1 \wedge \xi_{kl})\,.
  \end{multline}
\end{lemma-def}
\begin{remark}
  An easier (but less precise) way of displaying~\eqref{eq:31} is
  to ignore the pastings with the identity 2-arrows, so that we
  have:
  \begin{equation*}
    \lambda_\ast(\nu_{ijkl})\circ m_{ijl} \circ m_{jkl}
    = m_{ikl}\circ m_{ijk}\,.
  \end{equation*}
\end{remark}
\begin{proof}
  The calculations are tedious, but entirely straightforward. We
  will content ourselves to note that one has to form the
  standard cube of morphisms $\xi_{ij}$, etc.\ starting from
  \begin{equation}
    \label{eq:32}
    \lambda_\ast \bigl( \stE_{ij} \cprod{\grstA}\!
    (\stE_{jk} \cprod{\grstA}\!\!
    \stE_{kl})\bigr) \cprod{\grstB}\!\!
    \stF_l
  \end{equation}
  and ending to $\stF_i$, \emph{modulo} the association
  isomorphisms for the contracted product, which have been
  ignored in eq.~\eqref{eq:31}.  Then~\eqref{eq:31} is the result
  of composing the \emph{faces} of this cube. Note that
  in~\eqref{eq:31} there are five terms, since one of the faces
  will be strictly commutative, namely the one corresponding to
  contracting the first two, and the second two terms
  in~\eqref{eq:32}.
\end{proof}
We complement the definition of 1-cocycle with the notion of
equivalence as follows:
\begin{definition}
  \label{def:9}
  Two 1-cocycles $(\stE_{ij},\stF_i)$ and $(\stE'_{ij},\stF'_i)$
  with values in~\eqref{eq:32} are \emph{equivalent} if there
  exist $\grstA\vert_{U_i}$-torsors $\stQ_i$ over $U_i$ such that
  there are equivalences:
  \begin{subequations}
    \label{eq:33}
    \begin{gather}
      \label{eq:34}
      \stE'_{ij}\cprod{\grstA}\!\stQ_j
      \lisoto
      \stQ_i\cprod{\grstA}\! \stE_{ij} \\
      \label{eq:35}
      \lambda_\ast (\stQ_i) \cprod{\grstB}\!
      \stF_i \lisoto \stF'_i\,.
    \end{gather}
  \end{subequations}
\end{definition}
The following is a mild extension of the statement
in~\cite[4.1.11]{MR93k:18019} in the braided case.
\begin{theorem}
  \label{thm:3}
  Equivalence classes of 2-$(\grstA,\grstB)$-gerbes are
  classified by the set
  \begin{equation*}
    \HHH^1(X,\tors (\grstA)\to \tors (\grstB))\,,
  \end{equation*}
  namely the (pointed) set of equivalence classes of 1-cocycles
  in Lemma-Definition \ref{lem-def:1} under the equivalence of
  Definition~\ref{def:9}.
\end{theorem}
\begin{proof}
  Let $\tgG$ be a 2-$(\grstA,\grstB)$-gerbe. Since it
  is in particular a 2-$\grstA$-gerbe, the choice of objects
  $x_i\in \Ob \tgG_{U_i}$ with respect to an open cover
  $\cover{U}_X=(U_i\to X)_{i\in I}$ will generate a 1-cocycle
  $\lbrace \stE_{ij}\rbrace$ with values in $\tors
  (\grstA)$, as in the proof of proposition~\ref{prop:5},
  eqns.~\eqref{eq:24}. This part and the rest of the cocycle
  analysis of the 2-gerbe $\tgG$ is as in~\cite{MR95m:18006},
  especially \S 4.7, with the additional hypothesis that we are
  in the braided case (so that we are in the ``decoupled''
  situation). Full details will be found in \loccit

  The new part is the one related to the extra structure given by
  the 2-functor
  \begin{equation*}
    J\colon \tgG \lto \tors (\grstB)\,,
  \end{equation*}
  as part of the definition of
  2-$(\grstA,\grstB)$-gerbe. Using $J$, for each object
  $x_i$ we obtain a $\grstB$-torsor $\stF_i\eqdef
  J(x_i)$.  Now, recall that $\stE_{ij} = \catHom
  (x_j\vert_{U_{ij}},x_i\vert_{U_{ij}})$. Objects and arrows of
  $\stE_{ij}$ over $U^\alpha_{ij}\to U_{ij}$ correspond to
  1-arrows between $x_j\vert_{U^\alpha_{ij}}$ and
  $x_i\vert_{U^\alpha_{ij}}$ and 2-arrows between them. Via $J$,
  we get equivalences and natural isomorphisms between the
  corresponding torsors $\stF_j$ and $\stF_i$. In
  short, there is an equivalence:
  \begin{equation*}
    \stE_{ij} \lisoto \catHom (\stF_j, \stF_i)\,,
  \end{equation*}
  where the $\catHom$ on the right hand side denotes the category
  of morphisms of torsors (defined e.g. as in \cite[\S
  6]{MR92m:18019}). That is, it is the $\catHom$ in $\tors
  (\grstB)$. In turn, this equivalence can be written in the
  form of eqn.~\eqref{eq:29}, using the correspondence
  \begin{equation*}
    f^\alpha_{ij} \lmto [ y \lmto \lambda_\ast(f^\alpha_{ij})(y) ]
    \;\iso\;
    \lambda_\ast(f^\alpha_{ij}) \wedge y \lmto
    \lambda_\ast(f^\alpha_{ij})(y)\,,
  \end{equation*}
  where $f^\alpha_{ij}$ is an object of $\stE_{ij}$, i.e.\
  1-morphism of $\tgG$, over $U^\alpha_{ij}$, and similarly for
  2-arrows. Here we have also used the fact that $J$ is a
  $\lambda$-morphism, therefore an $\grstA$-torsor
  $\stP$ corresponds to $\lambda_\ast(\stP) =
  \stP\cprod{\grstA}\!\!\grstB$.

  The inverse correspondence is obtained by generalizing the
  standard gluing of local trivial 2-gerbes
  \begin{equation*}
    \tors (\grstA\vert_{U_i})
  \end{equation*}
  in a way analogous to the proof of Thm.~\ref{thm:2}. Namely,
  given a 1-cocycle $(\stE_{ij},\stF_i)$, first we glue $\tors
  (\grstA\vert_{U_j})\vert_{U_{ij}}$ with $\tors
  (\grstA\vert_{U_i})\vert_{U_{ij}}$ via $\stE_{ij}$ by
  \begin{equation*}
    \stP \lmto \stP\cprod{\grstA}\!\! \stE_{ij}\,,
  \end{equation*}
  and verify that this is coherent thanks to
  eqns.~\eqref{eq:24}. Thus we obtain a 2-$\grstA$-gerbe $\tgG$,
  and, as a byproduct, this procedure gives a collection of
  objects $x_i$ providing the labeling with respect to which the
  newly obtained 2-gerbe $\tgG$ is represented by the cocycle
  $\stE_{ij}$. We then define $J$ as:
  \begin{equation*}
    J\vert_{U_i} \colon \tgG_{U_i}\iso \tors (\grstA\vert_{U_i})
    \lto \tors (\grstB\vert_{U_i})
  \end{equation*}
  by assigning to $x_i$ the $\grstB$-torsor $\stF_i$. More
  generally, to any object of $\tgG_{U_i}$, i.e.\ to any
  $\grstA\vert_{U_i}$-torsor $\stP$, we assign the
  $\grstB\vert_{U_i}$-torsor
  \begin{equation*}
    \lambda_\ast (\stP)\cprod{\grstB}\!\! \stF_i\,.
  \end{equation*}
  We leave to the reader the task to verify that the two
  constructions are inverse of one another.

  Finally, given a 2-$(\grstA,\grstB)$-gerbe, a second collection
  of objects $\lbrace y_i\rbrace$ subordinated to the same cover
  determines a new cocycle $(\stE'_{ij},\stF'_i)$. Moreover, for
  each $i\in I$ we have the $\grstA\vert_{U_i}$-torsor
  $\stQ_i=\catHom(x_i,y_i)$. It is easily verified that the
  collection $\lbrace \stQ_i\rbrace$ satisfies both
  eqns.~\eqref{eq:33}.
\end{proof}
When the coefficient complexes of braided stacks come from
complexes of abelian groups the previous theorem can be rephrased
in terms of ordinary hypercohomology. More precisely, we have the
following statement.
\begin{theorem}
  \label{thm:4}
  If the braided \gr-stacks $\grstA$ and $\grstB$ are
  strict and correspond to abelian crossed modules $A\to G$ and
  $B\to H$, respectively, then equivalence classes of
  2-$(\grstA,\grstB)$-gerbes are classified by the
  (ordinary) hypercohomology group
  \begin{equation*}
    \HHH^3(X,A\to B\oplus G\to H)\,,
  \end{equation*}
  namely the coefficient complex is the \emph{cone} (shifted by
  1) of the abelian crossed square~\eqref{eq:22}.
\end{theorem}
\begin{proof}
  We will need to show how to extract an ordinary cocycle with
  value in the cone of \eqref{eq:22} from the abstract cocycle of
  Thm.~\ref{thm:3}.
  
  Let $\grstA=\tors (A,G)$ and $\grstB=\tors(B,H)$ with complexes
  $\delta\colon A\to G$ and $\sigma \colon B\to H$ and
  homomorphisms $f\colon A \to B$ and $u\colon G\to H$ arranged
  to make the square~\eqref{eq:22}. The corresponding (sheaf of)
  crossed module(s) is:
  \begin{equation*}
    \xymatrix@C+1em{%
      A\times G \ar@<0.7ex>[d]^s \ar@<-0.7ex>[d]_t
      \ar[r]^{(f,u)}&
      B\times H \ar@<0.7ex>[d]^s \ar@<-0.7ex>[d]_t \\
      G \ar[r]_u & H
    }
  \end{equation*}
  where in both cases the source and target maps $s$ and $t$ are
  as in the proof of Proposition~\ref{prop:5},
  page~\pageref{eq:24}. Thus the additive functor $\lambda \colon
  \grstA \to \grstB$ is induced (after having taken the associate
  stack functor) by the pair $(f,u)$.

  After having gone through these recollections, let us consider
  a 2-$(\grstA,\grstB)$-gerbe $\tgG$, and let us once again
  choose a cover $\cover{U}_X=(U_i\to X)$, and objects $x_i\in
  \ob \tgG_{U_i}$.  By Theorem~\ref{thm:3}, we obtain a 1-cocycle
  $(\stE_{ij},\stF_i)$ with values in the complex~\eqref{eq:27}
  satisfying eqns.~\eqref{eq:24} and~\eqref{eq:28}.  Our first
  task is to complement the proof of Proposition~\ref{prop:5},
  and obtain a 1-cocycle with values in the complex
  \begin{math}
    \lambda\colon \grstA \lto \grstB
  \end{math}
  itself.
  
  To this end, we will need to decompose the torsors $\stE_{ij}$
  as well as $\stF_i$ with respect to some choice of objects, and
  then apply the reasoning preceding eq.~\eqref{eq:23}.  More
  precisely, consider objects $f^\alpha_{ij}$ and
  ${{y_i}^\alpha}_j$, of $\stE_{ij}$ and $\stF_i$, respectively,
  given $(U^\alpha_{ij}\to U_{ij})_{\alpha \in
    A^\alpha_{ij}}$. (Similarly, we denote by ${{y_j}^\alpha}_i$
  an object of $\stF_j$ over $U^\alpha_{ij}$.) Then, since
  $\stF_i$ is a $\grstB$-torsor, the morphism $\xi_{ij}$ in
  eq.~\eqref{eq:29} translates into
  \begin{equation}
    \label{eq:36}
    (f^\alpha_{ij})_\ast ({{y_j}^\alpha}_i)\iso
    {{y_i}^\alpha}_j\cdot h^\alpha_{ij}\,,
  \end{equation}
  where $h^\alpha_{ij}$ is an object of $\grstB$ over
  $U^\alpha_{ij}$. (Here we have used the notation
  $(f^\alpha_{ij})_\ast=J(f^\alpha_{ij})$.)  Moreover,
  ${{y_i}^\alpha}_j$ and ${{y_i}^\beta}_k$ are related by:
  \begin{equation}
    \label{eq:37}
    {{y_i}^\alpha}_j \iso {{y_i}^\beta}_k \cdot
    q^{\beta\alpha}_{kij}\,,
  \end{equation}
  with $q^{\beta\alpha}_{kij}$ an object of $\grstB$ over
  $U^{\alpha\beta}_{ijk}$. It easily seen that these new objects
  satisfy the identity (up to isomorphism):
  \begin{equation}
    \label{eq:38}
    q^{\beta\alpha}_{kij} \cdot q^{\alpha\gamma}_{jil} \iso
    q^{\beta\gamma}_{kil}\,.
  \end{equation}
  For the part of the cocycle involving the $\stE_{ij}$'s alone,
  subject to eqns.~\eqref{eq:24}, our choice of objects
  determines an object $g^{\alpha\beta\gamma}_{ijk}$ of $\grstA$
  obtained from eqn.~\eqref{eq:25} in the standard way:
  \begin{equation*}
    f^\alpha_{ij}\wedge f^\beta_{jk} \lmto
    f^\alpha_{ij}\circ f^\beta_{jk}
    \iso g^{\alpha\beta\gamma}_{ijk} \circ f^\gamma_{ik}\,.
  \end{equation*}
  (Recall that the map $g_{ijk}$ is just composition of 1-arrows
  of $\tgG$.) Moreover, still using the arguments
  in~\cite{MR95m:18006}, starting from eqn.~\eqref{eq:26} we
  arrive at the morphism in $\grstA$:
  \begin{subequations}\label{eq:42}
  \begin{equation}
    \label{eq:39}
    \nu^{\alpha\beta\delta\gamma\eta\epsilon}_{ijkl} \colon
    g^{\alpha\beta\gamma}_{ijk}\cdot
    g^{\gamma\delta\epsilon}_{ikl}
    \lisoto
    g^{\beta\delta\eta}_{jkl}\cdot g^{\alpha\eta\epsilon}_{ijl}\,.
  \end{equation}
  To translate eqn.~\eqref{eq:30}, compute the composition over
  $U^{\alpha\beta\gamma}_{ijk}$:
  \begin{equation*}
    (f^\alpha_{ij} \circ f^\beta_{jk})_\ast ({{y_k}^\beta}_j)
  \end{equation*}
  in the two possible ways. A standard calculation, where we
  use~\eqref{eq:36} and~\eqref{eq:37}, yields the sought-after
  arrow in $\grstB$:
  \begin{equation}
    \label{eq:40}
    m^{\alpha\beta\gamma}_{ijk} \colon
    h^\alpha_{ij}\, q^{\alpha\beta}_{ijk} \, h^\beta_{jk} \lisoto
    \lambda (g^{\alpha\beta\gamma}_{ijk})\,
    q^{\alpha\gamma}_{jik} \, h^\gamma_{ik} \, q^{\gamma\beta}_{ikj}\,.
  \end{equation}
  This arrow in turn satisfies a cocycle condition, which is the
  translation of eqn.~\eqref{eq:31}. We arrive at it by
  considering the expression
  \begin{equation*}
    h^\alpha_{ij}\, q^{\alpha\beta}_{ijk} \, h^\beta_{jk} \,
    q^{\beta\delta}_{jkl} \, h^\delta_{kl}\,,
  \end{equation*}
  which would correspond to $\xi_{ij}\circ (1\wedge\xi_{jk})
  \circ (1\wedge 1 \wedge \xi_{kl})$, and computing it in the two
  possible obvious ways using~\eqref{eq:40}, the braiding of
  $\grstB$---and the help of~\eqref{eq:38}. The calculation
  itself proceeds according to the techniques expounded
  in~\cite{MR95m:18006}, therefore we will not reproduce it
  here. The result is that the arrows
  $m^{\alpha\beta\gamma}_{ijk}$ satisfy the cocycle condition:
  \begin{equation}
    \label{eq:41}
    \lambda (\nu^{\alpha\beta\delta\gamma\eta\epsilon}_{ijkl})
    \circ
    m^{\gamma\delta\epsilon}_{ikl}
    \circ
    m^{\alpha\beta\gamma}_{ijk} =
    m^{\alpha\eta\epsilon}_{ijl} \circ
    m^{\beta\delta\eta}_{jkl}\,.
  \end{equation}
  \end{subequations}
  Of course this identity holds \emph{modulo} the obvious
  isomorphisms arising from the association and braiding functors
  in $\grstB$, which we have silently ignored, as well as the
  pull-back functors between different fiber categories.

  The cocycle with values in the complex $\lambda \colon \grstA
  \to \grstB$ we have obtained comprises the quintuple:
  \begin{equation}
    \label{eq:47}
    \bigl( h^\alpha_{ij}\,, q^{\alpha\beta}_{ijk}\,,
    m^{\alpha\beta\gamma}_{ijk}\,,
    g^{\alpha\beta\gamma}_{ijk}\,,
    \nu^{\alpha\beta\delta\gamma\eta\epsilon}_{ijkl}
    \bigr)
  \end{equation}
  subject to eqns.~\eqref{eq:42} plus the cocycle condition on
  the terms $\nu^{\alpha\beta\delta\gamma\eta\epsilon}_{ijkl}$
  arising from the coherence condition on the
  maps~\eqref{eq:26}. We refrain from displaying such condition
  here.
  
  Now let us use the fact that both the \gr-stacks $\grstA$ and
  $\grstB$ are strict and in fact associated to crossed modules.
  From the recollections at the beginning we have that in the
  above quintuple $g^{\alpha\beta\gamma}_{ijk}$ will be a section
  of the abelian group sheaf $G$, $h^\alpha_{ij}$ and
  $q^{\alpha\beta}_{ijk}$ are both sections of $H$, whereas the
  arrows $m^{\alpha\beta\gamma}_{ijk}$ and
  $\nu^{\alpha\beta\delta\gamma\eta\epsilon}_{ijkl}$ will
  correspond to sections of $B$ and $A$, respectively denoted
  $b^{\alpha\beta\gamma}_{ijk}$ and
  $a^{\alpha\beta\delta\gamma\eta\epsilon}_{ijkl}$, satisfying
  the (strict) identities:
  \begin{subequations}
    \label{eq:43}
    \begin{align}
      \delta (a^{\alpha\beta\delta\gamma\eta\epsilon}_{ijkl})
      \cdot g^{\alpha\beta\gamma}_{ijk}\cdot
      g^{\gamma\delta\epsilon}_{ikl} & =
      g^{\beta\delta\eta}_{jkl}\cdot
      g^{\alpha\eta\epsilon}_{ijl}\,, \label{eq:44}\\
      \sigma (b^{\alpha\beta\gamma}_{ijk})\cdot
      h^\alpha_{ij}\cdot q^{\alpha\beta}_{ijk} \cdot h^\beta_{jk}
      & =
      u (g^{\alpha\beta\gamma}_{ijk})\cdot
      q^{\alpha\gamma}_{jik} \cdot h^\gamma_{ik} \cdot
      q^{\gamma\beta}_{ikj}\,, \label{eq:45} \\
      f (a^{\alpha\beta\delta\gamma\eta\epsilon}_{ijkl})
      \cdot
      b^{\gamma\delta\epsilon}_{ikl}
      \cdot
      b^{\alpha\beta\gamma}_{ijk} & =
      b^{\alpha\eta\epsilon}_{ijl} \cdot
      b^{\beta\delta\eta}_{jkl}\,.\label{eq:46}
    \end{align}
  \end{subequations}
  To these equations we have to add the condition satisfied by
  the $a^{\alpha\beta\delta\gamma\eta\epsilon}_{ijkl}$ as
  consequence of the identity satisfied by the arrows
  $\nu^{\alpha\beta\delta\gamma\eta\epsilon}_{ijkl}$.

  It is just a matter of using the definition of the mapping cone
  of a complex to realize that~\eqref{eq:43}  express the
  condition for the quintuple
  \begin{equation}
    \label{eq:48}
    \bigl( h^\alpha_{ij}\,, q^{\alpha\beta}_{ijk}\,,
    b^{\alpha\beta\gamma}_{ijk}\,,
    g^{\alpha\beta\gamma}_{ijk}\,,
    a^{\alpha\beta\delta\gamma\eta\epsilon}_{ijkl}
    \bigr)
  \end{equation}
  to define a cocycle of degree 3 with values in the complex
  \begin{equation}
    \label{eq:49}
    A \xrightarrow{(f,\delta)} B\oplus G \xrightarrow{\sigma\cdot
      u^{-1}} H\,,
  \end{equation}
  with $A$ placed in degree 0.  This finishes the proof.
\end{proof}
\begin{remark}
  Ignoring the intimidating upper indices relative to the
  hypercover used in the proof allows to set
  $q^{\alpha\beta}_{ijk}=1$ so that eqns.~\eqref{eq:43}, plus the
  cocycle identity on $a_{ijkl}$, will assume the standard form
  for a \cech cocycle of degree 3 with values in~\eqref{eq:49}.
\end{remark}
\begin{remark}
  The \emph{proof} of Theorem~\ref{thm:4} actually gives slightly
  more, in that it gives the 3-cocycle with values in the complex
  $\lambda \colon \grstA \to \grstB$ corresponding to the torsor
  1-cocycle with values in~\eqref{eq:27}, regardless of whether
  the involved (braided) \gr-stacks are associated to crossed
  modules.
\end{remark}
\begin{remark}
  The statement (but not the proof) of Theorem~\ref{thm:4}
  subsumes those of Theorem~\ref{thm:2} and
  Proposition~\ref{prop:5}.
\end{remark}
\begin{remark}
  \label{rem:2}
  The cocycle identities~\eqref{eq:43} satisfied by the
  quintuple~\eqref{eq:48} are symmetric under the exchange
  \begin{equation*}
    b^{\alpha\beta\gamma}_{ijk}
    \longleftrightarrow
    g^{\alpha\beta\gamma}_{ijk}\,,
  \end{equation*}
  and the corresponding exchanges $f\leftrightarrow\delta$ and
  $\sigma \leftrightarrow u$.  This symmetry rests upon that of
  the crossed square~\eqref{eq:22} determined by the crossed
  module of strict \gr-categories under consideration. Thus,
  calling $\mathcal{P}$ the crossed square~\eqref{eq:22}, a
  2-gerbe $\tgG$ satisfying the hypotheses of Theorem~\ref{thm:4}
  ought be more properly called a 2-$\mathcal{P}$-gerbe.
\end{remark}
Let us also observe that the situation described by the
hypotheses of Theorem~\ref{thm:4} has another interesting
subcase. Namely, we can consider a complex of length 3 as it was
done is sect.~\ref{sec:complexes-length-3}, and then define the
notion of a 2-gerbe bound by this complex. This is clearly
possible using Theorem~\ref{thm:4} by setting $G=1$ (or $B=1$).
Thus we can state the following definition, generalizing
Definition~\ref{def:4}.
\begin{definition}
  \label{def:10}
  Let $A\overset{\delta}{\lto} B\overset{\sigma}{\lto} C$ be a
  complex of (sheaves of) abelian groups on $\cat{C}/X$.  A
  2-$(A,B,C)$-gerbe is a 2-$A$-gerbe $\tgG$ equipped with a
  structure of 2-$(\grstA,\grstB)$-gerbe where $\grstA=\tors (A)$
  and $\grstB=\tors(B,C)$.
\end{definition}
In the previous definition $\grstA$ is the \gr-stack associated
to the abelian group $A$ viewed as a crossed module $A\to 1$. The
additive functor $\lambda$ is thus determined by the pair
$(\delta, 1)$.  Of course, up to a trivial isomorphism on the
resulting cohomology group, we could have chosen the combination
$\grstA=\tors(A,B)$, $\grstB=\tors(1,C)$ due to the symmetry of
the two resulting crossed squares.

In the end, one outcome of the material expounded in this section
is that the theory of 2-$(\grstA,\grstB)$-gerbes can account for
2-gerbes bound by complexes of abelian \emph{groups} which are in
fact of length 3.  It is particularly relevant, as we will see in
the applications to Hermitian Deligne cohomology further below,
that hypercohomology groups with values in the cone of a square
can naturally be obtained in this framework.

Two issues however suggest to push this circle of ideas a little
further.  On one hand, it is natural to ask whether
Definition~\ref{def:10} admits a ``naive'' generalization by
simply replacing groups with \gr-stacks. On the other hand,
capturing the geometric meaning of the
hypercohomology groups with values in the complex~\eqref{eq:11}
requires that we have a theory of 2-gerbes bound by complexes of
the appropriate length, which cannot be obtained from what we
have right now.

We will address the issue in section~\ref{sec:2-gerbes-bound-2}.

\subsection{Examples}
\label{sec:examples-1}
We review here a few fairly standard examples to illustrate the
foregoing theory.  In fact, the following examples are the
2-gerbe counterpart of the examples presented in
sect.~\ref{sec:examples} and~\ref{sec:further-examples}.  The
analysis of more interesting examples will be deferred until the
last section dedicated to the interpretation of certain Deligne
cohomology groups.

$X$ is an algebraic manifold, and we work with the standard site
determined by $X\an$ (see above).

\subsubsection{Connective structures (or ``concept of connectivity'').}
\label{sec:conn-struct-or}

This is the classical example due to Brylinski and McLaughlin
(see~\cite{brymcl:deg4:I,brymcl:deg4:II} and~\cite{bry:quillen}).

Let $\tgG$ be a 2-gerbe over $X$.  As expected, a connective
structure (or ``concept of connectivity'' as it was originally
called) on $\tgG$ is a structure of 2-gerbe bound by the complex
\begin{equation*}
  \sho{X}^\unit \xrightarrow{\d\log}  \shomega[1]{X}
\end{equation*}
in the sense of Definition~\ref{def:7} and Lemma~\ref{lem:2}.
Thus we retrieve Brylinski and McLaughlin's original definition,
wherein the connective structure is seen as a 2-functor assigning
to each local object of $\tgG$ over $U$ a corresponding $
\shomega[1]{U}$-gerbe.  In light of Proposition~\ref{prop:5} and
Theorem~\ref{thm:1} $\tgG$ can just as well be considered as a
2-gerbe bound by the \gr-stack of
$(\sho{X}^\unit,\shomega[1]{X})$-torsors.

From the classification results (see \loccit for the original
arguments) we have that 2-gerbes with this connective structure
are classified by the hypercohomology group:
\begin{equation*}
  \HHH^3(X,\sho{X}^\unit \xrightarrow{\d\log}  \shomega[1]{X})
  \iso
  \delH[4]{X}{\ZZ}{2}\,.
\end{equation*}

\subsubsection{Hermitian structures.}
\label{sec:hermitian-structures-1}

This version of the idea of hermitian structure was introduced
in~\cite{MR2142353} by analogy with the notion of connective
structure in the above mentioned works by Brylinski and
McLaughlin. Thus, a 2-$\sho{X}^\unit$-gerbe $\tgG$ over $X$ with
hermitian structure is a 2-gerbe bound by the complex:
\begin{equation*}
  \sho{X}^\unit \xrightarrow{\abs{\cdot}^2} \she[+]{X}\,,
\end{equation*}
or, alternatively, by the \gr-stack of $(\sho{X}^\unit,
\she[+]{X})$-torsors.  Equivalence classes of such 2-gerbes are
classified by the Hermitian Deligne cohomology group of weight 1:
\begin{equation*}
  \HHH^3(X,\sho{X}^\unit \xrightarrow{\abs{\cdot}^2} \she[+]{X})
  \iso 
  \dhhH[4]{X}{1}\,,
\end{equation*}
where we use the same quasi-isomorphism as in
sect.~\ref{sec:hermitian-structures}.

It is easy to continue the list of examples by promoting those of
sect.~\ref{sec:further-examples} to the realm of 2-gerbes.  We
will not do so here, and leave this task to the interested
reader.  We will examine finer examples of geometric structures
on 2-gerbes in sect.~\ref{sec:geom-appl}.

\section{2-Gerbes bound by complexes of higher degree}
\label{sec:2-gerbes-bound-2}

So far, we have outlined a theory of 2-gerbes bound (in the
appropriate sense) by a two-step complex of braided \gr-stacks.
We have found that this theory is powerful enough to provide an
interpretation in geometric terms of the elements of degree three
hypercohomology groups with values in (cones of) crossed squares
of abelian groups.  However, as pointed out above, we need to
address the case where the coefficient complexes have degree
higher than 3, where the degree loosely corresponds to the
length.  We set out to accomplish this goal by generalizing the
concept of $(A,B,C)$-gerbe, introduced in sect.~\ref{sec:a-b-c},
to the case of 2-gerbes by promoting the coefficient groups to be
\gr-stacks instead.  We will ultimately be interested in the case
of \gr-stacks associated to abelian crossed modules, therefore
the general style for this section will be slightly more
descriptive---and perhaps informal---compared to the preceding
ones.

\subsection{$(\grstB,\grstC)$-torsors}
\label{sec:grstb-grstc-torsors}

Consider a complex (i.e.\ a morphism) of two (braided, as usual)
\gr-stacks $\mu \colon \grstB \lto \grstC$ on $\cat{C}/X$. By
analogy with sect.~\ref{sec:b-c-torsors}, define a
$(\grstB,\grstC)$-torsor to be a pair $(\stP,\sigma)$, where
$\stP$ is a $\grstB$-torsor, and $\sigma$ is an equivalence:
\begin{equation*}
  \sigma \colon \stP\cprod{\grstB} \grstC \lisoto
  \grstC 
\end{equation*}
where on the right-hand side $\grstC$ is considered as a trivial
torsor.  Equivalently, we require that there be a morphism:
\begin{equation*}
  \sigma \colon \stP \lto \grstC\,,
\end{equation*}
namely a global object (over $\cat{C}/X$) of the fibered category
$\catHom(\stP,\grstC)$. Yet another equivalent point of view is
to regard $\sigma$ as a global object of the torsor
$\stP\cprod{\grstB} \grstC$. The latter point of view is useful
to arrive at a description in terms of cocycles. Suppose indeed
that $\stP$ is decomposed as in
sect.~\ref{sec:grstackb-torsors-b}, with associated 1-cocycle
$(b_{ij},\beta_{ijk})$ with values in $\grstB$
satisfying~\eqref{eq:23}. By the stack condition, an object of
$\stP \cprod{\grstB}\! \grstC$ is equivalent to a collection of
pairs
\begin{equation*}
  (x_i,c_i)\in \ob \bigl. ( \stP \cprod{\grstB}\! \grstC)
  \bigr\rvert_{U_i}
\end{equation*}
satisfying the descent condition on objects. Using the
description of contracted product found
in~\cite[\S 6.7]{MR92m:18019}, we find that the objects $c_i\in
\ob \grstC\vert_{U_i}$ satisfy the condition
\begin{subequations}
  \label{eq:51}
\begin{equation}
  \label{eq:50}
  \rho_{ij} \colon c_j \lisoto \mu(b_{ij}^\ast)\cdot c_i
\end{equation}
(where $b^\ast$ is a quasi-inverse of $b$).  This essentially
follows from the fact that a morphism
\begin{math}
  (x_j,c_j)\vert_{U_{ij}}\to
  (x_i,c_i)\vert_{U_{ij}}
\end{math}
in $\stP \cprod{\grstB}\! \grstC$ corresponds to the triple
\begin{equation*}
  \bigl(
  x_j\cdot b_{ji} \lisoto x_i\,,
  b_{ji}\,,
  c_j \lisoto \mu (b_{ji})\cdot c_i \bigr)
\end{equation*}
modulo an equivalence explained in \loccit The $\rho_{ij}$ are
morphisms in $\grstC_{U_{ij}}$ which then satisfy the coherence
condition:
\begin{equation}
  \label{eq:52}
  \mu (\beta_{ijk})\circ \rho_{ij}\circ \rho_{jk} = \rho_{ik}\,.
\end{equation}
\end{subequations}
This and~\eqref{eq:23} ensure, via the above mentioned
equivalence relation, that $b_{ki}$ and $b_{kj}\cdot b_{ji}$
correspond the same morphism, thereby ensuring that the cocycle
condition in the descent condition is indeed satisfied.
\begin{definition}
  \label{def:11}
  The triple $(b_{ij},\beta_{ijk},\rho_{ij})$ satisfying
  equations~\eqref{eq:51}, plus~\eqref{eq:23} and the coherence
  condition on the $\beta_{ijk}$ is a \emph{1-cocycle with values
    in the complex} $\mu\colon \grstB \lto \grstC$.
\end{definition}
Given the square of \gr-stacks
\begin{equation*}
  \xymatrix@+1pc{
    \grstB \ar[r]^\mu \ar[d]_\psi &
    \grstC \ar[d]^\pi_{}="pi" \\
    \grstB'\ar[r]_{\mu'}^{}="mupr"
    \ar@{=>}@/_0.3pc/ _k "pi";"mupr" &
    \grstC'}
\end{equation*}
we obtain a morphism
\begin{equation}
  \label{eq:56}
  (\psi,\pi)_\ast \colon \tors (\grstB,\grstC) \lto
  \tors (\grstB',\grstC')
\end{equation}
by sending a $\grstB$-torsor $\stP$ to $\stP\cprod{\grstB}
\grstB'$ and the morphism $\sigma$ to $\pi\circ\sigma$.

A morphism from a $(\grstB,\grstC)$-torsor $(\stP,\sigma)$ to a
$(\grstB',\grstC')$-torsor $(\stP',\sigma')$ consists of a square
\begin{equation}
  \label{eq:58}
  \vcenter{\xymatrix@+1pc{
      \stP \ar[r]^\sigma \ar[d]_\xi &
      \grstC \ar[d]^\pi_{}="pi" \\
      \stP'\ar[r]_{\sigma'}^{}="sigmapr"
      \ar@{=>}@/_0.3pc/ _t "pi";"sigmapr" &
      \grstC'}}
\end{equation}
In particular, for $\grstB'=\grstB$, $\grstC'=\grstC$, it reduces
to a triangle
\begin{equation}
  \label{eq:60}
  \vcenter{\xymatrix@+1pc{
    \stP \ar@/^/[dr]^\sigma_{}="pi" \ar[d]_\xi \\
    \stP'\ar[r]_{\sigma'}^{}="sigmapr"
    \ar@{=>}@/_0.1pc/ _t "pi";"sigmapr" &
    \grstC}}  
\end{equation}
Actually, any morphism~\eqref{eq:58} can be factored as the
canonical morphism~\eqref{eq:56} followed by a morphism of
$(\grstB',\grstC')$-torsors.  A morphism will be called an
equivalence if so is the underlying functor $\xi$.

In summary, a $(\grstB,\grstC)$-torsor $\stP$ determines (and it
is determined by, up to equivalence) an equivalence class of
1-cocycles as in the definition.  The equivalence relation being
the obvious one, we obtain the following
\begin{proposition}\hfill\par
  \label{prop:3}
  \begin{enumerate}
  \item\label{item:11} Equivalence classes of
    $(\grstB,\grstC)$-torsors are classified by the cohomology
    set:
    \begin{equation*}
      \HHH^1(X,\grstB \lto \grstC)\,.
    \end{equation*}
  \item\label{item:12} Moreover, if $\mu\colon \grstB\to \grstC$
    comes from the crossed square of abelian groups:
    \begin{equation*}
      \xymatrix{%
        B \ar[d]_\sigma \ar[r]^g & C \ar[d]^\tau \\
        H \ar[r]_v & K
      }
    \end{equation*}
    then the above cohomology set can be identified with the
    hypercohomology group
    \begin{equation*}
      \HHH^2(X,B\to C\oplus H\to K)\,.
    \end{equation*}
  \end{enumerate}
\end{proposition}
\begin{proof}
  Repeats previous arguments, hence omitted.
\end{proof}
\begin{remark}
  We can use the statement in the above proposition to obtain
  another characterization of gerbes bound by length~3-complex,
  specifically, the cone of the above crossed square. This gives
  an alternative point of view for the discussion in
  sect.~\ref{sec:a-b-c}.
\end{remark}
Since by definition $(\grstB,\grstC)$-torsors are
$\grstB$-torsors which become trivial as $\grstC$-torsors, the
following alternative characterization of
$(\grstB,\grstC)$-torsors coming from a crossed square as in
Proposition~\ref{prop:3}--\ref{item:12} is an immediate
consequence of Theorem~\ref{thm:6}:
\begin{proposition}
  \label{prop:6}
  Let $\mu\colon \grstB\to \grstC$ arise from a crossed square as
  in Proposition~\ref{prop:3}--\ref{item:12}.  The 2-functor $F$
  of Theorem~\ref{thm:6} induces an equivalence
  \begin{equation*}
    \tors (\grstB,\grstC) \lisoto
    \gerbes (B,H)_{(\tors (C),\tau_\ast)} 
  \end{equation*}
  where the right hand side denotes the ``fiber'' of the
  canonical morphism
  \begin{equation*}
    (g,v)_\ast \colon \gerbes (B,H)\to \gerbes (C,K)
  \end{equation*}
  over the neutral $(C,K)$-gerbe, that is $\tau_\ast \colon \tors
  (C)\to \tors (K)$.
\end{proposition}
\begin{proof}
  If $\stP$ is a $(\grstB,\grstC)$-torsor, by definition there is
  a morphism $\sigma \colon \stP \to \grstC$, and the diagram
  \begin{equation*}
    \xymatrix@C+1pc{%
      \tors (\grstB) \ar[r]^(0.4){F_\grstB} \ar[d]_{\mu_\ast} & 
      \gerbes (B,H) \ar[d]^{(g,v)_\ast} \\
      \tors (\grstC) \ar[r]_(0.4){F_\grstC} & \gerbes (C,K)}
  \end{equation*}
  from remark~\ref{rem:3} gives
  \begin{equation*}
    \xymatrix{%
      \stP \ar@{|->}[r] \ar[d]_\sigma &
      \tors(B)\cprod{\grstB}\! \stP
      \ar[d]^{g_\ast\wedge \sigma} \\
      \grstC \ar@{|->}[r] & \tors (C)}
  \end{equation*}
  and the lower right corner gives the neutral $(C,K)$-gerbe.
\end{proof}

\subsection{Complexes of braided \gr-stacks}
\label{sec:complexes-braided-gr}

Let $\grstA$, $\grstB$, and $\grstC$ be braided \gr-stacks over
$\cat{C}/X$, and let $\lambda \colon \grstA \to \grstB$ and $\mu
\colon \grstB \to \grstC$ be additive functors. We define the
composition
\begin{equation}
  \label{eq:53}
  \grstA \overset{\lambda}{\lto} \grstB
  \overset{\mu}{\lto} \grstC
\end{equation}
a \emph{complex of \gr-stacks} if $\mu\circ\lambda$ is isomorphic
to the ``null'' functor
\begin{math}
  \grstA\lto \mathbf{1}\,,
\end{math}
to the punctual category determined by the unit object $o_\grstC$
of $\grstC$.

As before, a situation of particular interest for us will be when
everything in sight is strict, and all the \gr-stacks above are
in fact associated to abelian crossed modules.  Building on what
we have already seen in sect.~\ref{sec:crossed-modules-gr},
assume that the morphisms $\lambda$ and $\mu$ are associated to
the squares
\begin{equation*}
  \xymatrix{%
    A \ar[d]_\delta \ar[r]^f & B \ar[d]^\sigma \\
    G \ar[r]_u & H
  }\qquad
  \xymatrix{%
    B \ar[d]_\sigma \ar[r]^g & C \ar[d]^\tau \\
    H \ar[r]_v & K
  }
\end{equation*}
respectively, which we splice together to obtain the map of
complexes:
\begin{equation}
  \label{eq:55}
  \begin{matrix}
  \xymatrix{%
    A \ar[d]_\delta \ar[r]^f & B \ar[d]^\sigma \ar[r]^g & C
    \ar[d]^\tau \\ 
    G \ar[r]_u & H \ar[r]_v & K}
  \end{matrix}
\end{equation}
In all the above we have of course assumed $\grstC$ to be
associated to the complex $\tau\colon C\to K$, the rest of the
notations being as in sect.~\ref{sec:crossed-modules-gr}.

\subsection{2-$(\grstA,\grstB,\grstC)$-gerbes}
\label{sec:2-grsta-grstb}

The main idea is to define 2-gerbes bound by the
complex~\eqref{eq:53} of braided \gr-stacks by analogy with what
was done for gerbes in sect.~\ref{sec:a-b-c}.
\begin{definition}
  \label{def:12}
  Let $\tgG$ be a 2-gerbe over $\cat{C}/X$.  We say that $\tgG$
  is bound by the complex~\eqref{eq:53}, or that is is a
  2-$(\grstA,\grstB,\grstC)$-gerbe, for short, if there is a
  2-functor 
  \begin{equation*}
    \Tilde J\colon \tgG \lto \tors (\grstB, \grstC)
  \end{equation*}
  such that $\tgG$ is a 2-$(\grstA,\grstB)$-gerbe for the
  $\lambda$-morphism defined by the composition of $\Tilde J$
  with the obvious morphism $\tors (\grstB, \grstC) \to \tors
  (\grstB)$.
\end{definition}
Next, we can consider the diagram of \gr-stacks:
\begin{equation*}
  \xymatrix@+1pc{
    \grstA \ar[r]^\lambda \ar[d]_\phi &
    \grstB \ar[r]^\mu \ar[d]^\psi_{}="psi" &
    \grstC \ar[d]^\pi_{}="pi" \\
    \grstA'\ar[r]_{\lambda'}^{}="lambdapr"
    \ar@{=>}@/_0.3pc/ _\jmath "psi";"lambdapr" &
    \grstB'\ar[r]_{\mu'}^{}="mupr"
    \ar@{=>}@/_0.3pc/ _k "pi";"mupr" &
    \grstC'
  }  
\end{equation*}
where the top and bottom rows are complexes in the sense
specified above in sect.~\ref{sec:complexes-braided-gr}.  Still
by analogy with sect.~\ref{sec:a-b-c}, where the corresponding
concept for gerbes was introduced, we define a morphism of a
2-$(\grstA,\grstB,\grstC)$-gerbe $\tgG$ to a
2-$(\grstA',\grstB',\grstC')$-gerbe $\tgG'$ to be a cartesian
2-functor
\begin{equation*}
  F \colon \tgG \lto \tgG' 
\end{equation*}
which is a $\phi$-morphism, supplemented by a 2-natural
transformation
\begin{equation*}
  \Tilde \alpha \colon (\psi,\pi)_\ast\circ \Tilde J
  \Longrightarrow \Tilde J'\circ F \colon
  \tgG \lto \tors (\grstB, \grstC)\,.
\end{equation*}
We require that composing (pasting) this with the obvious
morphism $\tors (\grstB, \grstC) \to \tors (\grstB)$ gives (up to
a modification) the natural morphism associated to the underlying
$(\phi,\psi)$-morphism.

\subsection{Classification III}
\label{sec:classification-iii}

Given the complex~\eqref{eq:53}, we obtain a corresponding
``complex'' of trivial 2-gerbes:
\begin{equation}
  \label{eq:54}
  \tors (\grstA) \overset{\lambda_\ast}{\lto}
  \tors (\grstB) \overset{\mu_\ast}{\lto}
  \tors (\grstC)
\end{equation}
where $\mu_\ast\circ \lambda_\ast\iso (\mu\circ\lambda)_\ast\iso
\mathbf{1}$.  
\begin{lemma-def}
  \label{lem-def:2}
  Given a cover $\cover{U}_X=(U_i\to X)_{i\in I}$, a
  \emph{1-cocycle with values in~\eqref{eq:54}} is given by
  the same data as those for a 1-cocycle with values
  in~\eqref{eq:27} stated in Lemma-Definition~\ref{lem-def:1},
  supplemented by the requirement that there exist morphisms
  \begin{equation}
    \label{eq:59}
    \sigma_i\colon \stF_i \lto \grstC\vert_{U_i}
  \end{equation}
  such that given the morphism $\xi_{ij}$ in~\eqref{eq:29} there is a
  morphism of $(\grstB,\grstC)$-torsors
  \begin{equation}
    \label{eq:61}
    (\xi_{ij},t_{ij}) \colon (\stF_j,\sigma_j)\vert_{U_{ij}} \lto
    (\stF_i,\sigma_i)\vert_{U_{ij}}
  \end{equation}
  satisfying a triangle analogous to~\eqref{eq:60}, namely:
  \begin{equation*}
    \xymatrix@+1pc{
      \stF_j \ar@/^/[dr]^{\sigma_j}_{}="pi" \ar[d]_{\xi_{ij}} \\
      \stF_i \ar[r]_{\sigma_i}^{}="sigmapr"
      \ar@{=>}@/_0.1pc/ _{t_{ij}} "pi";"sigmapr" &
      \grstC
    }      
  \end{equation*}
\end{lemma-def}
\begin{proof}
  We need only observe that a morphism
  \begin{equation*}
    \lambda_{\ast}(\stE_{ij})\cprod{\grstB}\!\stF_j \lto
    \grstC\vert_{U_{ij}}
  \end{equation*}
  can equivalently be seen as a morphism of $\grstC$-torsors:
  \begin{equation*}
    (\lambda_{\ast}(\stE_{ij})\cprod{\grstB}\!\stF_j)
    \cprod{\grstB}\! \grstC\vert_{U_{ij}}
    \lto
    \grstC\vert_{U_{ij}}\,.
  \end{equation*}
  But we have
  \begin{equation*}
    (\lambda_{\ast}(\stE_{ij})\cprod{\grstB}\!\stF_j)
    \cprod{\grstB}\! \grstC\vert_{U_{ij}} \iso
    \lambda_{\ast}(\stE_{ij})\cprod{\grstB}(\stF_j
    \cprod{\grstB}\! \grstC\vert_{U_{ij}}) \iso
    \stF_j \cprod{\grstB}\! \grstC\vert_{U_{ij}}
  \end{equation*}
  since $\mu_\ast\circ \lambda_\ast\iso
  (\mu\circ\lambda)_\ast\iso \mathbf{1}$.
\end{proof}
The argument of the proof also implies that two 1-cocycles
$(\stE_{ij},\stF_i,\sigma_i)$ and
$(\stE'_{ij},\stF'_i,\sigma'_i)$ with values in~\eqref{eq:54}
ought to be considered equivalent if the same conditions of
Definition~\ref{def:9} are satisfied, with the additional
requirement that the morphism~\eqref{eq:35} induces a morphism of
$(\grstB,\grstC)$-torsors
\begin{equation*}
  (\stF_i,\sigma_i) \lto (\stF'_i,\sigma'_i)\,.
\end{equation*}
We leave to the reader the task of spelling out the rest of the
details.

The next results combines the generalizations of
Theorems~\ref{thm:3} and~\ref{thm:4} to the present case. Large
parts of the proof can be simply carried over, therefore we will
be sketchy.
\begin{theorem}
  \label{thm:5}
  \hfill\par
  \begin{enumerate}
  \item\label{item:13} Equivalence classes of
    2-$(\grstA,\grstB,\grstC)$-gerbes are classified by the
    (pointed) set
    \begin{equation*}
      \HHH^1(X,\tors (\grstA)\to \tors (\grstB) \to \tors(\grstC))
    \end{equation*}
    of equivalence classes of 1-cocycles with values in the
    complex~\eqref{eq:54}, according to the
    Lemma-Definition~\ref{lem-def:2}.
  \item\label{item:14} If the braided \gr-stacks are all strict
    and associated to abelian cross modules as in
    sect.~\ref{sec:complexes-braided-gr}, then the above pointed
    set of equivalence classes is actually in 1--1 correspondence
    with the hypercohomology group
    \begin{equation*}
      \HHH^3(X,A\to B\oplus G\to C\oplus H\to K)
    \end{equation*}
    where we recognize the cone (shifted by 1) of the
    morphism~\eqref{eq:55}.
  \end{enumerate}
\end{theorem}
\begin{proof}
  Let $(\tgG,\Tilde J)$ be a 2-gerbe over $\cat{C}/X$ bound by
  the complex~\eqref{eq:53}. Let us make the usual choice of a
  cover $\cover{U}_X$, to be enhanced to a hypercover below. The
  proof of Part~\ref{item:13} rests upon the choice of a
  decomposition of $\tgG$ with respect to a collection of objects
  $x_i\in \ob\tgG_{U_i}$. By applying $\Tilde J$ we obtain
  $(\grstB,\grstC)$-torsors $\Tilde J (x_i)=\Hat\stF_i\coin
  (\stF_i,\sigma_i)$ and morphisms
  \begin{equation*}
    \stE_{ij} \lto
    \catHom(\Hat\stF_j\vert_{U_{ij}},\Hat\stF_i\vert_{U_{ij}}).
  \end{equation*}
  Forgetting the morphisms into $\grstC$ gives the underlying
  functor in $\tors (\grstB)$, therefore Part~\ref{item:13}
  follows from Thm.~\ref{thm:3} (or rather, its proof) and the
  argument made in the proof of~\ref{lem-def:2} to handle the
  extra morphisms into $\grstC$.
  
  The proof of Part~\ref{item:14} is more laborious, but only
  computationally so.  Fortunately everything that was done in
  the proof of Thm.~\ref{thm:4} can be transported verbatim here,
  so that we only have to deal with the extra data ensuing from
  the $(\grstB,\grstC)$-torsor.
  
  Our first task is to rewrite the classifying 1-cocycle with
  values in~\eqref{eq:54} from Part~\ref{item:13} in terms of a
  cocycle with values in the complex of \gr-stacks~\eqref{eq:53}.
  As before, this is accomplished by decomposing the cocycle
  $(\stE_{ij},\stF_i,\sigma_i)$ with respect to a choice of
  objects subordinated to a given hypercover. As in the proof of
  Thm.~\ref{thm:4}, we refine $\cover{U}_X$ by $(U^\alpha_{ij}\to
  U_{ij})_{\alpha \in A^\alpha_{ij}}$.  We also keep all the
  choices and notations made there.
  
  Recall that we had obtained the quintuple~\eqref{eq:47} which
  we rewrite here for convenience:
  \begin{equation*}
    \bigl( h^\alpha_{ij}\,, q^{\alpha\beta}_{ijk}\,,
    m^{\alpha\beta\gamma}_{ijk}\,,
    g^{\alpha\beta\gamma}_{ijk}\,,
    \nu^{\alpha\beta\delta\gamma\eta\epsilon}_{ijkl}
    \bigr)
  \end{equation*}
  where $h^\alpha_{ij}\,, q^{\alpha\beta}_{ijk}$ are objects of
  $\grstB$, $m^{\alpha\beta\gamma}_{ijk}$ are morphisms of
  $\grstB$, and $g^{\alpha\beta\gamma}_{ijk}$ and
  $\nu^{\alpha\beta\delta\gamma\eta\epsilon}_{ijkl}$ are objects
  and morphisms of $\grstA$, respectively.  They satisfy the
  cocycle conditions given by the equations~\eqref{eq:42}
  and~\eqref{eq:26}.
  
  Since the morphism $\sigma_i \colon \stF_i\to
  \grstC\vert_{U_i}$ are global over $U_i$, the arguments in
  sect.~\ref{sec:grstb-grstc-torsors} imply that there are
  objects ${{z_i}^\alpha}_j\in\ob\grstC\vert_{U^\alpha_{ij}}$ and
  morphisms $t^\alpha_{ij}$ and $\rho^{\alpha\beta}_{jik}$ in
  $\grstC\vert_{U^\alpha_{ij}}$ and
  $\grstC\vert_{U^{\alpha\beta}_{ijk}}$ such that:
  \begin{subequations}
    \label{eq:64}
    \begin{gather}
      \label{eq:65}
      t^\alpha_{ij} \colon {{z_j}^\alpha}_i \lisoto
      \mu (h^\alpha_{ij})\cdot {{z_i}^\alpha}_j\\
      \intertext{and}
      \label{eq:62}
      \rho^{\beta\alpha}_{kij} \colon
      {{z_i}^\alpha}_j \lisoto \mu (q^{\alpha\beta}_{jik})\cdot
      {{z_i}^\beta}_k\,.
    \end{gather}
  \end{subequations}
  Both equations~\eqref{eq:64} are obtained by applying the
  morphisms $\sigma_i$, $\sigma_j$, etc., namely the triangle
  right after~\eqref{eq:61}, to eqns. \eqref{eq:36} and
  \eqref{eq:37}, respectively. We have used the relation
  $q^{\alpha\beta}_{jik} \iso (q^{\beta\alpha}_{kij})^\ast$,  easily
  derived from~\eqref{eq:37}, 
  where $(\cdot)^\ast$ denotes the quasi-inverse.
  The final piece of
  the cocycle condition is a relation for the morphisms
  $t^\alpha_{ij}$ and $\rho^{\beta\alpha}_{kij}$ which is
  computed by passing from ${{z_k}^\gamma}_i$ to
  ${{z_i}^\alpha}_j$ in two different ways. Either as:
  \begin{equation}
    \label{eq:63}
    \rho^{\alpha\gamma}_{jik}\circ t^\gamma_{ik} \colon
    {{z_k}^\gamma}_i \lisoto
    \mu (h^\gamma_{ki})\cdot \mu (q^{\gamma\alpha}_{kij})\cdot
    {{z_i}^\alpha}_j\,,
  \end{equation}
  or as:
  \begin{equation}
    \label{eq:66}
    t^\alpha_{ij} \circ \rho^{\alpha\beta}_{ijk} \circ
    t^\beta_{jk} \circ \rho^{\beta\gamma}_{jki} \colon
    {{z_k}^\gamma}_i \lisoto
    \mu (q^{\gamma\beta}_{ikj}) \cdot \mu (h^\beta_{kj}) \cdot
    \mu (q^{\beta\alpha}_{kji})\cdot \mu (h^\alpha_{ji} )\cdot
    {{z_i}^\alpha}_j\,,
  \end{equation}
  where, as before, we are ignoring the various associator
  isomorphisms and natural transformations associated with $\mu$.
  
  If we replace the three middle terms in the
  right hand side of~\eqref{eq:66} using \eqref{eq:40} and the relations
  $\mu\circ\lambda (g^{\beta\alpha\gamma}_{kji})\iso o_\grstC$
  and $q^{\gamma\beta}_{ikj}\cdot q^{\beta\gamma}_{jki}\iso
  o_\grstC$, where $o_\grstC$ is the unit element of $\grstC$, we
  find
  \begin{equation*}
    \mu (m^{\beta\alpha\gamma}_{kji}) \circ
    t^\alpha_{ij} \circ \rho^{\alpha\beta}_{ijk} \circ
    t^\beta_{jk} \circ \rho^{\beta\gamma}_{jki} \colon
    {{z_k}^\gamma}_i \lisoto
    \mu (h^\gamma_{ki})\cdot \mu (q^{\gamma\alpha}_{kij})\cdot
    {{z_i}^\alpha}_j\,.
  \end{equation*}
  Comparing with~\eqref{eq:63}, we obtain the desired relation:
  \begin{equation}
    \label{eq:68}
    \mu (m^{\beta\alpha\gamma}_{kji}) \circ
    t^\alpha_{ij} \circ \rho^{\alpha\beta}_{ijk} \circ
    t^\beta_{jk} \circ \rho^{\beta\gamma}_{jki}
    = \rho^{\alpha\gamma}_{jik}\circ t^\gamma_{ik}\,.
  \end{equation}
  Thus, starting from the cocycle $(\stE_{ij},\stF_i,\sigma_i)$
  with values in~\eqref{eq:54}, the
  corresponding cocycle with values in the complex~\eqref{eq:53}
  is the 8-tuple
  \begin{equation}
    \label{eq:69}
    \bigl({{z_i}^\alpha}_j\,, t^\alpha_{ij}\,,
    \rho^{\alpha\beta}_{kij}\,,
    h^\alpha_{ij}\,, q^{\alpha\beta}_{ijk}\,,
    m^{\alpha\beta\gamma}_{ijk}\,,
    g^{\alpha\beta\gamma}_{ijk}\,,
    \nu^{\alpha\beta\delta\gamma\eta\epsilon}_{ijkl}
    \bigr)
  \end{equation}
  satisfying the conditions~\eqref{eq:42}, \eqref{eq:26},
  \eqref{eq:64}, and~\eqref{eq:68}.
  
  The proof will be complete when we specialize~\eqref{eq:69} and
  the relations it satisfies to the case where all the involved
  \gr-stacks are Picard and associated to the abelian crossed
  modules introduced in sect.~\ref{sec:complexes-braided-gr}.
  This means that $\grstC=\tors (C,K)$, where the underlying
  groupoid $C\times K\rightrightarrows K$ has source and target
  maps given by $(c,z)\to z$ and $(c,z) \to \tau (c)z$,
  respectively, and similarly for $\grstA$ and $\grstB$ with the
  appropriate notations and relations, which we can lift directly
  from the proof of Thm.~\ref{thm:4}, eqns.~\eqref{eq:43}.  Thus,
  the objects ${{z_i}^\alpha}_j$ will be identified with sections
  of the group $K\vert_{U^\alpha_{ij}}$, and we also need to
  introduce sections $c^\alpha_{ij}$ of $C\vert_{U^\alpha_{ij}}$
  and $l^{\alpha\beta}_{ijk}$ of $C\vert_{U^{\alpha\beta}_{ijk}}$
  to account for the morphisms $t^\alpha_{ij}$ and
  $\rho^{\alpha\beta}_{ijk}$, respectively.

  With these provisions, the 8-tuple~\eqref{eq:69} becomes
  \begin{equation}
    \label{eq:67}
    \bigl({{z_i}^\alpha}_j\,,c^\alpha_{ij}\,,l^{\alpha\beta}_{kij}\,,
    h^\alpha_{ij}\,, q^{\alpha\beta}_{ijk}\,,
    b^{\alpha\beta\gamma}_{ijk}\,,
    g^{\alpha\beta\gamma}_{ijk}\,,
    a^{\alpha\beta\delta\gamma\eta\epsilon}_{ijkl}
    \bigr)\,,    
  \end{equation}
  and equations~\eqref{eq:64} and~\eqref{eq:68} become
  \begin{subequations}
    \label{eq:70}
    \begin{gather}
      \label{eq:71}
      \tau (c^\alpha_{ij}) \cdot  {{z_j}^\alpha}_i =
      v (h^\alpha_{ij})\cdot {{z_i}^\alpha}_j\\
      \label{eq:72}
      \tau (l^{\beta\alpha}_{kij}) \cdot {{z_i}^\alpha}_j =
      v (q^{\alpha\beta}_{jik})\cdot {{z_i}^\beta}_k \\
      \label{eq:73}
      g (b^{\beta\alpha\gamma}_{kji}) \cdot
      c^\alpha_{ij} \cdot l^{\alpha\beta}_{ijk} \cdot
      c^\beta_{jk} \cdot l^{\beta\gamma}_{jki}
      = l^{\alpha\gamma}_{jik}\cdot c^\gamma_{ik}\,.      
    \end{gather}
  \end{subequations}
  The full cocycle condition for the 8-tuple~\eqref{eq:67} is
  then given by eqns.~\eqref{eq:70} \emph{plus}
  eqns.~\eqref{eq:43}, which we rewrite here:
  \begin{align*}
    \delta (a^{\alpha\beta\delta\gamma\eta\epsilon}_{ijkl})
    \cdot g^{\alpha\beta\gamma}_{ijk}\cdot
    g^{\gamma\delta\epsilon}_{ikl} & =
    g^{\beta\delta\eta}_{jkl}\cdot
    g^{\alpha\eta\epsilon}_{ijl}\,,\\
    \sigma (b^{\alpha\beta\gamma}_{ijk})\cdot
    h^\alpha_{ij}\cdot q^{\alpha\beta}_{ijk} \cdot h^\beta_{jk}
    & =
    u (g^{\alpha\beta\gamma}_{ijk})\cdot
    q^{\alpha\gamma}_{jik} \cdot h^\gamma_{ik} \cdot
    q^{\gamma\beta}_{ikj}\,,\\
    f (a^{\alpha\beta\delta\gamma\eta\epsilon}_{ijkl})
    \cdot
    b^{\gamma\delta\epsilon}_{ikl}
    \cdot
    b^{\alpha\beta\gamma}_{ijk} & =
    b^{\alpha\eta\epsilon}_{ijl} \cdot
    b^{\beta\delta\eta}_{jkl}\,.
  \end{align*}
  Finally we need also to add the cocycle condition on the
  elements $a^{\alpha\beta\delta\gamma\eta\epsilon}_{ijkl}$.
  
  The amount of typographical decoration provided by the upper
  indices related to the hypercover can be quite daunting.
  Ignoring these indices (that is, reducing everything to the
  \cech case), although potentially less precise from the
  cohomological point of view (cf.\ the discussion
  in~\cite{MR95m:18006}) does shed some light on how the various
  parts are organized.  Without upper indices we need to set
  $q^{\alpha\beta}_{ijk}=1$ and $l^{\alpha\beta}_{ijk}=1$ in the
  above formulas. Thus, the 8-tuple~\eqref{eq:67} becomes a
  sextuple
  \begin{equation*}
    \bigl({z_i}\,,c_{ij}\,, h_{ij}\,, 
    b_{ijk}\,, g_{ijk}\,, a_{ijkl}
    \bigr)
  \end{equation*}
  satisfying the cocycle condition:
  \begin{align*}
    \tau (c_{ij}) \cdot z_j & = v (h_{ij})\cdot z_i\\
    g (b_{kji}) \cdot c_{ij} \cdot c_{jk} & = c_{ik}\,,\\
    \delta (a_{ijkl}) \cdot g_{ijk}\cdot g_{ikl}
    & = g_{jkl}\cdot g_{ijl}\,,\\
    \sigma (b_{ijk})\cdot h_{ij}\cdot h_{jk}
    & =  u (g_{ijk})\cdot h_{ik} \,\\
    f (a_{ijkl}) \cdot b_{ikl} \cdot b_{ijk}
    & = b_{ijl} \cdot  b_{jkl}\,.
  \end{align*}
  Now write the cone of the the morphism of
  complexes~\eqref{eq:55} in the form:
  \begin{equation*}
    A\xrightarrow{%
      \left(\begin{smallmatrix}
          f\\\delta
        \end{smallmatrix}\right)}
    B\oplus G \xrightarrow{%
      \left(\begin{smallmatrix}
          g&1\\\sigma&u^{-1}
        \end{smallmatrix}\right)}
    C\oplus H\xrightarrow{%
      \left(\begin{smallmatrix}
          \tau & v^{-1}
        \end{smallmatrix}\right)}
    K
  \end{equation*}
  It can now be seen in a direct way that the
  8-tuple~\eqref{eq:67} (or its simplified \cech version) indeed
  defines a 3-cocycle with values in the cone of~\eqref{eq:55}.
  This is straightforward and left to the reader. We will also
  omit the verification that passing to an equivalent torsor
  1-cocycle $(\stE'_{ij},\stF'_i,\sigma'_i)$ representing
  $(\tgG,\Tilde J)$, we obtain an equivalent 3-cocycle.
\end{proof}
An even more special case of Theorem~\ref{thm:5}--~\ref{item:14}
is when the diagram~\eqref{eq:55} reduces to the complex $A
\overset{f}{\lto} B\overset{g}{\lto} C$. Let $(\tgG,\Tilde J)$
be a 2-gerbe over $\cat{C}/X$ bound by
$\tors (A)\to \tors (B) \to \tors (C)$.
By comparing the classifying cocycles we immediately obtain the
following
\begin{corollary}
  \label{cor:1}
  $(\tgG,\Tilde J)$ is equivalent to a 2-$(A,B,C)$-gerbe in the
  sense of Definition~\ref{def:10}.
\end{corollary}

\section{Applications}
\label{sec:geom-appl}

In this section we will address a few questions about the
correspondence between certain Hermitian Deligne Cohomology
groups and equivalence classes of 2-gerbes equipped with various
geometric structures of the type described in the previous
sections.

For consistency with the results of~\cite{MR2142353} and previous
work in Deligne cohomology we will be placing $\ZZ(p)_X$ in
degree zero, therefore all cohomology degrees will be shifted up
in comparison with those appearing in the previous sections.

\subsection{Truncated hermitian Deligne complexes}
\label{sec:trunc-herm-deligne}

Beside the Hermitian Deligne complexes recalled in
sect.~\ref{sec:vari-deligne-compl}, we need two more complexes we
introduced the in~\cite{MR2142353}, namely
\begin{equation*}
  \complex{\Gamma (2)} =  \cone
  \begin{pmatrix}
    \xymatrix{%
      \ZZ(2) \ar[r] & \sho{X} \ar[r] \ar[d]
      & \shomega[1]{X} \ar[d]\\
      & \she[0]{X}(1) \ar[r] & \she[1]{X}(1)}
  \end{pmatrix}[-1]\,,
\end{equation*}
plus the truncation
\begin{equation*}
  \complex{\Tilde\Gamma (2)} =  \cone
  \begin{pmatrix}
    \xymatrix{%
      \ZZ(2) \ar[r] & \sho{X} \ar[r] \ar[d]
      & \shomega[1]{X} \ar[d]\\
      & \she[0]{X}(1) \ar[r] & 0}
  \end{pmatrix}[-1]\,,
\end{equation*}
where the maps are the same as in the corresponding places in the
diagram defining $\ndhh{X}{2}$.  (It is convenient to pass, from
now on, to an additive notation.)  Note that $\complex{\Gamma
  (2)}$ is an obvious truncation of the Hermitian Deligne complex
$\ndhh{X}{2}$, while $\complex{\Tilde\Gamma (2)}$ is in turn a
truncation of $\complex{\Gamma (2)}$.  These two complexes were
introduced as part of the effort to analyze the interplay and
compatibility of different types of differential geometric
structures on 2-gerbes.  Indeed, it can be shown that
$\complex{\Gamma (2)}$ arises from the diagram of complexes:
\begin{equation*}
  \deligne{X}{\ZZ}{2} \lto \complex{C (2)} \longleftarrow
  \tate\otimes \ndhh{X}{1}
\end{equation*}
in the sense of \cite{bei:hodge:coho}, namely as the cone of the
difference of the two maps.  Here $\complex{C (2)}$ is the
complex
\begin{equation*}
  \ZZ (2)_X  \lto \sho{X} \lto \she[1]{X}(1) \,.
\end{equation*}
Similarly, $\complex{\Tilde\Gamma (2)}$ arises in the same way
from the diagram:
\begin{equation*}
  \deligne{X}{\ZZ}{2} \lto
  \deligne{X}{\ZZ}{1} \longleftarrow
  \tate\otimes \ndhh{X}{1}  \,,
\end{equation*}
where the two maps are just the forgetful maps.  We have
repeatedly seen how the complexes $\deligne{X}{\ZZ}{2}$ (resp.\ 
$\ndhh{X}{1}$) intervene in the definition of connective (resp.\ 
hermitian) structures. Note, however, that the above complexes
and their geometric role was introduced rather informally in the
context of \cite{MR2142353}.  The results of
sect.~\ref{sec:geom-interpr} provide a more rigorous footing.

We quote from~\cite{MR2142353} the following exact
sequences. From the definitions we immediately have:
\begin{gather*}
   0\lto \she[1]{X}(1)[-3] \lto
   \complex{\Gamma (2)} \lto
   \complex{\Tilde\Gamma (2)}\lto 0\\
   \intertext{and}
   0 \lto \she[2]{X}(1)\cap \sha[1,1]{X} [-4]
  \lto \ndhh{X}{2} \lto
  \complex{\Gamma (2)} \lto 0\,.
\end{gather*}
Furthermore, using the standard arguments, as well as the
softness of $\she[1]{X}(1)$, $\she[2]{X}(1)$, and $\sha[1,1]{X}$,
we obtain:
\begin{gather*}
  \dotsb \lto E^1_X(1) \lto
  \HHH^3(X,\complex{\Gamma (2)}) \lto
  \HHH^3(X,\complex{\Tilde\Gamma (2)}) \lto 0 \\
  \dotsb \lto E^2_X(1)\cap A^{1,1}_X \lto \dhhH[4]{X}{2} \lto
  \HHH^4 (X, \complex{\Gamma (2)})\lto 0
\end{gather*}
and the isomorphism
\begin{equation*}
  \HHH^k(X,\complex{\Gamma (2)}) \iso
  \HHH^k(X,\complex{\Tilde\Gamma (2)})\,,
  \quad k\geq 4\,.
\end{equation*}

\subsection{Geometric interpretation of some cohomology groups}
\label{sec:geom-interpr}

Observe that using $\sho{X}/\ZZ(2)_X\iso \sho{X}^\unit$, the
complex $\complex{\Gamma (2)}$ can be identified (modulo the
index shift) with the cone of the square
\begin{equation}
  \label{eq:74}
  \vcenter{%
  \xymatrix{%
    \sho{X}/\ZZ(2)_X \ar[r] \ar[d]
    & \shomega[1]{X} \ar[d]\\
    \she[0]{X}(1) \ar[r] & \she[1]{X}(1)}}
\end{equation}
and similarly for $\complex{\Gamma (2)}$ by replacing
$\she[1]{X}(1)$ with $0$:
\begin{equation}
  \label{eq:76}
  \vcenter{%
  \xymatrix{%
    \sho{X}/\ZZ(2)_X \ar[r] \ar[d]
    & \shomega[1]{X} \ar[d]\\
    \she[0]{X}(1) \ar[r] & 0}}  
\end{equation}
Both cases correspond to the diagram~\eqref{eq:22}.

To make contact with the contents of
sect.~\ref{sec:2-gerbes-bound}, let us set
\begin{equation*}
  \grstA = \tors (\sho{X}/\ZZ(2)_X, \she[0]{X}(1))\,,
  \quad
  \grstB = \tors (\shomega[1]{X}, \she[1]{X}(1))
\end{equation*}
so that we have the equivalences
\begin{gather*}
  \tors (\grstA) \lisoto
  \gerbes (\sho{X}/\ZZ(2)_X, \she[0]{X}(1))\\
  \intertext{and}
  \tors (\grstB) \lisoto
  \gerbes (\shomega[1]{X}, \she[1]{X}(1))\,.
\end{gather*}
Using Theorem~\ref{thm:6} and Proposition~\ref{prop:6} we find
the following alternative characterization of
$\sho{X}^\unit$-gerbes with compatible hermitian and connective
structure:
\begin{corollary}
  \label{cor:4}
  The group $\HHH^3(X,\complex{\Gamma (2)})$ classifies
  equivalence classes of \\ $(\sho{X}/\ZZ(2)_X,
  \she[0]{X}(1))$-gerbes, that is, hermitian gerbes in the sense
  of~\ref{sec:hermitian-structures}, which become neutral as
  $(\shomega[1]{X}, \she[1]{X}(1))$-gerbes.
\end{corollary}
Of course, the other possible but entirely equivalent statement
would have been that the cohomology group under scrutiny
classifies $(\tors (\sho{X}/\ZZ(2)_X, \she[0]{X}(1)),$ $\tors
(\shomega[1]{X}, \she[1]{X}(1)))$-torsors.  We leave to the
reader the task of formulating a similar statement for the
complex $\complex{\Tilde\Gamma (2)}$.
\begin{remark}
  A short remark is in order about other possible ways of
  interpreting the same cohomology group.  As noted, we can take
  advantage of the symmetry of the square~\eqref{eq:74} in the
  sense explained in Remark~\ref{rem:2}, and modify things
  accordingly. This preserves the cone, namely $\complex{\Gamma
    (2)}$, and does not alter the classifying group. It
  \emph{does} change the \gr-stacks $\grstA$ and $\grstB$, but
  ultimately not the fact that we are dealing with
  $\sho{X}^\unit$-gerbes.
\end{remark}
\begin{remark}
  The above characterization (and the general theory it descends
  from) provides a finer description of the geometric objects
  corresponding whose equivalence classes correspond to the group
  elements when the coefficient complex come from a cone. Had we
  just used the complex $\complex{\Gamma (2)}$ as it stands, we
  would have been in the rather awkward position of calling
  something with values in $\shomega[1]{X}\oplus \she[0]{X}(1)$ a
  ``connective structure,'' a fact that does not seem to sit well
  with the degrees.
\end{remark}

The corresponding result for 2-gerbes provides a similar
interpretation for the group of equivalence classes of
2-$\sho{X}^\unit$-gerbes with compatible hermitian and connective
structure defined in \cite{MR2142353}. It is an immediate
consequence of Theorem~\ref{thm:4} as follows:
\begin{corollary}
  \label{cor:2}
  Elements of the hypercohomology group $\HHH^4(X,
  \complex{\Gamma (2)})$ are in 1--1 correspondence with
  equivalence classes of 2-gerbes on $X$ bound by the
  square \eqref{eq:74} (in the sense of remark~\ref{rem:2}).  A
  similar conclusion holds by replacing $\complex{\Gamma (2)}$
  with $\complex{\Tilde\Gamma (2)}$.
\end{corollary}
Note that a remark concerning the square similar to the one just
made for gerbes holds in this case as well.

In a similar vein to what was just done for the complex
$\complex{\Gamma (2)}$, we can identify $\ndhh{X}{2}$ defined in
eq.~\eqref{eq:11} with the cone of
\begin{equation}
  \label{eq:75}
  \vcenter{%
  \xymatrix{%
    \sho{X}/\ZZ(2)_X \ar[r] \ar[d]
    & \shomega[1]{X} \ar[r] \ar[d] & 0 \ar[d]\\
    \she[0]{X}(1) \ar[r] & \she[1]{X}(1) \ar[r]
    & \she[2]{X}(1)\cap \sha[1,1]{X}}}
\end{equation}
which will correspond to the diagram~\eqref{eq:55}. We have
explicitly written the last column as $0\to \she[2]{X}(1)\cap
\sha[1,1]{X}$ in order to emphasize the correspondence. To take
the point of view of sect.~\ref{sec:2-gerbes-bound-2}, let us
introduce the \emph{discrete} \gr-stack
\begin{equation*}
  \grstC = \tors(0,\she[2]{X}(1)\cap \sha[1,1]{X})
  \iso \she[2]{X}(1)\cap \sha[1,1]{X}\,,
\end{equation*}
namely the only morphisms are the identity maps.  Note since
$\grstC$ is discrete, then the corresponding 2-gerbe is discrete
as well, that is we have:
\begin{equation*}
  \tors (\grstC) \iso \tors (\she[2]{X}(1)\cap \sha[1,1]{X})\,.
\end{equation*}
In other words, it has only identity 2-arrows, and it corresponds
to the neutral gerbe of torsors.

Now, as a consequence of Theorem~\ref{thm:5} we obtain the
following general geometric interpretation for the hermitian
Deligne cohomology group:
\begin{corollary}
  \label{cor:3}
  Elements of the hermitian Deligne cohomology group
  $\dhhH[4]{X}{2}$ are in 1--1 correspondence with equivalence
  classes of 2-gerbes on $X$ bound by the diagram~\eqref{eq:75},
  that is, by the complex~\eqref{eq:53} of \gr-stacks associated
  to the columns of~\eqref{eq:75}.
\end{corollary}

\subsection{Geometric construction of some cup products}
\label{sec:geom-interpr-some}

\subsubsection{}
\label{sec:cup-coho}

If $(\sheaf{L},\rho)$ and $(\sheaf{M},\sigma)$ are two metrized
line bundles (invertible sheaves) over $X$, their isomorphism
classes determine elements of $\dhhH[2]{X}{1}\iso \Pichat X$.
According to the last paragraph of
sect.~\ref{sec:vari-deligne-compl}, the cup product
\begin{math}
  [\sheaf{L},\rho] \cup [\sheaf{M},\sigma]
\end{math}
in hermitian Deligne cohomology will land in $\dhhH[4]{X}{2}$.

It is known from the works of Brylinski and McLaughlin
(\cite{brymcl:deg4:I,brymcl:deg4:II,bry:quillen}) that the
corresponding problem in standard Deligne cohomology has a
geometric interpretation: there is a 2-gerbe
$\tame{\sheaf{L}}{\sheaf{M}}$ bound by $\deligne{X}{\ZZ}{2}$
whose class is the cup product $[\sheaf{L}]\cup [\sheaf{M}] \in
\delH[4]{X}{\ZZ}{2}$ of the elements in $\Pic X$ determined by
$\sheaf{L}$ and $\sheaf{M}$.  Similarly, in \cite{MR2142353} we
constructed a modified cup product
\begin{equation*}
  \Pic X\otimes \Pic X \lto \dhhH[4]{X}{1}
\end{equation*}
and a corresponding ``tame symbol,'' namely a 2-gerbe
$\tamehh{\sheaf{L}}{\sheaf{M}}$ bound by $\ndhh{X}{1}$. It turns
out that both symbols have the ``same'' (in the sense of
equivalent) underlying 2-gerbe, obtained by applying a suitable
forgetful functor to both sides. In other words we have a lift
\begin{equation*}
  \Pic X\otimes \Pic X \lto \HHH^4(X,\complex{\Tilde\Gamma (2)})
\end{equation*}
and it follows from the material recalled in
sect.~\ref{sec:trunc-herm-deligne} that at the level of
cohomology the latter lift can be arranged to take values in
$\HHH^4(X,\complex{\Gamma (2)})$. Thus, from a pair of invertible
sheaves $\sheaf{L}$ and $\sheaf{M}$ we obtain (canonically) a
2-gerbe bound by the square~\eqref{eq:76}, and (non canonically)
by way of softness of one of the sheaves involved, a 2-gerbe
bound by the square~\eqref{eq:74}.

The cohomology exact sequences recalled in
sect.~\ref{sec:trunc-herm-deligne}, and the fact that truncation
will map the diagram~\eqref{eq:75} to the square~\eqref{eq:74},
and then to the square~\eqref{eq:76}, show that the 2-gerbe bound
by~\eqref{eq:75}  corresponding to the cup product 
\begin{math}
  [\sheaf{L},\rho] \cup [\sheaf{M},\sigma]
\end{math}
will provide the required lift.

\subsubsection{}
\label{sec:cup-geom}

We will denote by $\tamehh{\sheaf{L}}{\sheaf{M}}\sphat$ the
2-gerbe bound by \eqref{eq:75} corresponding to the cup product
of the two metrized line bundle.  Let us sketch the geometric
construction of such 2-gerbe borrowing on the corresponding
constructions of \cite{brymcl:deg4:II} and \cite{MR2142353}.

If we work locally with respect to some cover $U\to X$ of $X$,
any 2-$\grstA$-gerbe $\tgG$ will be a 2-gerbe of torsors, namely
there is an equivalence:
\begin{equation*}
  \tgG_U \lisoto \tors (\grstA\vert_U)\lisoto
  \gerbes (\sho{X}/\ZZ(2)_X\vert_U,\she[0]{X}(1)\vert_U)\,,
\end{equation*}
where the latter equivalence follows from Theorem~\ref{thm:6}.
Thus if $\tgG$ is bound by the complex of \gr-stacks determined
by the diagram~\eqref{eq:75}, with $\grstA$, $\grstB$, and
$\grstC$ as in sect.~\ref{sec:geom-interpr}, then locally it has
the form
\begin{equation*}
  \tors (\grstA\vert_U) \lto \tors(\grstB\vert_U,\grstC\vert_U)\,.
\end{equation*}
Note that, thanks to~\ref{prop:3}-\ref{item:12},
Proposition~\ref{prop:6}, and to the  fact that in the relevant
diagram one of the group is zero, we have an equivalence:
\begin{equation*}
  \tors(\grstB\vert_U,\grstC\vert_U) \lisoto
  \gerbes (\shomega[1]{X}\vert_U,\she[1]{X}(1)\vert_U,
  \she[2]{X}(1)\cap \sha[1,1]{X} \vert_U)\,.
\end{equation*}
Let $\tameangle{\sheaf{L}}{\sheaf{M}}$ denote the underlying
2-gerbe of both $\tame{\sheaf{L}}{\sheaf{M}}$ and
$\tamehh{\sheaf{L}}{\sheaf{M}}$. The local objects of
$\tameangle{\sheaf{L}}{\sheaf{M}}$ over $U$ are in 1--1
correspondence with the non-vanishing sections of
$\sheaf{L}\vert_U$. We may denote such a section $s$, which
thought of as an object, by $\tameangle{s}{\sheaf{M}}$.

The choice of $s$ will determine an $\grstA\vert_U$-torsor as
follows. Given any other non-vanishing section $s'$, write
$s=s'\cdot g$ where $g\in \sho{X}/\ZZ(2)$. The
$\grstA\vert_U$-torsor $\catHom (s,s')$ can be
 identified with the
$(\sho{X}/\ZZ(2)_X\vert_U,\she[0]{X}(1)\vert_U)$-gerbe
$\tamehh{g}{\sheaf{M}}$ by the above equivalence.  Let us denote
by $\tameangle{g}{\sheaf{M}}$ the underlying
$\sho{X}/\ZZ(2)_X$-gerbe. Recall from
\cite{brymcl:deg4:II,MR2142353} that its objects over $U$ are in 1--1
correspondence with the non-vanishing sections $t$ 
of $\sheaf{M}\vert_U$, denoted $\tameangle{g}{t}$, and that an
arrow $\phi \colon \tameangle{g}{t} \to \tameangle{g}{t'}$ is
identified with a section of Deligne's torsor $\tame{g}{g'}$,
where $t=t'\cdot g'$, for $g'$ a section of $\sho{X}/\ZZ(2)_X$
over $U$, see \cite{del:symbole}. We reserve the notation
$\tame{g}{g'}$ for the same torsor equipped with the connection
defined in \loccit, whereas the notation $\tamehh{g}{g'}$ denotes
the same underlying torsor equipped with the hermitian structure
defined in \cite{MR2142353}.

To summarize, to define $\tamehh{\sheaf{L}}{\sheaf{M}}\sphat$ we
have to define a 2-functor $\Tilde J_U$ from $\tgG_U$ to the fibered
2-category of gerbes bound by
\begin{equation*}
  \shomega[1]{X}\vert_U \overset{\pi_1}{\lto}
  \she[1]{X}(1)\vert_U \overset{\pi\circ\d}{\lto}
  \she[2]{X}(1)\cap \sha[1,1]{X}\,.
\end{equation*}
To begin with, let us define a 2-functor $J_U$ to $\gerbes
(\shomega[1]{X}\vert_U,\she[1]{X}(1)\vert_U)$ as follows.  To an
object $\tameangle{s}{\sheaf{M}}$ assign the trivial
$\grstB\vert_U$-torsor $T(\grstB\vert_U)\iso \tors
(\shomega[1]{U},\she[1]{U}(1))$. To a 1-arrow
\begin{equation*}
  \tameangle{g}{t} \colon \tameangle{s}{\sheaf{M}} \lto
  \tameangle{s'}{\sheaf{M}}
\end{equation*}
the functor $\tameangle{g}{t}_\ast \colon T(\grstB\vert_U) \to
T(\grstB\vert_U)$ defined as follows: an object of $T
(\grstB\vert_U)$ is identified with an object $(C,\xi)$ of $\tors
(\shomega[1]{U},\she[1]{U}(1))$, where $C$ is a
$\shomega[1]{U}$-torsor which becomes trivial as a
$\she[1]{U}$-torsor by way of $\xi$, which in turn can be
identified with a section of $\she[1]{U}$.  Then we define
$\tameangle{g}{t}_\ast$ by
\begin{equation}
  \label{eq:77}
  \tameangle{g}{t}_\ast  \colon (C,\xi) \lmto (C, \xi + \xi_t)\,,
\end{equation}
where the underlying map on $\tors (\shomega[1]{U})$ is the
identity, and $\xi_t$ is the imaginary 1-form:
\begin{equation}
  \label{eq:78}
  \xi_t = -\onehalf \log\abs{g}\cdot \dc\log \sigma (t)
  +\onehalf \dc\log \abs{g}\cdot \log \sigma (t)\,.
\end{equation}
Here we have used the notation $\sigma (t) =
\abs{t}^2_\sigma$. It is straightforward to verify that this is
compatible with morphisms in $T(\grstB\vert_U)$ and with the
action of $\grstB\vert_U$: if $(D,\eta)$ is an object of $\tors
(\shomega[1]{U},\she[1]{U}(1))$, then 
\begin{equation*}
  (C,\xi)\cdot (D,\eta) = (C\otimes D, \xi+\eta)\,,
\end{equation*}
and obviously this commutes with~\eqref{eq:77}, making it a
morphism of torsors.

Now, if $\phi$ is a section of $\tameangle{g}{g'}$, the
corresponding object of $\tamehh{g}{g'}$ is $(\phi,\norm{\phi})$
where $\norm{\cdot}$ is the hermitian structure given in
\cite{MR2142353}. To it we assign the natural transformation
given by the morphism in $T(\grstB\vert_U)$:
\begin{equation}
  \label{eq:79}
  (\phi,\norm{\phi})_\ast \colon (C, \xi + \xi_t) \lto
  (C, \xi + \xi_{t'})\,,
\end{equation}
which is defined by the underlying map
\begin{equation}
  \label{eq:80}
  \begin{aligned}
    \phi \colon  C &\lto C \\
     c & \lmto c + \phi^{-1}\nabla \phi
  \end{aligned}
\end{equation}
where $\nabla$ is the connection on $\tame{g}{g'}$. From
\cite{del:symbole} we have that locally it has the form
\begin{math}
  -\log g\,\d\log g'\,.
\end{math}
Therefore the section $\xi+\xi_t$ will map to
$\xi+\xi_t+\pi_1(\phi^{-1}\nabla\phi)$ and notice that this
differs from $\xi_{t'}$ by $\tate \d\log \norm{\phi}$, using the
fact that locally $\norm{\cdot}$ is given by $\pi_1(\log
g)\,\log\abs{g'}$. Note that the addition of $\d\log\norm{\phi}$
is just the action of $(\shomega[1]{U},\tate\,\d\log\norm{\phi})$
as an object of $\grstB\vert_U$.

Finally, in order to get the functor $\Tilde J_U$, we need one
more prescription. Namely we define it by assigning to
$\tameangle{s}{\sheaf{M}}$ the
$(\grstB\vert_U,\grstC\vert_U)$-torsor defined as follows. It is
the trivial $\grstB\vert_U$-torsor defined as above equipped with
the morphism
\begin{equation*}
  \tors(\shomega[1]{U},\she[1]{U}(1)) \lto
  \she[2]{U}(1)\cap \sha[1,1]{U}
\end{equation*}
defined by the assignment
\begin{equation}
  \label{eq:81}
  (C,\xi) \lmto \pi(\d\xi)
  -\onefourth\log\rho(s)\,\d\dc\log\sigma(t)
\end{equation}
for every object $(C,\xi)$ of
$\tors(\shomega[1]{U},\she[1]{U}(1))$. Observe that
$\d\dc\log\sigma(t) = c_1(\sheaf{M})$, hence there is no
dependence on $t$. Now, a calculation shows that
\begin{equation*}
  \pi (\d\xi_t) = -\onehalf \log\abs{g} \, \d\dc\log\sigma(t) 
\end{equation*}
so that it is immediately verified that the
assignment~\eqref{eq:81} commutes with the
morphism~\eqref{eq:77}.

With these provisions we have:
\begin{theorem}
  \label{thm:7}
  The class of the 2-gerbe $\tamehh{\sheaf{L}}{\sheaf{M}}\sphat$
  in the cohomology group $\dhhH[4]{X}{2}$ is the cup product
  $[\sheaf{L},\rho] \cup [\sheaf{M},\sigma]$ in hermitian Deligne
  cohomology.
\end{theorem}
\begin{proof}
  It follows immediately from Theorem~\ref{thm:5}, the form of
  the maps in diagrams~\eqref{eq:11} and \eqref{eq:75}, and the
  cup product map
  \begin{equation*}
    \ndhh{X}{1} \otimes \ndhh{X}{1} \lto \ndhh{X}{2}
  \end{equation*}
  given in~\cite{MR2145708}, where the explicit cup-product in
  \cech cohomology is computed.
\end{proof}

\section*{Conclusions}
\addcontentsline{toc}{section}{Conclusions}
\label{sec:conclusions}

We have generalized the concept of ``abelian gerbe bound by a
complex'' to the case of longer coefficient complexes, and to
2-gerbes, where we have used complexes of \gr-stacks of length 3.
We have verified that these 2-gerbes are classified by cohomology
sets of degree 1 with values in the associated complexes of
torsors over these \gr-stacks. We have also obtained, by choosing
appropriate decompositions and hypercovers, that in the strictly
abelian situation the general classification reduces to degree 3
cohomology groups with values in cones of crossed squares, and
other similar diagrams. In all cases we have obtained explicit
cocycles, where we have given their expression in terms of
hypercover, rather than simply in terms of \cech cocycles.

As an application, we have dealt with differential geometric
structures on gerbes and 2-gerbes and questions of geometric
constructions of certain cup-products in hermitian Deligne
cohomology. In particular, we have put certain by now standard
constructions of the concept of connection and curvature in the
general context of gerbe (or 2-gerbe) bound by a complex. We have
further clarified the reason why there seem to exist different
possibilities in defining what a ``hermitian gerbe'' should be
(cf.\ remark~\ref{rem:4}).  Finally, in the last section we have
geometrically constructed a 2-gerbe bound by the hermitian
Deligne complex $\ndhh{X}{2}$ corresponding to the cup product of
two metrized line bundles in hermitian Deligne cohomology.

There are several possible extensions and generalizations of the
work carried out in this paper.  In the case of gerbes, it would
be interesting to remove the abeliannes assumption and work in
the same framework as~\cite{MR0480515} to study extended
structures as coefficients, beyond crossed modules: crossed
squares, 2-crossed complexes, etc.\ come to mind. In particular,
it would be interesting to see whether the idea of phrasing the
notion of connection and curvature in terms of gerbes bound by
complexes extends to the non-abelian case, and how it compares
with other existing approaches (see, e.g.\
\cite{math.AG/0106083}).  In~\cite{MR0480515} a compelling
motivation was to obtain a theory of non-abelian $\H^2$ which
behaved better than Giraud's with respect to group exact
sequences.  Pursuing some these ideas in the case of 2-gerbes
would also be quite interesting.  We hope to return to some of
these issues in future publications.

\bibliography{general} \bibliographystyle{halpha}

\def\polhk#1{\setbox0=\hbox{#1}{\ooalign{\hidewidth
  \lower1.5ex\hbox{`}\hidewidth\crcr\unhbox0}}}
\begin{thebibliography}{Ald05b}

\bibitem[Ald04]{math.CV/0211055} Ettore Aldrovandi.  \newblock On
  hermitian-holomorphic classes related to uniformization, the
  dilogarithm and the liouville action.  \newblock
  \emph{Communications in Mathematical Physics}, 251:27--64,
  2004, math.CV/0211055

\bibitem[Ald05a]{MR2142353} Ettore Aldrovandi.  \newblock
  Hermitian-holomorphic (2)-gerbes and tame symbols.  \newblock
  \emph{J. Pure Appl. Algebra}, 200(1-2):97--135, 2005,
  math.CT/0310027

\bibitem[Ald05b]{MR2145708} Ettore Aldrovandi.  \newblock
  Hermitian-holomorphic {D}eligne cohomology, {D}eligne pairing
  for singular metrics, and hyperbolic metrics.  \newblock
  \emph{Int. Math. Res. Not.}, (17):1015--1046, 2005,
  math.AG/0408118

\bibitem[Be{\u\i}84]{MR86h:11103} Alexander~A. Be{\u\i}linson.
  \newblock Higher regulators and values of {$L$}-functions.
  \newblock In \emph{Current problems in mathematics, Vol. 24},
  Itogi Nauki i Tekhniki, pages 181--238. Akad. Nauk SSSR
  Vsesoyuz. Inst. Nauchn. i Tekhn.  Inform., Moscow, 1984.

\bibitem[Be{\u\i}86]{bei:hodge:coho} Alexander~A. Be{\u\i}linson.
  \newblock Notes on absolute {H}odge cohomology.  \newblock In
  \emph{Applications of algebraic $K$-theory to algebraic
    geometry and number theory, Part I, II (Boulder, Colo.,
    1983)}, pages 35--68. Amer.  Math. Soc., Providence, RI,
  1986.

\bibitem[BG97]{MR99d:14015} Jos{\'e}~I. Burgos~Gil.  \newblock
  Arithmetic {C}how rings and {D}eligne-{B}eilinson cohomology.
  \newblock \emph{J. Algebraic Geom.}, 6(2):335--377, 1997.

\bibitem[BGKK]{math.AG/0404122} Jos{\'e}~I. Burgos~Gil,
  J.~Kramer, and Ulf K{\"u}hn.  \newblock Cohomological
  arithmetic {C}how rings, math.AG/0404122

\bibitem[BM]{math.AG/0106083} Lawrence Breen and William Messing.
  \newblock Differential geometry of gerbes,
  math.AG/0106083

\bibitem[BM94]{brymcl:deg4:I} Jean-Luc Brylinski and
  Dennis~A. McLaughlin.  \newblock The geometry of degree-four
  characteristic classes and of line bundles on loop spaces. {I}.
  \newblock \emph{Duke Math. J.}, 75(3):603--638, 1994.

\bibitem[BM96]{brymcl:deg4:II} Jean-Luc Brylinski and
  Dennis~A. McLaughlin.  \newblock The geometry of degree-$4$
  characteristic classes and of line bundles on loop
  spaces. {I}{I}.  \newblock \emph{Duke Math. J.},
  83(1):105--139, 1996.

\bibitem[Bre90]{MR92m:18019} Lawrence Breen.  \newblock
  Bitorseurs et cohomologie non ab\'elienne.  \newblock In
  \emph{The Grothendieck Festschrift, Vol.\ I}, volume~86 of {\em
    Progr. Math.}, pages 401--476. Birkh\"auser Boston, Boston,
  MA, 1990.

\bibitem[Bre92]{MR93k:18019} Lawrence Breen.  \newblock Th\'eorie
  de {S}chreier sup\'erieure.  \newblock \emph{Ann. Sci. \'Ecole
    Norm. Sup. (4)}, 25(5):465--514, 1992.

\bibitem[Bre94a]{MR95m:18006} Lawrence Breen.  \newblock On the
  classification of {$2$}-gerbes and {$2$}-stacks.  \newblock
  \emph{Ast\'erisque}, (225):160, 1994.

\bibitem[Bre94b]{MR95b:18009} Lawrence Breen.  \newblock
  Tannakian categories.  \newblock In \emph{Motives (Seattle, WA,
    1991)}, volume~55 of \emph{Proc. Sympos.  Pure Math.}, pages
  337--376. Amer. Math. Soc., Providence, RI, 1994.

\bibitem[Bry93]{bry:loop} Jean-Luc Brylinski.  \newblock
  \emph{Loop spaces, characteristic classes and geometric
    quantization}.  \newblock Birkh\"auser Boston Inc., Boston,
  MA, 1993.

\bibitem[Bry94]{MR1317118} Jean-Luc Brylinski.  \newblock
  Holomorphic gerbes and the {B}e\u\i linson regulator.
  \newblock \emph{Ast\'erisque}, (226):8, 145--174, 1994.
  \newblock $K$-theory (Strasbourg, 1992).

\bibitem[Bry99]{bry:quillen} Jean-Luc Brylinski.  \newblock
  Geometric construction of {Q}uillen line bundles.  \newblock In
  \emph{Advances in geometry}, pages 107--146. Birkh\"auser
  Boston, Boston, MA, 1999.

\bibitem[Deb77]{MR0480515} R.~Debremaeker.  \newblock Non abelian
  cohomology.  \newblock \emph{Bull. Soc. Math. Belg.},
  29(1):57--72, 1977.

\bibitem[Del79]{MR546620} Pierre Deligne.  \newblock Vari\'et\'es
  de {S}himura: interpr\'etation modulaire, et techniques de
  construction de mod\`eles canoniques.  \newblock In
  \emph{Automorphic forms, representations and $L$-functions
    (Proc.  Sympos. Pure Math., Oregon State Univ., Corvallis,
    Ore., 1977), Part 2}, Proc. Sympos. Pure Math., XXXIII, pages
  247--289. Amer. Math. Soc., Providence, R.I., 1979.

\bibitem[Del87]{MR89b:32038} Pierre Deligne.  \newblock Le
  d\'eterminant de la cohomologie.  \newblock In \emph{Current
    trends in arithmetical algebraic geometry (Arcata, Calif.,
    1985)}, volume~67 of \emph{Contemp. Math.}, pages
  93--177. Amer. Math.  Soc., Providence, RI, 1987.

\bibitem[Del91]{del:symbole} Pierre Deligne.  \newblock Le
  symbole mod\'er\'e.  \newblock \emph{Inst. Hautes \'Etudes
    Sci. Publ. Math.}, (73):147--181, 1991.

\bibitem[Esn88]{esn:char} H{\'e}l{\`e}ne Esnault.  \newblock
  Characteristic classes of flat bundles.  \newblock
  \emph{Topology}, 27(3):323--352, 1988.

\bibitem[EV88]{esn-vie:del} H{\'e}l{\`e}ne Esnault and Eckart
  Viehweg.  \newblock Deligne-{B}e\u\i linson cohomology.
  \newblock In \emph{Be\u\i linson's conjectures on special
    values of $L$-functions}, pages 43--91. Academic Press,
  Boston, MA, 1988.

\bibitem[Gir71]{MR49:8992} Jean Giraud.  \newblock
  \emph{Cohomologie non ab\'elienne}.  \newblock Springer-Verlag,
  Berlin, 1971.  \newblock Die Grundlehren der mathematischen
  Wissenschaften, Band 179.

\bibitem[Gon04]{gonch:JAMS} Alexander~B. Goncharov.  \newblock
  Polylogarithms, regulators, and arakelov motivic complexes.
  \newblock \emph{Journal of the American Mathematical Society},
  2004, math.AG/0207036  \newblock Posted online November
  2004.

\bibitem[Hak72]{MR0364245} Monique Hakim.  \newblock \emph{Topos
    annel\'es et sch\'emas relatifs}.  \newblock Springer-Verlag,
  Berlin, 1972.  \newblock Ergebnisse der Mathematik und ihrer
  Grenzgebiete, Band 64.

\bibitem[Hit01]{MR1876068} Nigel Hitchin.  \newblock Lectures on
  special {L}agrangian submanifolds.  \newblock In \emph{Winter
    School on Mirror Symmetry, Vector Bundles and Lagrangian
    Submanifolds (Cambridge, MA, 1999)}, volume~23 of
  \emph{AMS/IP Stud. Adv. Math.}, pages
  151--182. Amer. Math. Soc., Providence, RI, 2001.

\bibitem[Lod82]{MR651845} Jean-Louis Loday.  \newblock Spaces
  with finitely many nontrivial homotopy groups.  \newblock
  \emph{J. Pure Appl. Algebra}, 24(2):179--202, 1982.

\bibitem[Mil80]{MR559531} James~S. Milne.  \newblock
  \emph{\'{E}tale cohomology}, volume~33 of \emph{Princeton
    Mathematical Series}.  \newblock Princeton University Press,
  Princeton, N.J., 1980.

\bibitem[Mil03]{math.AG/0301304} James~S. Milne.  \newblock
  {Gerbes and abelian motives}, 2003, math.AG/0301304

\bibitem[SR72]{MR0338002} Neantro Saavedra~Rivano.  \newblock
  \emph{Cat\'egories {T}annakiennes}.  \newblock Springer-Verlag,
  Berlin, 1972.  \newblock Lecture Notes in Mathematics,
  Vol. 265.

\end{thebibliography}
\end{document}